\documentclass[a4paper, 11pt]{amsart}

\newcommand{\usecustompackage}[2][]{%
    \IfFileExists{Packages/#2.sty}{%
        \usepackage[#1]{./Packages/#2}%
    }{%
        \usepackage[#1]{#2}%
    }%
}
\usepackage{silence}
\usepackage[T1]{fontenc}
\usepackage[utf8]{inputenc}
\usepackage[english]{babel}
\usepackage[left=2.2cm, right=2.2cm, top=2.2cm, bottom=3.3cm]{geometry}

\usecustompackage{common_packages}

\usetikzlibrary{backgrounds, arrows.meta}

\usecustompackage{macros_generic}
\usecustompackage[cnt1=section,cnt2=subsection]{macros_theorems}
\usecustompackage{macros_bib}

\fancypagestyle{nonumber}{
    \fancyhf{}
}

\fancypagestyle{numbered}{
    \fancyhf{}
    \fancyhead[L]{}
    \fancyhead[C]{}
    \fancyhead[R]{}
    \fancyfoot[L]{}
    \fancyfoot[C]{\arabic{page}}
    \fancyfoot[R]{}
}

\newcolumntype{L}{>{$}l<{$}}
\newcolumntype{C}{>{$}c<{$}}
\newcolumntype{R}{>{$}r<{$}}

\pgfplotsset{compat=1.17}

\makeatletter
\newcommand*{\tocsamepage}{\@starttoc{toc}}
\makeatother

\title{Analytic lattice cohomology of isolated curve singularities}

\author{Tamás Ágoston}
\address{
    Alfréd Rényi Institute of Mathematics, Reáltanoda utca 13-15, H-1053, Budapest, Hungary \newline
    \hspace*{4mm} ELTE - Eötvös Loránd University, Dept. of Geometry, Budapest, Hungary
}
\email{agoston.tamas@renyi.hu}

\author{András Némethi}
\address{
    Alfréd Rényi Institute of Mathematics,
    Reáltanoda utca 13-15, H-1053, Budapest, Hungary \newline
    \hspace*{4mm} ELTE - Eötvös Loránd University, Dept. of Geometry, Budapest, Hungary \newline
    \hspace*{4mm}
    BCAM - Basque Center for Applied Math.,
    Mazarredo, 14 E48009 Bilbao, Basque Country – Spain
}
\email{nemethi.andras@renyi.hu}

\newcommand*\hh{\Got h}
\newcommand*\rank{\mathop{\Rm{rank}}\nolimits}

\newcommand*\sq{\mathord{\square}}
\let\oldchi\chi
\def\chi{\raisebox{.12em}{$\oldchi$}}

\let\longto\longrightarrow

\newcommand*\grverts{\Cal V}
\newcommand*\gredges{\Cal E}

\newcommand*\laurpols{\Cal T}
\newcommand*\fnring{\Cal O}
\newcommand*\semigp{\Cal S}

\newcommand*\cw{\Got X}
\newcommand*\subcw{\Got Y}

\newcommand*\lcoh{\bbh}
\newcommand*\lcohred[1][]{\lcoh_{red\if\relax#1\relax\else,#1\fi}}

\newcommand*\cell{C}
\newcommand*\cellface{D}
\newcommand*\cellset{\Cal Q}
\newcommand*\cellgp{\Scr C}
\newcommand*\celldualgp{\Cal F}

\newcommand*\maxpoint{\ell^*}

\newcommand*\contspace{\Cal O}
\newcommand*\spsystem{\Scr A}
\newcommand*\combcell{\ti C}
\newcommand*\combcellface{\ti D}
\newcommand*\topcell{\sq}

\newcommand*\cellroot{R}
\newcommand*\celltip{T}
\newcommand*\cellgens{\Cal G}
\newcommand*\cellgenss{\Cal H}
\newcommand*\cellgen{G}
\newcommand*\levelcells{\Scr S}
\newcommand*\levelset{S}

\newcommand*\filt{\Scr F}

\newcommand\Boxed[2][]{\tikz[overlay]{\node[anchor=text, draw=red, shape=rectangle, inner sep=3pt, #1] {$#2$};}\phantom{\TextOrMath{$#2$}{#2}}}
\newcommand\circled[2][]{\tikz[overlay]{\node[anchor=text, draw=red, shape=circle, inner sep=2pt, #1] {#2};}\phantom{#2}}

\newcommand*\refeqn[1]{\xref{#1}{(\num)}}

\begin{document}

\subjclass[2010]{Primary 32S05, 32S10, 32S25, 32S50; Secondary 14Bxx, 57K10}
\thanks{The authors are partially supported by NKFIH Grants ``Élvonal (Frontier)'' KKP 126683 and KKP 144148.}

\begin{abstract}
    We construct a lattice cohomology $\bbh^*(C,o)=\oplus_{q\geq 0}\bbh^q(C,o)$ and a graded root ${\Got R}(C,o)$ to any complex isolated curve singularity $(C,o)$. Each $\bbh^q(C,o)$ is a $\bbz$-graded $\bbz[U]$-module. The Euler characteristic of $\bbh^*(C,o)$ is the delta-invariant of $(C,o)$. The construction is based on the multivariable Hilbert series of the multifiltration provided by valuations of the normalization.

    Several examples are discussed, e.g. Gorenstein curves (where an additional symmetry is established), plane curves (in particular, Newton non-degenerate ones), ordinary $r$-tuples.

    We also prove that a flat deformation $(C_t,o)_{t\in (\bbc,0)}$ of isolated curve singularities induces an explicit degree zero graded $\bbz[U]$-module morphism $\lcoh^0(C_{t=0},o)\to \lcoh^0(C_{t\not=0},o)$, and a graded (graph) map of degree zero at the level of graded roots ${\Got R}(C_{t\not=0},o)\to{\Got R}(C_{t=0},o)$.

    In the treatment of the deformation functor we need a second construction of the lattice cohomology in terms of the system of linear subspace arrangement associated with the above filtration.
\end{abstract}

\maketitle

\section{Introduction}

\subsection{}
In this note we define a categorification of the delta invariant of an isolated curve singularity. This means that for any complex isolated curve singularity $(C,o)$ we define a cohomology module, the \It{lattice cohomology} $\bbh^*(C,o)$ of $(C,o)$, whose Euler characteristic ${\rm eu}(\bbh^*(C,o))$ is the delta invariant. The lattice cohomology $\bbh^*(C,o)$ can be written as a direct sum $\oplus_{q\geq 0}\bbh^q(C,o)$ of $\bbz[U]$-modules, where each summand $\bbh^q(C,o)$ is a $\bbz$-graded $\bbz[U]$-module, the $\bbz$-rank of any $d$-homogeneous component $\bbh^q(C,o)_d$ is finite, and the $U$-action is homogeneous of degree $-2$.

The construction of the lattice cohomology provides automatically an enhancement of $\bbh^0$, the \It{graded root}. In particular, in the case of curve singularities, we will also define a graded root ${\Got R}(C,o)$ associated with $(C,o)$.

In order to define a lattice cohomology one usually needs a lattice $\bbz^r$ and a weight function $w:\bbz^r\to \bbz$ (see section \ref{ss:latweight}). In different geometric situations the lattice appears rather naturally, however the weight function can be hidden in the geometrical structure.

Lattice cohomologies were already defined in different other situations. The first version was associated with the topological type of a complex normal surface singularity $(X,o)$ (whenever the link is a rational homology sphere). The link of such a singularity is a plumbed (graph) 3-manifold, where the plumbing graph can be chosen as one of the resolution graphs of the germ. In particular, it is connected and negative definite. In such a situation, using the combinatorics of the graph, a lattice cohomology $\bbh^*_{top}(X,o)$ was constructed, with an extra grading (the spin$^c$-structures of the link), whose Euler characteristic was the Seiberg-Witten invariant of the link \cite{Nlattice}. The lattice was the lattice generated by the vertices of the graph (i.e., integral divisors supported on the exceptional curve) and the weight function was a Riemann-Roch expression.

It was verified that the output $\bbh^*_{top}(X,o)$ is independent of the resolution \cite{Nlattice,Nkonyv}. In \cite{Nlattice} was conjectured (and in some cases verified) and in \cite{Zemke} proved that it coincides with the Heegaard Floer cohomology of the link. For some more connections with topology (e.g. applications of the graded roots) see \cite{DM,K2,KarakurtLidman,K1}.

The analytic version $\bbh^*_{an}(X,o)$ of the lattice cohomology associated with a normal surface singularity was defined in \cite{AgostonNemethiI,AgostonNemethiIII}, and even a graded $\bbz[U]$-module morphism $\bbh^*_{an}(X,o)\to \bbh^*_{top}(X,o)$ was provided. In this analytic case the Euler characteristic equals the geometric genus of the germ.

Later, in \cite{AgostonNemethiIV}, for any $n$-dimensional complex isolated singularity $(n\geq 2)$ an analytic lattice cohomology was constructed. In this higher dimensional case its topological pair is (still) missing (conjecturally it is the Embedded Contact Homology of the link).

In any dimension $n\geq 2$ the construction of $\bbh^*_{an}$ is based on the multivariable divisorial filtration provided by a resolution: the lattice is the free $\bbz$-module generated by the exceptional divisors, and the weight function was given by some sheaf-cohomology ranks.

\vspace{2mm}

In parallel with the lattice cohomology and graded roots in each situation (with given lattice and weight function) the \It{path lattice cohomology} is/can also be considered. Its major role in the case of study of analytic and topological types of normal surface singularities is that it provides a natural topological upper bound for a key analytic invariant, for the geometric genus (and for several families this bound is the sharp optimal one), cf. \cite{NSig,Nkonyv}. In some special links of surface singularities the Euler characteristic of the path lattice cohomologies associated with well--chosen pathes can coincide with the Euler characteristic of the topological lattice cohomology, however in general they behave rather differently (see e.g. \cite{Nkonyv}).

\subsection{}
In this note we analyse the case $n=1$. In this case the rank of the lattice is the number of irreducible components of $(C,o)$, and the weight function is provided by the Hilbert function associated with the normalization of the irreducible components. In order to prove that the Euler characteristic of the cohomology equals the delta invariant, we will develop a combinatorial setup for the lattice cohomology. In this general combinatorial situation we provide a formula for the Euler characteristic whenever the (combinatorial) Hilbert function satisfies a `matroid type inequality' and a `combinatorial duality property'. In our original analytic case of curve singularities the matroid inequality follows from the fact that the Hilbert function is given by valuations. Furthermore, the duality is also satisfied by analytic duality properties of curve singularities.

In the body of the paper we compute and analyse several examples. The Gorenstein and the plane singularity cases have distinguished roles. In the plane curve case the Hilbert function can also be determined from the embedded topological type (via multivariable Alexander polynomials, cf. \cite{GorNem2015}). Hence, a posteriori, for plane curves, the lattice cohomology $\bbh^*(C,o)$ is an invariant of the embedded topological type of the corresponding link in $S^3$. (In this case the new theory can be compared with the Link Heegaard Floer theory of the link, a program, which will be developed in the forthcoming notes.) However, when the embedding dimension of $(C,o)$ is larger than two, then similar topological interpretations do not exist, however our analytic setup runs independently of the magnitude of the embedding dimension; this emphasizes once more the power of our new analytic construction.

We exemplify certain higher embedding dimensional case as well, e.g. by the case of ordinary $r$-tuples, or by curves associated with numerical semigroups.

\vspace{2mm}

In the note we will consider the corresponding \It{path lattice cohomologies} associated with $(C,o)$ (and with a `lattice path') as well. As opposed to the topological case of surface singularities (cf. \cite{NSig,Nkonyv}), and in agreement with the analytic lattice cohomologies for $n\geq 2$ (cf. \cite{AgostonNemethiI,AgostonNemethiIII,AgostonNemethiIV}), in this case we will have the coincidence of the Euler characteristic of any path lattice cohomology and of the lattice cohomology, cf. Corollary \thmref{cor:EUcurves}.

\vspace{2mm}

If $(C,o)$ is a plane curve singularity then the above defined lattice cohomology has several deep connections with the Link Heegaard Floer theory of Ozsváth and Szabó (see e.g. \cite{OSzHol}) and also with local lattice cohomology considered in \cite{GorNem2015}. In order to realize them one has to introduce some additional filtration of $\bbh^*(C,o)$. All these will be the subject of a forthcoming manuscript.

\subsection{}
The second main goal is to analyse the behaviour of the new invariant under flat deformations.

Recall that in the case of Heegaard Floer theory (as indicated by any topological quantum field theory), to any 3-manifold one associates a graded $\bbz[U]$-module, and to any cobordism of 3-manifolds a graded $\bbz[U]$-module morphism of the corresponding modules. In our case, to any analytic type of curve germ we associate a cohomology module $\bbh^*(C,o)$ and to any analytic flat deformation of isolated curve singularities $\setof{(C_t,o)}_{t\in(\bbc,0)}$ a functor of the corresponding cohomology modules.

In Theorem \thmref{th:DEF2} we prove the following statement.

\begin{thm}[nonum, nostep]
    Let $\setof{(C_t,o)}_{t\in(\bbc,0)}$ be a flat deformation of isolated curve singularities. Then we construct an explicit degree zero graded $\bbz[U]$--module morphism $\lcoh^0(C_{t=0},o)\to \lcoh^0(C_{t\not=0},o)$, and similarly a graded (graph) map of degree zero at the level of graded roots ${\Got R}(C_{t\not=0},o)\to{\Got R}(C_{t=0},o)$.
\end{thm}

The fact that the constructed geometric map induces a morphism $\lcoh^{\geq 1}(C_{t=0},o)\to \lcoh^{\geq 1}(C_{t\not=0},o)$ is formulated as a conjecture.

In fact, if $(C,o)$ has $r$ irreducible components then $\bbh^{\geq r}(C,o)=0$. In particular, if along a flat deformation either $(C_{t=0},o)$ or $(C_{t\not=0},o)$ is irreducible, then the above theorem covers all the non-trivial cases.

In order to define the deformation functor, we needed to provide a second construction of the lattice cohomology. This is based on the structure of the linear subspaces of the local algebra provided by the ideal filtration of the valuations mentioned above. In this way we create the possibility to extend our present construction of the lattice cohomology to any situation where some linear subspace arrangements appear.

In subsequent articles we plan to apply the present results in the classification of adjacency relations of (plane) curve singularities.

\section{Preliminary. Basic properties of the lattice cohomology.}

\subsection{The lattice cohomology associated with a system of weights}\label{ss:latweight}\cite{Nlattice}\hfill\smallbreak

We consider a free $\bbz$-module, with a fixed basis $\setof{E_v}_{v\in\grverts}$, denoted by $\bbz^s$. It is also convenient to fix a total ordering of the index set $\grverts$, which will be denoted by $\setof{1,\ldots,s}$ henceforth.

The (original) construction of the \It{lattice cohomology} associates a graded $\bbz[U]$-module to the pair $(\bbz^s,\setof{E_v}_v)$ endowed with a set of weights. For the definition see \cite{Nlattice,NJEMS,Nkonyv}. First we set some notations regarding $\bbz[U]$-modules.

\begin{ntn}[label=9zu1]
    Consider the graded $\bbz[U]$-module $\laurpols:=\bbz[U,U^{-1}]$, and (following \cite{OSzP}) denote by $\laurpols_0^+$ its quotient by the submodule $U\cdot\bbz[U]$. This has a grading in such a way that $\deg(U^{-d})=2d$ ($d\geq 0$). Similarly, for any $n\geq 1$, the quotient of $U^{-(n-1)}\cdot\bbz[U]$ by $U\cdot\bbz[U]$ (with the same grading) defines the graded module $\laurpols_0(n)$. Hence, $\laurpols_0(n)$, as a $\bbz$-module, is freely generated by $1,U^{-1},\ldots,U^{-(n-1)}$, and has finite $\bbz$-rank $n$.

    More generally, for any graded $\bbz[U]$-module $P$ with $d$-homogeneous elements $P_d$, and for any $m\in\bbz$, we denote by $P[m]$ the same module graded in such a way that $P[m]_{d+m}=P_{d}$. Then set $\laurpols^+_m:=\laurpols^+_0[m]$ and $\laurpols_m(n):=\laurpols_0(n)[m]$. Hence, for $m\in\bbz$, $\laurpols_{2m}^+=\bbz\gen{U^{-m}, U^{-m-1},\ldots}$ as a $\bbz$-module.
\end{ntn}
There are two ways to define the lattice cohomology. The first one constructs a cochain complex from which the lattice cohomology arises. Classically, we consider the space $\bbz^s\otimes\bbr=\bbr^s$ (or, $(\bbr_{\geq 0})^s$) and its natural cellular decomposition into cubes determined by the lattice points and $\setof{E_v}_{v\in\grverts}$. This means that to each lattice point $\ell\in\bbz^s$ and $I\subseteq\grverts$ determines a cube $\topcell=(\ell, I)$ that has its vertices in the points $\ell+\sum\limits_{v\in I'}E_v$ ($I'\subseteq I$), cf. \cite{Nlattice}.

Here, we instead do this in slightly greater generality. Though the statements may be more general, the proofs are completely analogous, and we leave their verification to the reader. For the original proofs see \cite{NOSz,NGr,Nlattice,NJEMS,Nkonyv}.

\begin{ntn}[label={ntn:cw}]
    Let $\cw$ be a CW complex (for definitions and properties see e.g.~\cite{Whitehead}). Let $\setof{\sk_q\cw}_{q\geq 0}$ be the skeleton decomposition of $\cw$. In this note we prefer to replace the $q$-simpleces by $q$-cubes $[0,1]^q$. If $\setof*{\kappa_{q,\alpha}}_\alpha$ are the characteristic maps for the $q$-cells of $\cw$ then the homology classes $\topcell_{q,\alpha}\in H_q(\sk_q\cw,\sk_{q-1}\cw)$ of the cycles $\kappa_{q,\alpha}$ constitute the set $\cellset_q=\cellset_q(\cw)$ of `$q$-dimensional cells', and they form a basis in $\cellgp_q=\cellgp_q(\cw)=\bbz\gen{\cellset_q}$, the free $\bbz$-module generated by them.

    On each cell $\topcell$ of $\cw$ we fix an orientation. Then for any pair of cells $\topcell_q$ and $\topcell_{q-1}$ (of dimension $q$ and $q-1$ respectively) there exist integers $[\topcell_q,\topcell_{q-1}]$, almost all zero, the \It{incidence numbers}, and a homological boundary operator $\partial :\cellgp_q\to \cellgp_{q-1}$ such that
    \[
        \partial \topcell_q=\sum _{\topcell_{q-1}} [\topcell_q,\topcell_{q-1}]\,\topcell_{q-1}.
    \]
    Then $\partial^2=0$ and the complex $(\cellgp_*(\cw),\partial)$ provides the homology of $\cw$.

    In the sequel we will use the same notation $\topcell_q$ for the images of the characteristic maps $[0,1]^q\to \sk_q\cw$ (with restrictions $\partial [0,1]^q\to \sk_{q-1}\cw$); they are also called the `closed cells/cubes in $\cw$'. We also write $\topcell^\circ_q:= \topcell_q\setminus \sk_{q-1}\cw$ for their relative interiors. We call them the `open cells'.

    In this article we assume that $\cw$ is \It{regular} (in the sense of \cite[II.6]{Whitehead}). This means that each closed $q$-cell of $\cw$ is homeomorphic with $[0,1]^q$, and if $\topcell_q$ is a closed $q$-cell in $\cw$ then $\topcell\setminus \topcell^\circ$ is a union of finitely many $(q-1)$-cells of $\cw$.
\end{ntn}

\begin{rem}
    In the original setting of the cubical CW complex $\bbz^s\otimes\bbr$ \cite{Nlattice}, we can determine the orientations of the cubes/cells uniformly by fixing an ordering of the index set $\grverts$. Then the ordering of the indices in $I$ gives us a natural orientation of $\topcell=(\ell, I)$. For a general CW complex, we do not have such a way to globally choose the orientations of each cell, but apiece choices of the cell-orientations will not affect the
    cohomology we define shortly (see \thmref{cw-cohom-orient-invar}).
\end{rem}

Using the above setting, in order to define a `more interesting' cohomology theory, we consider a set of compatible \It{weight functions} $\setof{w_q}_q$.

\begin{defn}[label={def:compat-weight}]
    A set of functions $w_q\colon\cellset_q\to\bbz$ is called a \It{set of compatible weight functions} if
    \begin{enumerate}[label={(\alph*)}]
        \item $w_0$ is bounded from below;
        \item for any $\topcell_q\in\cellset_q$ and any point $p\in \topcell_q\setminus \topcell^\circ_q$, consider $r<q$ such that $p\in\sk_r\cw\setminus \sk_{r-1}\cw$, and $\topcell^\circ_r$, the unique open cell of $\sk_r\cw$ with $p\in \topcell_r^\circ$.
        Then in any such case we require $w_q(\topcell_q)\geq w_{r}(\topcell_{r})$.
    \end{enumerate}
    The index $q$ of $w_q$ may be omitted henceforth if it causes no confusion, i.e.~we set $w=\cup_{q} w_q$. Such a pair $(\cw,w)$ is called a \It{weighted CW complex}.
\end{defn}

\begin{rem}[label={rem:level-set-assumptions}]
    In the original setting of the cubical CW complex $\bbz^s$ in \cite{Nlattice}, we have required that $w_0^{-1}\bigl((-\infty,k]\bigr)$ is finite for all $k$, a stronger condition than $w_0$ being bounded from below. That will also be satisfied in all applications of this note (however, sometimes is useful to relax that condition).
\end{rem}

In the presence of any fixed set of compatible weight functions $w$ we define $\celldualgp^q=\celldualgp^q(\cw,w)$ as the set of morphisms $\Hom_{\bbz}(\cellgp_q,\laurpols^+_0)$ with finite support on $\cellset_q$.

Notice that $\celldualgp^q$ is a $\bbz[U]$-module by $(p*\oldphi)(\topcell_q):=p(\oldphi(\topcell_q))$ ($p\in\bbz[U]$, $\oldphi\in \celldualgp^q$). Moreover, $\celldualgp^q$ has a $\bbz$-grading: $\oldphi\in\celldualgp^q$ is homogeneous of degree $\deg(\oldphi)=d\in\bbz$ if for each $\topcell_q\in\cellset_q$ with $\oldphi(\topcell_q)\neq 0$, $\oldphi(\topcell_q)$ is a homogeneous element of $\laurpols^+_0$ of degree $d-2\cdot w(\topcell_q)$. (In fact, the grading is $2\bbz$-valued; hence, the reader interested only in the present construction may divide all the degrees by two. Nevertheless, we prefer to keep the present form in our presentation because of its resonance and certain comparisons with the Heegaard Floer theory. This convention was preserved in the case of the `topological' lattice cohomology associated with links of normal surface singularities, and the present theory will also be compared with all these cohomology theories.)

Next, we define $\delta_w\colon\celldualgp^q\to\celldualgp^{q+1}$. For this, fix $\oldphi\in\celldualgp^q$ and we show how $\delta_w\oldphi$ acts on a cell $\topcell_{q+1}\in\cellset_{q+1}$. First write $\partial \topcell_{q+1}=\sum _{\topcell_{q}} [\topcell_{q+1},\topcell_{q}]\,\topcell_{q}$, then set
\[
    (\delta_w\oldphi)(\topcell_{q+1}):=\sum_{\topcell_q} \ [\topcell_{q+1}, \topcell_{q}]\,
    U^{w(\topcell_{q+1})-w(\topcell_q)}\, \oldphi(\topcell_q).
\]
One can show that $\delta_w\circ\delta_w=0$, i.e. $(\celldualgp^*,\delta_w)$ is a cochain complex.

\medskip

In fact, $(\celldualgp^*,\delta_w)$ has a natural \It{augmentation} too. Indeed, set $m_w:=\min\limits_{x\in\cellset_0}w_0(x)$ and choose $x_w\in\cellset_{0}$ such that $w_0(x_w)=m_w$. Then one defines the $\bbz[U]$-linear map
\[
    \epsilon_w\colon\laurpols^+_{2m_w}\longto\celldualgp^0
\]
such that $\epsilon_w(U^{-m_w-s})(x)$ is the class of $U^{-m_w+w_0(x)-s}$ in $\laurpols^+_0$ for any $x\in \cellset_0$ and $s\geq 0$. Then $\epsilon_w$ is injective, $\delta_w\circ\epsilon_w=0$, and $\epsilon_w,\delta_w$ are homogeneous morphisms between $\bbz[U]$-modules of degree $0$.

\begin{defn}[label=9def12]
    The homology of the cochain complex $(\celldualgp^*,\delta_w)$ is called the \It{lattice cohomology} of the pair $(\cw,w)$, and it is denoted by $\lcoh^*(\cw,w)$. The homology of the augmented cochain complex
    \[
        0\longrightarrow
        \laurpols^+_{2m_w}\stackrel{\epsilon_w}{\longrightarrow}
        \celldualgp^0\stackrel{\delta_w}{\longrightarrow}
        \celldualgp^1\stackrel{\delta_w}{\longrightarrow}
        \ldots
    \]
    is called the \It{reduced lattice cohomology} of the pair $(\cw,w)$, and it is denoted by $\lcohred^*(\cw,w)$.
\end{defn}

If the pair $(\cw,w)$ is clear from the context, we omit it from the notation.

For any fixed $q\geq 0$, the $\bbz$-grading of $\celldualgp^q$ induces a $\bbz$-grading on $\lcoh^q$ and $\lcohred^q$; the homogeneous part of degree $d$ is denoted by $\lcoh^q_d$ and $\lcohred[d]^q$ respectively. Moreover, both $\lcoh^q$ and $\lcohred^q$ admit an induced graded $\bbz[U]$-module structure and $\lcoh^q=\lcohred^q$ for $q>0$.

\begin{lem}[label=9lemma3]
    One has a graded $\bbz[U]$-module isomorphism $\lcoh^0=\laurpols^+_{2m_w}\oplus\lcohred^0$.
\end{lem}

Next, we present another realization of the modules $\lcoh^*$.

\begin{defn}[label=9rem]
    For each $n\in\bbz$ define $S_n=S_n(w)\subseteq\cw$ as the union of all the closed cells $\topcell_q$ (of any dimension) with $w(\topcell_q)\leq n$. Clearly, $S_n=\emptyset$, whenever $n<m_w$. For any $q\geq 0$, set
    \[
        \bbs^q(\cw,w):=\bigoplus_{n\geq m_w}\, H^q(S_n,\bbz).
    \]
    Then $\bbs^q$ is $\bbz$ (in fact, $2\bbz$)-graded: the $d=2n$-homogeneous elements $\bbs^q_d$ consist of $H^q(S_n,\bbz)$.

    Furthermore, $\bbs^q$ has a $\bbz[U]$-module structure: the $U$-action is given by the restriction map $r_{n+1}\colon H^q(S_{n+1},\bbz)\to H^q(S_n,\bbz)$. Namely, $U*(\alpha_n)_n=(r_{n+1}\alpha_{n+1})_n$. Moreover, for $q=0$, the fixed base point $x_w\in S_n$ provides an augmentation (splitting) $H^0(S_n,\bbz)=\bbz\oplus\ti{H}^0(S_n,\bbz)$, hence an augmentation of the graded $\bbz[U]$-modules
    \[
        \bbs^0=
        \laurpols^+_{2m_w}\oplus\bbs^0_{red}=
        \Bigl(\bigoplus_{n\geq m_w}\bbz\Bigr)\oplus\Bigl(\bigoplus_{n\geq m_w}\ti{H}^0(S_n,\bbz)\Bigr).
    \]
\end{defn}

\begin{thm}[label=9STR1]
    There exists a graded $\bbz[U]$-module isomorphism, compatible with the augmentations: $\lcoh^*(\cw,w)=\bbs^*(\cw,w)$.
\end{thm}
\begin{cor}[label=cw-cohom-orient-invar]
    The cohomology $\lcoh^*(\cw,w)$ is invariant on the choice of orientation for the cells.
\end{cor}

We also wish to define the Euler characteristic of the cohomology $\lcoh^*(\cw,w)$. First note that even should $\lcohred^*(\cw,w)$ have finite $\bbz$-rank in any fixed homogeneous degree (which is guaranteed in the classical case, but not in our more permissive definition, see \thmref{rem:level-set-assumptions}), it is not, in general, finitely generated over $\bbz$ --- in fact, not even over $\bbz[U]$.

\begin{ex}[label={ex:notFinGen}]
    Set $\cw=\bbr$ with the standard decomposition given by $\bbz$, and define $w_0$ as
    \[
        w_0(-n)=w_0(n)=[n/2]+4\{n/2\}\text{ for any }n\in\bbz_{\geq 0},
    \]
    where $[\ ]$ and $\{\ \}$ denote the integral and fractional parts respectively; and let $w_1$ on the segment $[n,n+1]$ take the value $\max\setof{w_0(n),w_0(n+1)}$. Then $\lcohred^0=\bigoplus\limits_{k\geq 1}\laurpols_k(1)^2$.
\end{ex}

However, in cases that interest us, $\lcohred^*$ will indeed be finitely generated over $\bbz$.

\begin{defn}[label={def:eu}, description={Euler characteristic}]
    Assume that $\lcohred^*(\cw,w)$ has finite $\bbz$-rank (e.g. when $\cw$ is a finite CW complex). We define the Euler characteristic of $\lcoh^*(\cw,w)$ as
    \[
        \eu\bigl(\lcoh^*(\cw,w)\bigr)=-\min\setdef{w(x)}{x\in\cellset_0}+\sum_q(-1)^q\rank_\bbz\bigl(\lcohred^q(\cw,w)\bigr).
    \]
\end{defn}

\begin{lem}[label={lem:eu-w}]
    If $\cw$ is a finite CW complex with (classical) topological Euler characteristic one (e.g. it is contractible) then
    \begin{equation}\label{eq:eu}
        \eu\bigl(\lcoh^*(\cw,w)\bigr)=\sum_{q\geq 0}\ \sum_{\topcell_q\in\cellset_q} (-1)^{q+1}w(\topcell_q).
    \end{equation}
\end{lem}

\begin{ex}[label={lat-coh-restriction}, description={Restrictions}]
    Let us consider a subcomplex $\subcw\subseteq\cw$ and the restriction of $w$ on $\subcw$.
    Then the restriction map $r\colon\celldualgp^q(\cw,w)\to\celldualgp^q(\subcw,w)$ induces a natural homogeneous graded $\bbz[U]$-module homomorphism of degree zero
    \[
        r^*\colon\lcoh^*(\cw,w)\to\lcoh^*(\subcw,w|_\subcw).
    \]

    E.g., in the case of $\cw=\bbr^s$ (with cell-decomposition induced by $\bbz^s$ and $\setof*{E_v}_v$), $\subcw$ can be one of the following:
    \begin{enumerate}[label={(\arabic*)}]
        \item the first quadrant $(\bbr_{\geq 0})^s$ with $\cellset_0=(\bbz_{\geq 0})^s$;
        \item the rectangle $R=R(0,c)=\setdef{x\in\bbr^s}{0\leq x\leq c}$ for some lattice point $c
        \in(\bbz_{\geq 0})^s$;
        \item a connected path composed of edges in the lattice, see \thmref{def:pathlat}.
    \end{enumerate}
\end{ex}

In the above cases (2) and (3) $\cw$ is finite and contractible, hence lemma \thmref{lem:eu-w} applies.

If $\cw=\sk_s\cw$ for some $s$, then $\lcoh^q(\cw,w)=0$ for $q>s$. If $\cw$ is any of the cases (1) or (2) above then $\lcoh^q(\cw,w)=0$ for $q\geq s$. (For the vanishing in case (3) see below.)

\begin{defn}[label={def:pathlat}, description={Path lattice cohomology}]
    Fix a weighted cubical CW complex $(\bbz^s\otimes\bbr,w)$ (with cell-decomposition induced by $\bbz^s$ and $\setof*{E_v}_v$). Consider a sequence $\gamma=\setof{x_i}_{i=0}^t$ so that $x_0=0$, and $x_{i+1}=x_i+ E_{v(i)}$ for $0\leq i<t$. We write $T$ for the union of 0-cubes marked by the points $\setof{x_i}_{i=0}^t$ and of the segments of type $[x_i,x_{i+1}]$. Then $\lcoh^*(T,w)$ is called the \It{path lattice cohomology} associated with the `path' $\gamma$ and weights $\setof{w_q}_{q=0,1}$, which we will also denote by $\lcoh^*(\gamma,w)$. Also, by \thmref{lat-coh-restriction} we get a $0$-homogeneous $\bbz[U]$-module homomorphism $\lcoh^*(\bbr^s,w)\to\lcoh^*(\gamma,w)$.

    From the realization in \thmref{9rem}, we immediately have $\lcoh^q(\gamma,w)=0$ for $q\geq 1$. Furthermore, we have an augmentation with sum decomposition
    \[
        \lcoh^0(\gamma,w)=\laurpols_{2m_\gamma}^+\oplus\lcohred^0(\gamma,w),
    \]
    where $m_\gamma=\min\limits_i\setof{w_0(x_i)}$, and $\lcohred^0(\gamma,w)$ is the \It{reduced path lattice cohomology}.
\end{defn}

Since the above $T$ is a finite CW complex, by \thmref{def:eu}, we can also consider its Euler characteristic:
\[
    \eu\bigl(\lcoh^*(\gamma,w)\bigr)=
    -m_\gamma+\rank_\bbz\,\bigl(\lcohred^0(\gamma,w)\bigr).
\]
Sometimes we also write $\eu\bigl(\lcoh^0(\gamma,w)\bigr)$ instead of $\eu\bigl(\lcoh^*(\gamma,w)\bigr)$. It can also be obtained as follows.

\begin{lem}[label=9PC2]
    We have
    \[
        \eu\bigl(\lcoh^*(\gamma,w)\bigr)=-w_0(0)+\sum_{i=0}^{t-1}w_1([x_{i},x_{i+1}])-w_0(x_{i+1}).
    \]
    If $w_1([x_{i},x_{i+1}])=\max\setof{w_0(x_{i}),w_0(x_{i+1})}$ then
    \begin{equation}\label{eq:pathweights}
        \eu\bigl(\lcoh^*(\gamma,w)\bigr)=-w_0(0)+\sum_{i=0}^{t-1}\max\setof{0, w_0(x_i)-w_0(x_{i+1})}.
    \end{equation}
\end{lem}

\subsection{Graded roots and their cohomologies}\label{s:grgen}\cite{NOSz,NGr}\hfill\smallbreak

\begin{defn}[label=def:gr-root]
    Let $\Got R$ be an infinite tree with vertices $\grverts$ and edges $\gredges$. We denote by $[u,v]$ the edge with end vertices $u$ and $v$. We say that $\Got R$ is a \It{graded root} with grading $\chi\colon\grverts\to\bbz$ if
    \begin{enumerate}[label={(\alph*)}]
        \item $\chi(u)-\chi(v)=\pm 1$ for any $[u,v]\in\gredges$;
        \item $\chi(u)>\min\setof{\chi(v),\chi(w)}$ for any $[u,v],[u,w]\in\gredges$, $v\neq w$;
        \item $\chi$ is bounded from below, $\chi^{-1}(k)$ is finite for any $k\in\bbz$, and $\abs*{\chi^{-1}(k)}=1$ if $k\gg 0$.
    \end{enumerate}
\end{defn}
A path connecting two vertices is \It{monotone} if $\chi$ restricted to the set of vertices on the path is strictly monotonous. If a vertex $v$ can be connected by another vertex $w$ by a monotone path and $\chi(v)>\chi(w)$, then we write $v\succ w$. This gives us a partial ordering on $\grverts$.

\begin{defn}[description={The cohomology $\bbz[U]$-module associated with a graded root}]
    For any graded root $(\Got R,\chi)$, let $\lcoh(\Got R,\chi)$ be the set of functions $\oldphi\colon\grverts\to\laurpols^+_0$ with the following property: whenever $[v,w]\in\gredges$ with $\chi(v)<\chi(w)$, then $U\cdot \oldphi(v)=\oldphi(w)$. Or, equivalently, one requires
    \begin{equation}\label{eq:bH}
        U^{\chi(w)-\chi(v)}\cdot\oldphi(v)=\oldphi(w)\text{ for any }w\succ v.
    \end{equation}
    Clearly $\lcoh(\Got R,\chi)$ is a $\bbz[U]$-module via $(U\oldphi)(v)=U\cdot\oldphi(v)$. Moreover, $\lcoh(\Got R,\chi)$ has a $\bbz$-grading: the element $\oldphi\in\lcoh(\Got R,\chi)$ is homogeneous of degree $d\in\bbz$ if for each $v\in\grverts$ with $\oldphi(v)\neq 0$, $\oldphi(v)\in\laurpols^+_0$ is homogeneous of degree $d-2\chi(v)$. Since $2\chi(v)+\deg\oldphi(v)=2\chi(w)+\deg\oldphi(w)$ in \refeqn{eq:bH}, $d$ is well-defined. Note also that any $\oldphi$ as above is automatically finitely supported.
\end{defn}
For a detailed concrete description of $\lcoh(\Got R)$ in terms of the combinatorics of the root $\Got R$ see \cite{NOSz}.

\begin{ex}
    \begin{enumerate}[label={(\alph*)}]
        \item For any integer $n\in\bbz$, let $\Got R_{(n)}$ be the tree with $\grverts=\setof{v^{k}}_{k\geq n}$ and $\gredges=\setof{[v^{k},v^{k+1}]}_{k\geq n}$, and the grading set to $\chi(v^{k})=k$. Then $\lcoh(\Got R_{(n)})=\laurpols_{2n}^+$.
        
        \item Consider the following two graded roots:
        
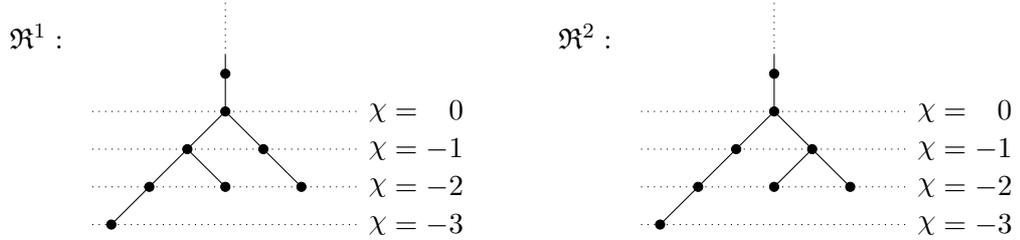
\begin{center}
            \begin{tikzpicture}[scale=0.5]
                \node[left] at (-4,2) {$\Got R^1:$};

                \draw[fill] (0,1) circle [radius=0.12];
                \draw[fill] (0,0) circle [radius=0.12];
                \draw[fill] (-1,-1) circle [radius=0.12];
                \draw[fill] (1,-1) circle [radius=0.12];
                \draw[fill] (-2,-2) circle [radius=0.12];
                \draw[fill] (0,-2) circle [radius=0.12];
                \draw[fill] (2,-2) circle [radius=0.12];
                \draw[fill] (-3,-3) circle [radius=0.12];

                \draw (0,0) -- (0,1.5);
                \draw[dotted] (0,1.5) -- (0,3);
                \draw (0,0) -- (-3,-3);
                \draw (0,0) -- (2,-2);
                \draw (-1,-1) -- (0,-2);

                \draw[dotted] (-3.5,0) -- (3.5,0);
                \draw[dotted] (-3.5,-1) -- (3.5,-1);
                \draw[dotted] (-3.5,-2) -- (3.5,-2);
                \draw[dotted] (-3.5,-3) -- (3.5,-3);

                \node[right] at (3.5,0) {$\chi=\phantom{-}0$};
                \node[right] at (3.5,-1) {$\chi=-1$};
                \node[right] at (3.5,-2) {$\chi=-2$};
                \node[right] at (3.5,-3) {$\chi=-3$};

                \node at (8,2) { };
            \end{tikzpicture}
            \begin{tikzpicture}[scale=0.5]
                \node[left] at (-4,2) {$\Got R^2:$};

                \draw[fill] (0,1) circle [radius=0.12];
                \draw[fill] (0,0) circle [radius=0.12];
                \draw[fill] (-1,-1) circle [radius=0.12];
                \draw[fill] (1,-1) circle [radius=0.12];
                \draw[fill] (-2,-2) circle [radius=0.12];
                \draw[fill] (0,-2) circle [radius=0.12];
                \draw[fill] (2,-2) circle [radius=0.12];
                \draw[fill] (-3,-3) circle [radius=0.12];

                \draw (0,0) -- (0,1.5);
                \draw[dotted] (0,1.5) -- (0,3);
                \draw (0,0) -- (-3,-3);
                \draw (0,0) -- (2,-2);
                \draw (1,-1) -- (0,-2);

                \draw[dotted] (-3.5,0) -- (3.5,0);
                \draw[dotted] (-3.5,-1) -- (3.5,-1);
                \draw[dotted] (-3.5,-2) -- (3.5,-2);
                \draw[dotted] (-3.5,-3) -- (3.5,-3);

                \node[right] at (3.5,0) {$\chi=\phantom{-}0$};
                \node[right] at (3.5,-1) {$\chi=-1$};
                \node[right] at (3.5,-2) {$\chi=-2$};
                \node[right] at (3.5,-3) {$\chi=-3$};
            \end{tikzpicture}
        \end{center}

        These graded roots $\Got R^1$ and $\Got R^2$ are not isomorphic but their $\bbz[U]$-modules are isomorphic: $\lcoh(\Got R^1)=\lcoh(\Got R^2)=\laurpols^+_{-6}\oplus \laurpols_{-4}(1)\oplus\laurpols_{-4}(2)$. Hence, in general, a graded root carries more information than its $\bbz[U]$-module.
    \end{enumerate}
\end{ex}

\begin{defn}[label={def:GRootW}, description={The graded root associated with a
weighted CW complex}]
    Fix a CW complex $\cw$ and a system of weights $\setof{w_q}_q$ as above. Consider the sequence of spaces $\setof{S_n}_{n\geq m_w}$ with $S_n\subset S_{n+1}$ as in \thmref{9rem}. Let $\pi_0(S_n)=\setof{\Cal C_n^1,\ldots,\Cal C_n^{p_n}}$ be the set of connected components of $S_n$.

    Then we define the graded graph $(\Got R_w,\chi_w)$ as follows. The vertex set $\grverts(\Got R_w)$ is $\cup_{n\geq m_w} \pi_0(S_n)$. The grading $\chi_w\colon\grverts(\Got R_w)\to\bbz$ is $\chi_w(\Cal C_n^j)=n$, that is, $\chi_w|_{\pi_0(S_n)}=n$.

    Furthermore, if $\Cal C_{n}^i\subset\Cal C_{n+1}^j$ for some $n$, $i$ and $j$, then we introduce an edge $[\Cal C_n^i,\Cal C_{n+1}^j]$. All the edges of $\Got R_w$ are obtained in this way.
\end{defn}

\begin{lem}[label=lem:GRoot]
    $(\Got R_w,\chi_w)$ satisfies all the required properties of the definition of a graded root, except possibly the last one: not all sets $\chi_w^{-1}(k)$ are finite, nor is $\abs*{\chi_w^{-1}(n)}=1$ guaranteed for $n\gg 0$ in case of an arbitrary CW complex $\cw$. It is, however, true when $\cw$ is a finite, contractible CW complex.
\end{lem}

For a case when the property $\abs*{\chi_w^{-1}(n)}=1$ for $n\gg 0$ is not satisfied see \thmref{ex:notFinGen}.

\begin{prop}[label={prop:HH}]
    If $\Got R$ is the graded root associated with $(\cw,w)$, all sets $\chi_w^{-1}(k)$ are finite, and $\abs*{\chi_w^{-1}(n)}=1$ for $n\gg 0$, then $\lcoh(\Got R)=\lcoh^0(\cw,w)$.
\end{prop}

\section{Combinatorial lattice cohomology} \label{sec:CombLattice}

\subsection{Combinatorial weight functions}\label{ss:CWF}\hfill\smallbreak

In the previous section, we introduced the cohomology $\lcoh^*$ for arbitrary weighted CW complexes. In this section, however, we will work with the classical case of a rectangle $R=R(0,c)=\setdef{x\in \bbr^s}{0\leq x\leq c}$ in $\bbr^s$, where the cubical structure is induced by $L=\bbz\gendef{E_v}{v\in\grverts}=\bbz^s\subset \bbr^s$. We allow $c=\infty$ as well, in that case $R=(\bbr_{\geq 0})^s$. Let $E_I=\sum\limits_{v\in I}E_v$ for $I\subset\grverts$ and $E=E_\grverts$.

\begin{defn}[label={def:comblattice}]
    For some $c\in\bbz^s$, $c\geq E$, consider the subcomplex $R=R(0,c)$ (see \thmref{lat-coh-restriction}), and assume that to each $\ell\in R\cap\bbz^s$ we assign
    \begin{enumerate}[label={(\roman*)}]
        \item an integer $h(\ell)$ such that $h(0)=0$ and $h(\ell+E_v)\geq h(\ell)$ for any $v$,
        \item an integer $h^\circ (\ell)$ such that $h^\circ(\ell+E_v)\leq h^\circ(\ell)$ for any $v$.
    \end{enumerate}
    \nni Once $h$ is fixed with (i), and $c\in (\bbz_{\geq 0})^s$ is finite, a possible choice for $h^\circ $ is $h^{sym}$, where $h^{sym}(\ell)=h(c-\ell)$. Clearly, the definition of $h^{sym}$ depends on $c$.

    We consider the set of cubes $\setof{\cellset_q}_{q\geq 0}$ of $R$ as in \thmref{ntn:cw} and the weight function
    \[
        w_0\colon\cellset_0\to\bbz\text{ by }w_0(\ell):=h(\ell)+h^\circ (\ell)-h^\circ(0).
    \]
    Clearly $w_0(0)=0$. Furthermore, we define
    \[
        w_q\colon\cellset_q\to\bbz,\quad
        w_q(\topcell_q)=\max\setdef{w_0(\ell)}{\ell\text{ is a vertex of }\topcell_q}.
    \]
    We will use the symbol $w$ for the system $\setof{w_q}_q$. The compatible weight functions define the lattice cohomology $\lcoh^*(R,w)$. Moreover, whenever $c$ is finite, for any increasing path $\gamma$ connecting 0 and $c$, we also have a path lattice cohomology $\lcoh^0(\gamma,w)$ as in \thmref{def:pathlat}. Accordingly, we have the numerical Euler characteristics $\eu\bigl(\lcoh^*(R,w)\bigr)$, $\eu\bigl(\lcoh^0(\gamma,w)\bigr)$ and $\min\limits_\gamma\eu\bigl(\lcoh^0(\gamma,w)\bigr)$ too.
\end{defn}

\begin{lem}[label={lem:comblat}] Assume that $c$ is finite. Then
    we have $0\leq\eu\bigl(\lcoh^0(\gamma,w)\bigr)\leq h^\circ (0)-h^\circ (c)$ for any increasing path $\gamma$ connecting 0 to $c$. The equality $\eu\bigl(\lcoh^0(\gamma,w)\bigr)=h^\circ (0)-h^\circ (c)$ holds if and only if for any $i$ the differences $h(x_{i+1})-h(x_i)$ and $h^\circ (x_{i})-h^\circ (x_{i+1})$ are not simultaneously positive.
\end{lem}
\begin{pf}
    Since $w_0(0)=0$ we have
    \[
        \eu\bigl(\lcoh^0(\gamma,w)\bigr)=-\min\, w_0+\rank\lcohred^0(\gamma,w)\geq 0.
    \]
    Next, by \refeqn{eq:pathweights},
    \[
        \eu\bigl(\lcoh^0(\gamma,w)\bigr)=\sum_{i=0}^{t-1}\max\setof{0,w(x_i)-w(x_{i+1})}
    \]
    On the other hand, via $\max\setof{0,x}+\min\setof{0,x}=x$,
    \[
        w(c)+\sum_i\max\setof{0,w(x_i)-w(x_{i+1})}=-\sum_i\min \setof{0,w(x_i)-w(x_{i+1})},
    \]
    hence, via $\max\setof{0,x}-\min\setof{0,x}=\abs{x}$,
    \begin{align*}
        w(c)+ 2\cdot\eu\bigl(\lcoh^0(\gamma,w)\bigr)
            & =\sum_i\abs{w(x_i)-w(x_{i+1}}
            =\sum_i\abs{h(x_{i})+h^\circ (x_i)-h(x_{i+1})-h^\circ (x_{i+1})}\leq \\
        & \leq\sum_i ( h(x_{i+1})-h(x_i))+\sum_i (h^\circ (x_i)-h^\circ (x_{i+1}))
            =h(c)+h^\circ (0)-h^\circ (c). \qedhere
    \end{align*}
\end{pf}

\begin{rem}[label={rem:ketpelda}]
    The inequality $0\leq\eu\bigl(\lcoh^*(R,w)\bigr)$ is not true in general. Take e.g. the following table ($s=2$, $c=(2,2)$) with $h^\circ =h^{sym}$, and where $(0,0)$ is the left-bottom corner,
    \begin{center}
        \setlength{\tabcolsep}{0.4em}
        \begin{tabular}{RRRRp{1cm}RRRRp{1cm}L}
                & 0 & 1 & 1 &   &       & 0 & 0 & 0 &   & \min w_0=0 \\
            h:  & 0 & 1 & 1 &   & w_0:  & 0 & 1 & 0 &   & \lcohred^0=0,\, \lcoh^1=\bbz \\
                & 0 & 0 & 1 &   &       & 0 & 0 & 0 &   & \eu=-1
        \end{tabular}
    \end{center}
    \nni
    The inequality $\eu\bigl(\lcoh^*(R,w)\bigr)\leq h^\circ (0)-h^\circ(c)$ might fail too, even with the choice $h^\circ=h^{sym}$:
    \begin{center}
        \setlength{\tabcolsep}{0.6em}
        \def\w#1 {\makebox[0pt][r]{$#1$}}
        \begin{tabular}{RRRRp{1cm}RRRRp{1cm}L}
                & 0 & 1 & 2 &   &       &\w-2   &\w-1   & 0     &   & \min w_0=-2 \\
            h:  & 0 & 1 & 1 &   & w_0:  &\w-1   & 0 &\w-1   &   & \rank\lcohred^0=2,\, \lcoh^1=0 \\
                & 0 & 0 & 0 &   &       & 0     &\w-1   &\w-2   &   & \eu=4
        \end{tabular}
    \end{center}
\end{rem}

In what follows, we are going to nail down what conditions we should impose on the functions $h$ and $h^\circ$ (based on certain specific examples), in which cases the cohomology does have nice properties --- it turns out that surprisingly little is required to formulate some interesting statements.

Basically we will use two types of restrictions, one of them is a `matroid type inequalities' (cf. \thmref{ex:latgrfilt}), the other is some kind of `duality property' (in applications it is the trace of genuine duality statements in algebraic geometry, for motivation for the next definition see Lemma \thmref{lem:comblat}, duality properties of Gorenstein curves in \ref{ss:GORduality} and Newton non-degenerate curves in \thmref{ex:latgrfilt2}).

\begin{defn}[label={def:COMPGOR}]
    Fix $(h,h^\circ,R)$ as in \thmref{def:comblattice}. We say that the pair $h$ and $h^\circ$ satisfy the `Combinatorial Duality Property' (CDP) if $h(\ell+E_v)-h(\ell)$ and $h^\circ(\ell+E_v)-h^\circ(\ell)$ simultaneously cannot be nonzero for $\ell,\ell+E_v\in R$. Furthermore, we say that $h$ satisfies the CDP if the pair $(h,h^{sym})$ satisfies it.
\end{defn}

\begin{ex}
    The `Combinatorial Duality Property' for the analytic lattice cohomology, case $n\geq 2$, was verified in \cite{AgostonNemethiI,AgostonNemethiIII,AgostonNemethiIV}. The `Combinatorial Duality Property' in the context of the topological lattice cohomology of normal surface singularities usually is not true. For a combinatorial setup when CDP holds see Example \thmref{ex:latgrfilt2}.
\end{ex}

\begin{ex}[label={ex:latgrfilt}, description={Graded and filtered vector spaces}]
    In several geometrical constructions we face the following situation: we have $\bbz^s$ and $c$ as above, and a finite dimensional vector space $M$ with a $\bbz^s$-grading such that $M_{\Bf{a}}=0$ whenever either $\Bf{a}\ngeq 0$ or $\Bf{a}\geq c$ (where $M_{\Bf a}$ denotes the space of homogeneous elements of degree $\Bf a$). Let $\hh$ be the Hilbert function $\hh(\ell)=\sum_{\Bf{a}\ngeq\ell}\,\dim M_{\Bf{a}}$ and let $h$ be its restriction to $R=R(0,c)$. Then $h(0)=0$ and $h(c)=\dim M$.

    More generally, assume that $M$ is a finite dimensional vector space endowed with a decreasing $\bbz^s$-filtration such that $\filt(0)=M$ and $\filt(c)=0$ and define $\hh(\ell)=\dim\bigl(M/\filt(\ell)\bigr)$ for any $\ell\geq 0$. Again, define $h$ as the restriction of $\hh$ to $R$.

    If the filtration satisfies $\filt(\max\setof{\ell_1,\ell_2})=\filt(\ell_1)\cap \filt(\ell_2)$ (e.g. it is induced by valuation or order functions) then the $h$-function associated with a filtration satisfies the following \It{matroid inequality}
    \begin{equation}\label{eq:matroid}
        h(\ell_1)+h(\ell_2)\geq h(\min\setof{\ell_1,\ell_2})+h(\max\setof{\ell_1,\ell_2}), \ \ \ell_1,\ell_2\in R.
    \end{equation}
    This implies the \It{`stability property'}, valid for any $\bar{\ell}\geq 0$ with $E_v$-multiplicity of $\bar{\ell}$ zero,
    \begin{equation}\label{eq:stability}
        h(\ell)=h(\ell+E_v)\ \ \Rightarrow\ \ h(\ell+\bar{\ell})=h(\ell+\bar{\ell}+E_v).
    \end{equation}
    Such properties of a $h$-function are very natural, and they will be used intensively in the sequel.

    Let us perform an Euler characteristic computation as well. Assume that $s=1$. In this case, from the point of view of $h$, we can assume that we are in the graded vector space case. Indeed, in this case $h$ is associated with $M=\bbc^{h(c)}$ graded as $M_\ell=\bbc^{h(\ell+1)-h(\ell)}$ for $\ell<c$ and $M_c=0$. Moreover, if $h^\circ=h^{sym}$, then
    \[
        w_0(\ell)=\sum_{s<\ell}\dim M_s+\sum_{s<c-\ell}\dim M_s-h(c).
    \]
    If $h$ satisfies the CDP then
    \[
        w_0(\ell+1)-w_0(\ell)=
        \begin{cases}
            \dim M_\ell & \text{if }h(\ell+1)>h(\ell), \\
            -\dim M_{c-1-\ell} & \text{if }h(c-\ell)>h(c-\ell-1), \\
            0 & \text{otherwise.}
        \end{cases}
    \]
    Therefore (compare with Lemma \thmref{lem:comblat})
    \begin{equation}\label{eq:eurankone}
        \eu\bigl(\lcoh^*(R,w)\bigr)=\sum_{\ell=0}^{c-1}\bigl(w_1([\ell,\ell+1])-w_0(\ell)\bigr)=\sum_{\ell=0}^{c-1}\dim M_\ell=\dim M=h^\circ (0)-h^\circ (c).
    \end{equation}
\end{ex}

\begin{defn}[label={def:comblat}]
    We say that the pair $(h, h^\circ) $ satisfies the
    \begin{enumerate}[label={(\alph*)}]
        \item \It{`path eu-coincidence'} if $\eu(\lcoh^0(\gamma,w))=h^\circ (0)-h^\circ (c)$ for any increasing path $\gamma$.
        \item \It{`eu-coincidence'} if $\eu(\lcoh^*(R,w))=h^\circ (0)-h^\circ (c)$.
    \end{enumerate}
\end{defn}

\begin{thm}[label={th:comblattice}]
    Assume that $h$ satisfies the stability property, and the pair $(h,h^\circ)$ satisfies the Combinatorial Duality Property. Then $(h,h^\circ)$ satisfies both the path eu-coincidence and the eu-coincidence properties: for any increasing $\gamma$,
    \[
        \eu(\lcoh^*(\gamma,w))=\eu(\lcoh^*(R,w))=h^\circ (0)-h^\circ(c).
    \]
\end{thm}
\begin{pf}
    The identity $\eu(\lcoh^*(\gamma,w))=h^\circ (0)-h^\circ (c)$ follows from Lemma \thmref{lem:comblat}. Next we focus on the second identity. We claim that for any $I\subset\grverts$ we have
    \begin{equation}\label{eq:KEYIDEN}
        w\bigl((\ell,I)\bigr)-w(\ell)=h(\ell+E_I)-h(\ell).
    \end{equation}
    We use induction over the cardinality $\abs{I}$ of $I$. If $I=\setof{v}$, then
    \[
        w\bigl((\ell,I)\bigr)-w(\ell)=\max\setof{0, w(\ell+E_v)-w(\ell)}
    \]
    But if $h(\ell+E_v)>h(\ell)$ then $h^\circ (\ell+E_v)=h^\circ (\ell)$ hence $w(\ell+E_v)-w(\ell)=h(\ell+E_v)-h(\ell)$. Otherwise $w(\ell+E_v)\leq w(\ell)$.

    Next, assume that $\abs{I}>1$ and $h(\ell+E_v)=h(\ell)$ for every $v\in\grverts$. Then by stability property of $h$, $h(\ell+E_J)=h(\ell)$ for any $J\subset I$. Moreover, $w\bigl((\ell,I)\bigr)-w(\ell)=0=h(\ell+E_I)-h(\ell)$.

    Finally, assume that $\abs{I}>1$ and $h(\ell+E_v)>h(\ell)$ for a certain $v\in\grverts$. This means that $w(\ell+E_v)>w(\ell)$, hence $w\bigl((\ell,I)\bigr)>w(\ell)$ ($\dag$). Assume that $w\bigl((\ell,I)\bigr)$ is realized for a certain $J\subset I$. Let $J$ be minimal by this property. By ($\dag$) we know that $J\neq\emptyset$.

    By minimality of $J$, $w(\ell+E_{J\setminus u})<w(\ell+E_J)$ for any $u\in J$. By CDP $h(\ell+E_{J\setminus u})<h(\ell+E_J)$ too, hence by the stability of $h$ we also have $h(\ell)<h(\ell+E_u)$.

    In particular we found $u\in I$ such that $h(\ell)<h(\ell+E_u)$ (hence $h^\circ(\ell)=h^\circ(\ell+E_u)$) and $w\bigl((\ell, I)\bigr)=w\bigl((\ell+E_u, I^*)\bigr)$, where $I^*:= I\setminus u$. Now, we use induction applied for the cube $(\ell+E_v, I^*)$. In particular,
    \[
        w\bigl((\ell,I)\bigr)-w(\ell)=w\bigl((\ell+E_v, I^*)\bigr)-w(\ell)= w(\ell+E_v)+h(\ell+E_I)-h(\ell+E_v)-w(\ell)=h(\ell+E_I)-h(\ell)
    \]
    This ends the proof of \eqref{eq:KEYIDEN}.

    Let us denote by $\cellset$ the set of cubes of $R$. By Lemma \thmref{lem:eu-w} $\eu(\lcoh^*(R,w))$ is gven by
    \[
        \eu\bigl(\lcoh^*(R,w)\bigr)=\sum_{(\ell,I)\in\cellset}\, (-1)^{\abs{I}+1} w\bigl((\ell,I)\bigr).
    \]
    We will subdivide the lattice points of $R$ into the following disjoint subsets. For any $0\leq t\leq s$ we denote by $R_t$ the set of those points $\ell\in R\cap L$ for which the cardinality of $\setdef{v\in\grverts}{\ell_v=c_v}$ is $t$. If $\ell\in R_t$ set $I(\ell):=\setdef{v}{\ell_v<c_v}$. In particular, the set $\cellset$ is also a disjoint union of the sets $\setdef{(\ell, I)}{\ell\in R_t, I\subset I(\ell)}_t$. Then
    \[
        \eu\bigl(\lcoh^*(R,w)\bigr)=\sum_t\sum_{\ell\in R_t}\sum_{I\subset I(\ell)}(-1)^{\abs{I}+1} w\bigl((\ell,I)\bigr).
    \]
    If $0\leq t<s$ and $\ell\in R_t$ then $I(\ell)\neq\emptyset$. In this case for any $I$-independent (but maybe $\ell$-dependent) `constant' $a(\ell) $ we have $\sum\limits_{I\subset I(\ell)} (-1)^{\abs{I}+1}a(\ell)=0$. Therefore, by \refeqn{eq:KEYIDEN}, for any such $\ell$,
    \begin{equation}\label{eq:SUMt}
        \sum_{I\subset I(\ell)}\, (-1)^{\abs{I}+1} w\bigl((\ell,I)\bigr)=\sum_{I\subset I(\ell)}\, (-1)^{\abs{I}+1} h(\ell+E_I).
    \end{equation}
    On the other hand, if $t=s$, then $R_t=\setof{c}$, $I(c)=\emptyset$, hence $\sum\limits_{I\subset I(c)}(-1)^{\abs{I}+1}w\bigl((c,I)\bigr)=-w(c)$. Hence, corresponding to $t=s$ \refeqn{eq:SUMt} fails, i.e. it must be corrected by $-w(c)+h(c)=h^\circ (0)-h^\circ(c)$. Therefore,
    \begin{equation}\label{eq:SUMtt}
        \eu(\lcoh^*(R,w))=h^\circ (0)-h^\circ(c)\, +\, \sum_{(\ell,I)\in\cellset}\, (-1)^{\abs{I}+1} h(\ell+E_I).
    \end{equation}
    For any fixed $\wti\ell\in R$, consider the following summation over $\setdef*[\big]{(\ell,I)}{\ell+E_I=\wti\ell}$:
    \[
        S(\wti\ell):=\sum \, (-1)^{\abs{I}+1}h(\ell+E_I)=h(\wti{\ell})\cdot \sum \, (-1)^{\abs{I}+1}.
    \]
    Then, whenever the cardinality of $\setdef*[\big]{(\ell,I)}{\ell+E_I=\wti\ell}$ is $>1$, the above sum $\ds\sum(-1)^{\abs{I}+1}=0$. Thus, $S(\wti\ell)=0$ except when $\wti\ell=0$. But, for $\wti\ell=0$, $S(0)=-h(0)=0$ too. Thus $\ds\eu\bigl(\lcoh^*(R,w)\bigr)=h^\circ (0)-h^\circ(c)+\textstyle{\sum_{\,\wti\ell\,}}S(\wti\ell)=h^\circ (0)-h^\circ(c)$.
\end{pf}

\begin{rem}[label={rem:P0AN}]
    The very same argument as in the previous proof provides the following fact. Assume that instead of the rectangle $[0,c]$ we take $L_{\geq 0}$ (i.e., $c=\infty$), and we also have two functions $h,\,h^\circ\colon L_{\geq 0}\to\bbz$ with the very same properties (i) and (ii) as in \thmref{def:comblattice}. Furthermore, assume that both $h$ satisfies the stability property, and the pair $(h,h^\circ)$ satisfies the Combinatorial Duality Property (as in Theorem \thmref{th:comblattice}). Then \refeqn{eq:KEYIDEN} (whose proof is independent of the choice of $c$) implies
    \[
        \sum_{\ell\geq 0}\,\sum _I\, (-1)^{\abs{I}+1} w\bigl((\ell,I)\bigr)\, \Bf t^{\ell}=\sum_{\ell\geq 0}\,\sum _I\, (-1)^{\abs{I}+1} h(\ell+E_I)\,\Bf t^{\ell}.
    \]
\end{rem}

\begin{rem}
    Concrete examples show the following facts (see e.g. the $w$-table of Example \thmref{ex:2552}).

    Even if $h$ satisfies the path eu-coincidence (and $h^\circ =h^{sym}$), in general it is not true that $\lcoh^0(\gamma,w)$ is independent of the choice of the increasing path. (This fact remains valid even if we consider only the symmetric increasing paths, where a path $\gamma=\setof{x_i}_{i=0}^t$ is symmetric if $x_{t-i}=c-x_i$ for any $i$.)

    Even if $h$ satisfies both the path eu-coincidence and the eu-coincidence, in general it is not true that $\lcoh^*(R,w)$ equals any of the path lattice cohomologies $\lcoh^0(\gamma,w)$ associated with a certain increasing path. However, amazingly, all the Euler characteristics might agree (see below).
\end{rem}

\begin{rem}[label={ex:WCGP}]
    In the next discussion assume that $h^\circ =h^{sym}$.
    \begin{enumerate}[label={(\alph*)}]
        \item The CDP of $h$ implies the `path eu-coincidence' of $h$, see Lemma \thmref{lem:comblat}. However, the CDP of $h$ does not imply the `eu-coincidence' of $h$. As an example consider the second case from \thmref{rem:ketpelda}. Note that in this case the matroid or stabilization properties are not satisfied by $h$.
        
        \item On the other hand, an eu-coincidence type property cannot be hoped without (some type of) CDP. Indeed, in the next example $s=2$, $c=(2,2)$, and $h$ is associated with a graded vector space of dimension 2 supported in $(0,1)$ and $(1,2)$. In this case $h(c)=2$, $\eu(\lcoh^0(R,w))=0$ and for any symmetric increasing path $\eu(\lcoh^0(\gamma,w))=0$ too, and for any non-symmetric paths $\eu(\lcoh^0(\gamma,w))=1$.
        \begin{center}
            \setlength{\tabcolsep}{0.4em}
            \begin{tabular}{RRRRp{1cm}RRRR}
                    & 1 & 1 & 2 & &      & 1 & 0 & 0 \\
                h:  & 0 & 1 & 2 & & w_0: & 0 & 0 & 0 \\
                    & 0 & 1 & 2 & &      & 0 & 0 & 1
            \end{tabular}
        \end{center}
    \end{enumerate}
\end{rem}

\section{The (analytic) lattice cohomology of curves}\label{sec:curveslattice}

\subsection{Classical invariants, filtrations}\label{ss:anlattice-prelim}\hfill\smallbreak

Let $(C,o)$ be an isolated curve singularity with local algebra $\fnring=\fnring_{C,o}$. Let $\bigcup\limits_{i=1}^r(C_i,o)$ be the irreducible decomposition of $(C,o)$, denote the local algebra of $(C_i,o)$ by $\fnring_i$. We denote the integral closure of $\fnring_i$ by $\ol{\fnring_i}=\bbc\setof{t_i} $, and we consider $\fnring_i$ as a subring of $\ol{\fnring_i}$. Similarly, we denote the integral closure of $\fnring$ by $\ol\fnring=\oplus_i \bbc\setof{t_i}$. Let $\delta=\delta(C,o)$ be the delta invariant $\dim _\bbc\, \ol\fnring/\fnring$ of $(C,o)$.

We denote by $\Got v_i\colon\ol{\fnring_i}\to\ol\bbz_{\geq 0}=\bbz_{\geq 0}\cup\setof{\infty}$ the discrete valuation of $\ol{\fnring_i}$, where $\Got v_i(0)=\infty$. This restricted to $\fnring_i$ reads as $\Got v_i(f)=\ord_{t_i}(f)$ for $f\in \fnring_i$.

We define the semigroup (monoid) of $(C,o)$ as $\semigp=\Got v(\fnring)\cap(\bbz_{\geq 0})^r$, or
\[
    \semigp=\setdef{\Got v(f):=(\Got v_1(f),\ldots,\Got v_r(f))}{f\text{ is a nonzero divisor in $\fnring$}}\subset(\bbz_{\geq 0})^r.
\]
Let $\Got c=(\fnring:\ol\fnring)$ be the conductor ideal of $\ol\fnring$, it is the largest ideal of $\fnring$, which is an ideal of $\ol\fnring$ too. It has the form $(t_1^{c_1}, \ldots, t_r^{c_r})\ol\fnring$. The element $c_{\semigp}=(c_1,\ldots , c_r)$ is called the conductor of $\semigp$. From definitions $c_{\semigp}+(\bbz_{\geq 0})^r\subset \semigp$ and $c_{\semigp}$ is the smallest lattice point with this property, cf. \cite{delaMata}.

Assume that $(C,o)$ is the union of two (not necessarily irreducible) germs $(C',o)$ and $(C'',o)$ without common components, and fix some embedding $(C,o)\subset (\bbc^n,0)$. We define the \It{Hironaka's intersection multiplicity} of $(C',o)$ and $(C'',o)$ by $(C',C'')_{Hir}:=\dim ( \fnring_{\bbc^n,o}/ I'+I'')$, where $I'$ and $I''$ are the ideals in $\fnring_{(\bbc^n,0)}$ of $(C',o)$ and $(C'',o)$. Then one has the following formula \cite{Hironaka,BG80,Stevens85}:
\[
    \delta(C,o)=\delta(C',o)+\delta(C'',o)+ (C',C'')_{Hir}.
\]

From this it follows inductively that $\delta(C,o)\geq r-1$. In fact, $\delta(C,o)= r-1$ if and only if $(C,o)$ is analytically equivalent with the union of the coordinate axes of $\bbc^r,0)$, called \It{ordinary $r$-tuple} \cite{BG80}.

For plane curve germs $(C',C'')_{Hir}$ agrees with the usual intersection multiplicity at $(\bbc^2,0)$. In this case, the conductor entries are $c_i=2\delta(C_i,o)+ (C_i,\cup_{j\neq i} C_j)$. For a formula of $\setof*{c_i}_i$ in the general case see \cite{D'Anna}.

Now we introduce certain filtrations. We will focus on two vector spaces: the first is the infinite dimensional local algebra $\fnring$, the second one is the finite dimensional $\ol\fnring/\fnring$.

Consider $\grverts=\setof{1,\ldots, r}$ and the lattice $L:=\bbz^r$ with its natural basis $\setof{E_i}_i$. If $\ell=(\ell_1, \ldots, \ell_r)\in\bbz^r$, we set $\abs{\ell}:=\sum_i \ell_i$. Then $\ol\fnring$ has a natural filtration indexed by $\ell\in(\bbz_{\geq 0})^r$ given by $\ol\filt(\ell):=\setdef{g}{\Got v(g)\geq\ell}$. This induces an ideal filtration of $\fnring$ by $\filt(\ell):=\ol\filt(\ell)\cap\fnring\subset\fnring$, and also a filtration
\[
    \filt^\circ(\ell)=\frac{\ol\filt(\ell)+\fnring}{\fnring}=\frac{\ol\filt(\ell)}{\ol\filt(\ell)\cap\fnring}\subset\quo*{\ol{\fnring}}{\fnring}.
\]
The first filtration of $\fnring$ was considered in \cite{cdg2,cdg,Moyano,NPo}. (The second one is probably new in the literature.) Starting from $\filt$ one can define a multivariable Hilbert series. The second filtration provides a `multivariable sum-decomposition' of $\delta$. Set
\begin{align*}
    \hh(\ell) & =\dim \fnring/\filt(\ell)=\dim (\ol\filt(\ell)+\fnring)/\ol\filt(\ell), \\
    \ol{\hh}(\ell) & =\codim(\filt^\circ(\ell)\subset \ol\fnring/\fnring)=\dim(\ol\fnring/(\ol\filt(\ell)+\fnring)), \\
    \hh^\circ(\ell) & =\delta-\ol{\hh}(\ell).
\end{align*}
Then clearly $\hh$ and $\ol\hh$ are increasing, and $\hh(0)=\ol\hh(0)=0$. Moreover, $\ol{\hh}(\ell)=\delta$ for any $\ell\geq c_{\semigp}$ (since $\ol\filt(\ell)\subset\Got c\subset\fnring$ for such $\ell$). As such, $\hh^\circ$ is decreasing, takes $\delta$ at $0$, and $0$ for $\ell\geq c_{\semigp}$. We also have
\begin{equation}\label{eq:DUAL1}
    \hh(\ell)+\ol{\hh}(\ell)=\dim\, \ol\fnring/\ol\filt(\ell)=\abs{\ell}.
\end{equation}
In particular, $\hh(\ell+E_i)-\hh(\ell)\in\setof{0,1}$ for any $\ell$ and $i$.

With the notation $\Bf t=(t_1,\ldots, t_s)$ and $\Bf t^{\ell}=\prod_i t_i^{\ell_i}$, the multivariable Hilbert series is defined as the generating function $H(\Bf t)=\sum_{\ell\in(\bbz_{\geq 0})^r}\hh(\ell)\cdot{\Bf t}^\ell$.

The semigroup $\semigp$ and the Hilbert function $\setof{\hh(\ell)}_{\ell\geq 0}$ determine each other. Indeed, they are related as follows:
\begin{equation}\label{eq:Sfromh}
    \semigp=\setdef{\ell}{\hh(\ell+E_i)>\hh(\ell)\text{ for all $i$}}.
\end{equation}
Furthermore, for any $i\in\grverts$,
\begin{equation}\label{eq:hfromS}
    \hh(\ell+E_i)>\hh(\ell)\text{ exactly when there exists }s\in\semigp\text{ with }s\geq \ell, s_i=\ell_i.
\end{equation}
E.g., if $r=1$ then $\hh(\ell)=\#\setdef{s\in\semigp}{s<\ell}$.

\subsection{The weights and the lattice cohomology. Finiteness properties.}\label{ss:ANweightsCurves}\hfill\smallbreak

Note that $\bbz^r$ with its fixed basis $\setof{E_i}_i$, and the functions $h$ and $h^\circ$ satisfy all the requirements of \thmref{def:comblattice}. Set (cf. \eqref{eq:DUAL1})
\begin{equation}\label{eq:weightFormula}
    w_0(\ell)=\hh(\ell)+\hh^\circ(\ell)-\hh^\circ(0)=\hh(\ell)-\ol{\hh}(\ell)= 2\cdot\hh(\ell)-\abs{\ell}.
\end{equation}
The cubical decomposition of $(\bbr_{\geq 0})^r$ with the induced system of weight functions as in \thmref{def:comblattice} defines a lattice cohomology $\lcoh^*((\bbr_{\geq 0})^r,w)$. It will be denoted by $\lcoh^*(C,o)$. In parallel, definition \thmref{def:GRootW} gives a root $\Got R(C,o)$ as well.

Similarly, for any $c\in (\bbz_{\geq 0})^r$ we can consider the cubical decomposition of the rectangle $R(0,c)$ and the induced weight function. The corresponding lattice cohomology will be denoted by $\lcoh^*(R(0,c),w)$ and the root by $\Got R(R(0,c))$.

\begin{lem}[label={lem:INDEPCURVE}]
    For any $c\geq c_{\semigp}$, $c\in (\bbz_{\geq 0})^r$, the homotopy types of the spaces $\setof*{S_n}_n$ are independent of the choice of $c$, hence the lattice cohomology $\lcoh^*(R(0,c),w)$ is independent of the choice of $c$. Therefore, the natural restrictions induce graded $\bbz[U]$-module isomorphisms $\lcoh^*(R(0,c_{\semigp}),w)\simeq\lcoh^*(R(0,c),w)\simeq \lcoh^*(C,o)$.

    Similarly, $\Got R(R(0,c))$ is independent of $c$
    ($c\geq c_{\semigp}$),
    and $\Got R(R(0, c_{\semigp}))\simeq \Got R(R(0,c))\simeq \Got R(C,o)$.

    In particular, $\lcohred^*(C,o)$ has finite $\bbz$-rank, and $\Got R(C,o)$ is a graded root with all the required finiteness properties. Hence, the statements of Lemma \thmref{lem:eu-w} and Proposition \thmref{prop:HH} also hold.
\end{lem}
\begin{pf}
    Fix some $c\geq c_{\semigp}$ and choose $E_i$ with $c-E_i\geq c_{\semigp}$ too. Then for any $\ell\in R(0,c)\cap\bbz^s$ with $\ell_i=c_i$, by \eqref{eq:hfromS} we have $\hh(\ell-E_i)<\hh(\ell)$. Then from \eqref{eq:DUAL1} $\ol{\hh}(\ell-E_i)=\ol{\hh}(\ell)$. Thus $w_0(\ell-E_i)< w_0(\ell)$.

    Then for any $n\in \bbz$, a strong deformation retraction in the direction $E_i$ realizes a homotopy equivalence of the spaces $S_n\cap R(0,c-E_i) \subset S_n\cap R(0,c)$. A retraction $r:S_n\cap R(0,c)\to S_n\cap R(0,c-E_i)$ can be defined as follows. If $\square = (l,I)$ belongs to $S_n\cap R(0,c-E_i)$ then $r$ on $\square$ is defined as the identity. If $(l,I)\cap R(0,c-E_i)=\emptyset$, then $l_i=c_i$, and we set $r(x)=x-E_i$. Else, $\square =(l,I)$ satisfies $i\in I$ and $l_i=c_i-1$. Then we retract $(l,I)$ to $(l, I\setminus i)$ in the $i$-direction. The strong deformation retraction is defined similarly.
\end{pf}

Since $\hh$ is associated with a valuative filtration, it satisfies the matroid inequality, hence the stability condition too. Furthermore, as a consequence of \refeqn{eq:DUAL1}, the pair $(\hh,\hh^\circ)$ satisfies the Combinatorial Duality Property. In fact,
\begin{equation}\label{lem:CDPcurves}
    \bigl(\hh(\ell+E_i)-\hh(\ell)\bigr)+ \bigl(\ol{\hh}(\ell+E_i)-\ol{\hh}(\ell)\bigr)=1.
\end{equation}

Then the combinatorial theorem \thmref{th:comblattice} implies the following.

\begin{cor}[label={cor:EUcurves}]
    Fix any increasing $\gamma$ connecting 0 with a certain $c$, which satisfies $c\geq c_{\semigp}$. Then we have $\eu(\lcoh^*(C,o))=\eu(\lcoh^*(\gamma,w))=\delta(C,o)$.

    In particular, $\lcoh^*(C,o)$ is a `categorification' of the delta invariant of $(C,o)$.
\end{cor}

\subsection{Duality properties. Differential forms.}\label{rem:Rosenlicht}\hfill\smallbreak

In the case of normal surface singularities, $\hh^\circ$ has a reinterpretation in terms of 2-forms (with poles along the exceptional curve) \cite{Nkonyv}. In the present curve case we have a similar connection.

Write $n\colon\ol{(C,o)}\to(C,o)$ for the normalization. Let $\Omega^1(*)$ be the germs of meromorphic differential forms on the normalization $\ol{(C,o)}$ with pole (of any order) at most in $\ol {o}=n^{-1}(o)$. Let $\Omega^1_{\ol{(C,o)}}$ be the germs of regular differential forms at $\ol{(C,o)}$. The \It{Rosenlicht's regular differential forms} are defined as
\[
    \omega^R_{C,o}:=\setdef*[\Big]{\alpha\in \Omega^1(*)}{\sum_{p\in\ol{o}}\Rm{res}_p(f\alpha)=0\text{ for all }f\in\fnring}.
\]
(In fact one shows that it is canonically isomorphic with the dualizing module of Grothendieck associated with $(C,o)$.) Then, by \cite{BG80,Serre}, one has a perfect duality between $\omega^R_{C,o}/ n_*\Omega^1_{\ol{(C,o)}}$ and $\ol\fnring/\fnring=n_*\fnring_{\ol{(C,o)}}/ \fnring_{C,o}$:
\begin{equation}\label{eq:ROS}
    n_*\fnring_{\ol{(C,o)}}/ \fnring_{C,o}\ \times \ \omega^R_{C,o}/ n_*\Omega^1_{\ol{(C,o)}}\ \longrightarrow\ \bbc,\ \ \ [f]\times [\alpha]\mapsto \sum_{p\in\ol{o}} \, \res_p(f\alpha).
\end{equation}
Moreover, one can define a $\bbz^r$-filtration in $\omega^R_{C,o}/ n_*\Omega^1_{\ol{(C,o)}}$ such that the duality is compatible with the filtrations in $n_*\fnring_{\ol{(C,o)}}/ \fnring_{C,o}= \ol\fnring/\fnring$ and $\omega^R_{C,o}/ n_*\Omega^1_{\ol{(C,o)}}$.

\subsection{The Gorenstein case}\label{ss:GORduality}\hfill\smallbreak

By Serre \cite{Serre} or Bass \cite{Bass} (see also \cite{Huneke}) $(C,o)$ is Gorenstein if and only if $\dim(\ol\fnring/\fnring)=\dim(\fnring/\Got c)$. On the other hand, Delgado in \cite{delaMata} proved that the condition $\dim (\ol\fnring/\fnring)=\dim(\fnring/\Got c)$ is equivalent with the symmetry of the semigroup of values $\semigp$. If $r=1$ then the symmetry can be formulated easily: $\ell\in\semigp\ \Leftrightarrow \ c_{\semigp}-1-\ell\not\in\semigp$. If $r\geq 2$ then the definition is the following \cite{delaMata}.

For any $\ell\in\bbz^r$ and $i\in\setof{1,\cdots, r}$ set
\[
    \Delta_i(\ell)=\setdef{s\in\semigp}{s_i=\ell_i\text{ and }s_j>\ell_j\text{ for all }j\neq i},
    \quad\text{and}\quad
    \Delta(\ell):=\cup_i\Delta_i(\ell).
\]
Then $\semigp$ is called symmetric if $\ell\in\semigp\ \Leftrightarrow \ \Delta(c_{\semigp}-\Bf{1}-\ell)=\emptyset$. (Here $\Bf{1}=(1,\ldots, 1)$.)

If $(C,o)$ is Gorenstein then in \cite{cdk} (see also \cite{Moyano}) is proved that
\[
    \bigl(\hh(\ell+E_i)-\hh(\ell)\bigr)+\bigl(\hh(c_{\semigp}-\ell)-
    \hh(c_{\semigp}-\ell-E_i)\bigr)=1,
\]
or
\[
    \hh(\ell)-\hh(c_{\semigp}-\ell)=\sum_i \ell_i-\delta.
\]
This combined with \refeqn{eq:DUAL1} gives $\hh(c_{\semigp}-\ell)=h^\circ(\ell)$, or $\hh^\circ=h^{sym}$ (where $h^{sym}$ is associated with $c_{\semigp}$).

In fact, we can show even more. Recall that $\omega^R_{C,o}$ is an $\fnring$-module by multiplication by functions. Since $\Got c\ol\fnring\subset\fnring$, we also have $\Got c\omega^R_{C,o}\subset \Omega^1_{\ol{C,o}}$, hence $\omega^R_{C,o}/\Omega^1_{\ol{C,o}}$ is an $\fnring/\Got c$-module. The Gorenstein property means that it is a free rank one module, that is, there exists $[\alpha_{Gor}]\in \omega^R_{C,o}/\Omega^1_{\ol{C,o}}$ so that
\[
    \omega^R_{C,o}/\Omega^1_{\ol{C,o}}=[\alpha_{Gor}]\cdot \fnring/\Got c.
\]
Using the perfect pairing from \refeqn{eq:ROS}, we get a perfect pairing
\begin{equation}\label{eq:ROS2}
    D\colon\ol\fnring/\fnring\ \times \ \fnring/\Got c\to\bbc,\ \ \ [f]\times [g]\mapsto\sum_{p\in\ol{o}}\,\res_p(fg\cdot\omega_{Gor}).
\end{equation}
On both $\ol\fnring/\fnring$ and $\fnring/\Got c$ we have filtrations indexed by $\ell\in[0, c_{\semigp}]\cap \bbz^r$. On $\ol\fnring/\fnring$, it is $\filt^\circ$ defined above. On $\fnring/\Got c$, it is induced by $\filt$ considered above, namely $\setof{\filt^c(\ell)=\filt(\ell)/\Got c}_{\ell\in R}$ (note that $\Got c=\filt(c_{\semigp})\subset\filt(\ell)$ for any $\ell\in R(0,c_{\semigp})$). Then we claim that $\filt^\circ$ and $\filt^c$ are dual in the following sense: for any $\ell\in R(0,c_{\semigp})$, the pairing $D(f,g)=0$ for any $f\in\filt^\circ(\ell)$ if and only if $g\in\filt^c(c_{\semigp}-\ell)$.

This implies that $\codim\filt^\circ(\ell)+\codim\filt^c(c_{\semigp}-\ell)=\delta$, or, again, $\hh^\circ=\hh^{sym}$.

Thus, in the Gorenstein case, the weight function on $[0,c_{\semigp}]$ is $w_0(\ell)=\hh(\ell)+\hh^{sym}(\ell)-\delta$, and
\begin{equation}\label{eg:Gorhsim}
    \lcoh^*(C,o)=\lcoh^*(R(0,c_{\semigp}),\hh+\hh^{sym}-\delta).
\end{equation}
Hence, for Gorenstein singularities, the involution of $R(0,c_{\semigp})$, $\ell\mapsto c_{\semigp}-\ell$, induces an involution --- a degree zero $\bbz[U]$-module morphism of order two --- on $\lcoh^*(C,o)$.

\begin{rem}
    \begin{enumerate}[label={(\alph*)}]
        \item In the non-Gorenstein case the pair $(\hh,\hh^{sym})$ does not satisfy the CDP. Take e.g. for the irreducible curve the image of the parametrization $(\bbc,0)\to (\bbc^3,0)$, $t\mapsto (t^3, t^5, t^7)$ (that is, the Hilbert function associated with the semigroup $\gen{3,5,7}\subset \bbz_{\geq 0}$).

        \item In the non-Gorenstein case the lattice cohomology (associated with the pair $(h, h^\circ)$ as above) usually does not admit any $\bbz_2$-involution (see the example from part (a)).
    \end{enumerate}
\end{rem}

\begin{ex}[label={ss:planecurves}, description={Plane curve singularities}]
    Recall that for any curve the lattice cohomology and graded root are determined by the Hilbert function $\hh$. On the other hand, in the case of plane curves, the Hilbert function can be recovered from the embedded topological type of the link $L_{(C,o)}=\sqcup_{i=1}^r S^1\subset S^3$ of $(C,o)\subset (\bbc^2,o)$. This is what we explain next.
    
    Recall that any embedded link with $r$ components in $S^3$ has a multivariable Alexander polynomial $\Delta(t_1, \ldots, t_r)$. It is convenient to renormalize it, we will call the resulting object `Poincaré series' (usually the Poincaré series are introduced via the Hilbert series of $(C,o)$, here it is convenient to connect them directly to the Alexander polynomials, for details see e.g. \cite{cdg2,cdg,Nkonyv}). The Poincaré series $P(\Bf t)=P_{(C,o)}(\Bf t)$ and $\Delta(\Bf t)$ are connected as follows. If $r=1$ then $P(t)(1-t)=\Delta(t)$, while if $r>1$ then $P(\Bf t)=\Delta(\Bf t)$. If $r=1$ then we also have the identification with the semigroup
    \begin{equation}\label{eq:semigroupDelta}
        P(t)=\frac{\Delta(t)}{1-t}=\sum_{s\in\semigp_{C,o}}t^s.
    \end{equation}
    
    For any nonempty subset $I=\setof{i_1,\ldots, i_{\abs{I}}}\subset\setof{1,\ldots,r}$ consider the curve germ $(C_I,o)=(\cup_{i\in I}C_i,o)$. It has an Alexander polynomial and a Poincaré series $P_{(C_I,o)}$ connected by the above principle. For $I=\emptyset$ we set $P_{\emptyset}=0$. Then the next formula recovers the Hilbert series $H(t_1,\ldots , t_r)$ of $(C,o)$ in terms of the Poincaré series $\setof*{P_{C_I}(\Bf t)}_I$ (or Alexander polynomials) (see \cite{GorNem2015}):
    \begin{equation}\label{hilbert}
        H(t_1,\ldots,t_r)=
            \frac{1}{\prod_{i=1}^{r}(1-t_i)}\cdot
            \sum_{I\subset\setof{1,\ldots,r}}(-1)^{\abs{I}-1}
                \Big(\prod_{i\in I}t_i\Big)\cdot P_{C_I}(t_{i_1},\ldots,t_{i_{|I|}}).
    \end{equation}
    For example, in the case of the plane curve singularity $A_{2n-1}$, defined by the equation $x^2-y^{2n}=0$, $P(t_1,t_2)=1+t_1t_2+\cdots+(t_1t_2)^{n-1}$, $P_{\setof{x\pm y^n=0}}(t_i)=1/(1-t_i)$, hence
    \[
        H(t_1,t_2)=\frac{1}{(1-t_1)(1-t_2)}\left(\frac{t_1}{1-t_1}+\frac{t_2}{ 1-t_2}-t_1t_2(1+\ldots+(t_1t_2)^{n-1})\right).
    \]
    
    Recall that for plane curve singularities, by \cite{Yamamoto} the multi-variable Alexander polynomial $\Delta$ of $(C,o)$ determines the embedded topological type of $(C,o)$, hence the multi-variable Alexander polynomials for all $C_I$ as well. In particular, by the above correspondence between the Alexander polynomials and Poincaré series, $\Delta$, or $P_{(C,o)}$, determines the collection $\setof*{P_{(C_I,o)}}_I$ as well. A direct reduction procedure from $P_C$ to $P_{C_I}$ is given by Torres formula \cite{torres} as well:
    \begin{equation*}\label{eq:redP}
        P_{C_2\cup\cdots \cup C_s}(t_2,\ldots,t_s)=
        P(t_1,\ldots,t_s)|_{t_1=1}\cdot \big(1-t_2^{(C_1,C_2)}\cdots t_s^{(C_1,C_s)}\big)^{-1},
    \end{equation*}
    where $(C_i,C_j)$ is the intersection multiplicity at the origin of $C_i$ and $C_j$, $i\not=j$.
    
    For a formula of $P(\Bf t)$ in terms of the embedded resolution graph of $(C,o)$ see \cite{EN}.
    
    We emphasize again, for plane curve singularities the invariants $H$, $P$, $\Delta$ and $\semigp$ determine each other and are complete embedded topological invariants of $(C,o)$.
\end{ex}

\vspace*{-3mm}
\subsection{Example. Newton non-degenerate isolated singularities}\label{ss:planecurves} \hfill\smallbreak

The following construction is valid for any hypersurface singularity in $(\bbc^{n+1},0)$ with non-degen\-er\-ate Newton principal part. For simplicity we present here the case of plane curve singularities (case $n=1$), for the general case see \cite{Nkonyv}.

\begin{ex}[label={ex:latgrfilt2}, description={Combinatorial
lattice cohomology determined by a Newton diagram}] \\
    Assume that $\Gamma$ is a convenient Newton diagram (boundary) in $\bbr_{\geq 0}^{2}$. (`Convenient' means that $\Gamma$ intersects both coordinate axes; for more details and notations regarding Newton diagrams see e.g. \cite{Nkonyv}). Let $\Gamma^{(1)}$ denote the $1$-dimensional faces (segments) of $\Gamma$ (which might contain interior lattice points). For each $\sigma\in \Gamma^{(1)}$ let $N_\sigma$ be the normal primitive integral vector of the corresponding face $F_\sigma$ with all entries positive. Set $m_\sigma$ in such a way that $F_\sigma$ is on the affine line $\gen{N_\sigma,p}=m_\sigma$. Let $M$ be the $\bbc$-vector space generated by lattice points $p\in (\bbz_{>0})^{2}$, which sit either on or below $\Gamma$. This means that $p\in (\bbz_{>0})^{2}$ should satisfy $\gen{N_\sigma,p}\leq m_\sigma$ for at least one $\sigma$. Then $M\neq 0$ if and only if $\Bf{1}=(1,1)$ is such a point; we will assume this fact from here on.

    We order $\Gamma^{(1)}$ as $\setof{\sigma_1,\ldots,\sigma_s}$. We introduce a $\bbz^s$-grading of $M$ by
    \[
        \deg(p)=\bigl(\gen{N_{\sigma_1}, p-\Bf{1}},\ldots,\gen{N_{\sigma_s}, p-\Bf{1}}\bigr)\in\bbz^s.
    \]
    Define $c_{\Gamma}=(c_1,\ldots , c_s)$ by $c_i:= m_{\sigma_i}+1-\gen{N_{\sigma_i}, \Bf{1}}$.

    Since $p\geq \Bf{1}$ for any such point $p\in M$, we also have $\deg(p)\in \bbz_{\geq 0}^s$. Moreover, for any $p\in M$ there exists $\sigma_i$ with $\gen{N_{\sigma_i}, p}\leq m_{\sigma_i}$, hence $\gen{N_{\sigma_i}, p-\Bf{1}} <c_i$. Therefore, $\deg(p)\ngeq c_{\Gamma}$.

    In particular, the conditions of \thmref{ex:latgrfilt} are satisfied
    with $c=c_{\Gamma}$.

    For any lattice point from $R(0,c)$ set $h(\ell)=\sum _{\Bf a\not\geq\ell}\dim\, M_{\Bf a}$ as in \thmref{ex:latgrfilt}, hence it satisfies the stability property. Consider also $h^\circ=h^{sym}$ associated with the above $c_{\Gamma}$.

    We claim that if there is no face $F_{\sigma}$ whose relative interior contains a lattice point $p$ then the pair $(h,h^{sym})$ satisfies the CDP. Indeed, assume that this is not the case. Then there exists $\ell$ and $1\leq i\leq s$ such that $h(\ell)\not\not=h(\ell+E_i)$ and $h(c_{\Gamma}-\ell)\not\not=h(c_{\Gamma}-\ell-E_i)$. By the first fact there exists $p\in (\bbz_{>0})^2$ such that $\inn{N_{\sigma_j}}{p-\Bf 1}\geq\ell_j$ for every $j$ and $\inn{N_{\sigma_i}}{p-\Bf 1}=\ell_i$. By the second one there exists $q\in(\bbz_{>0})^2$ such that $\inn{N_{\sigma_j}}{q-\Bf 1}\geq (c_{\Gamma})_j-\ell_j$ for every $j\not =i$ and $\inn{N_{\sigma_i}}{q-\Bf 1} = (c_{\Gamma})_i-1-\ell_i$. In particular, $\inn{N_{\sigma_j}}{p+q-\Bf 1}>m_{\sigma_j}$ for every $j\neq i$ and $\inn{N_{\sigma_i}}{p+q-\Bf 1}=m_{\sigma_i}$. But then $p+q-\Bf 1$ is a relative interior lattice point of $F_{\sigma_i}$, a fact which contradicts the assumption.

    Let us call the lattice cohomology associated with the rectangle $R(0, c_{\Gamma})$ and the weight function $\ell\mapsto w_0(\ell)=h(\ell)+h^{sym}(\ell)-h^{sym}(0)$ the \It{combinatorial lattice cohomology associated with the Newton diagram $\Gamma$}. It will be denoted by $\lcoh^*_{\Gamma}(R,w)$. Similarly, for any increasing path $\gamma$ connecting $0$ and $c_{\Gamma}$ we consider the \It{combinatorial path lattice cohomology associated with the Newton diagram $\Gamma$}, denoted by $\lcoh^*_{\Gamma}(\gamma,w)$.

    Then, by the above discussion and by Theorem \thmref{th:comblattice}, if in the relative interior of the faces $F_\sigma$ there are no lattice points, we have
    \begin{equation}\label{eq:NGamma}
        \eu(\lcoh^*_{\Gamma}(\gamma,w)) = \eu(\lcoh^*_{\Gamma}(R,w)) = \dim(M).
    \end{equation}
    Now, any such Newton boundary determines (an equisingularity type) of an isolated plane curve singularity $(C,o)$. E.g., a possible equation is a generic linear combination of monomials sitting on $\Gamma$, and one can also add with any coefficient any monomial above $\Gamma$. The embedded topological type of the germ does not depend on the choice of the coefficients (chosen in this way), and it can be determined from $\Gamma$. The number of irreducible components of $(C,o)$ is the combinatorial length of $\Gamma$, and all the components have exactly one Puiseux pair. The condition that `the relative interior of the faces $F_\sigma$ contain no lattice points' reads as all the irreducible components have different first Puiseux pair, and their number is in bijection with the number of faces $\setof{F_{\sigma}}_\sigma$. (Hence $r=s$.) Furthermore, $\dim (M)=\delta_{(C,o)}$ \cite{SaitoNN}, and $c_{\Gamma}$ is the conductor $c_{\semigp}$ \cite{Nkonyv}.
\end{ex}

\begin{ex}[label={ex:2552}]
    Consider the situation of \thmref{ex:latgrfilt2} with $\Gamma$ generated by the lattice points $(0,7), (2,2)$, and $(7,0)$. In fact, this is the Newton diagram of the Newton non-degenerate plane curve singularity $f=(x^2+y^5)(y^2+x^5)$. The normal vectors $(2,5)$ and $(5,2)$ of the two faces provide the degrees. The next diagram shows the points $p$, which generate $M$ as base elements, and their degrees.
    
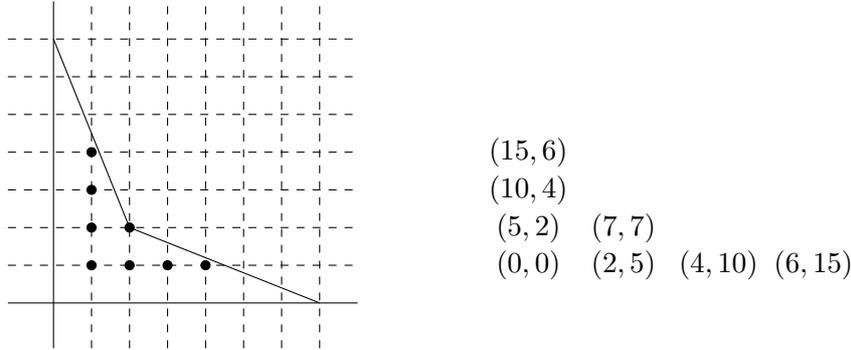
\begin{center}
        \begin{tikzpicture}[scale=0.5]
            \draw[fill] (1,1) circle [radius=0.12];
            \draw[fill] (2,1) circle [radius=0.12];
            \draw[fill] (3,1) circle [radius=0.12];
            \draw[fill] (4,1) circle [radius=0.12];
            \draw[fill] (1,2) circle [radius=0.12];
            \draw[fill] (2,2) circle [radius=0.12];
            \draw[fill] (1,3) circle [radius=0.12];
            \draw[fill] (1,4) circle [radius=0.12];

            \draw (-1.2,0) -- (8,0);
            \draw (0,-1.2) -- (0,8);
            \draw[dashed] (-1.2,1) -- (8,1);
            \draw[dashed] (-1.2,2) -- (8,2);
            \draw[dashed] (-1.2,3) -- (8,3);
            \draw[dashed] (-1.2,4) -- (8,4);
            \draw[dashed] (-1.2,5) -- (8,5);
            \draw[dashed] (-1.2,6) -- (8,6);
            \draw[dashed] (-1.2,7) -- (8,7);
            \draw[dashed] (1,-1.2) -- (1,8);
            \draw[dashed] (2,-1.2) -- (2,8);
            \draw[dashed] (3,-1.2) -- (3,8);
            \draw[dashed] (4,-1.2) -- (4,8);
            \draw[dashed] (5,-1.2) -- (5,8);
            \draw[dashed] (6,-1.2) -- (6,8);
            \draw[dashed] (7,-1.2) -- (7,8);

            \draw (0,7) -- (2,2) -- (7,0);

            \node at (12.5,1) {$(0,0)$};
            \node at (15.0,1) {$(2,5)$};
            \node at (17.5,1) {$(4,10)$};
            \node at (20.0,1) {$(6,15)$};
            \node at (12.5,2) {$(5,2)$};
            \node at (15.0,2) {$(7,7)$};
            \node at (12.5,3) {$(10,4)$};
            \node at (12.5,4) {$(15,6)$};
        \end{tikzpicture}
    \end{center}

    The vector $c_{\Gamma}$ is $(8,8)$. It is the conductor of $f$. Furthermore, $\dim (M)=8$ is the delta-invariant $\delta(f)$ of $f$. The next tables show the function $h$ and the weight function $w_0$ (induced by $\Gamma$) on $R$.
    \begin{center}
        \setlength{\tabcolsep}{0.6em}
        \def\w#1 {\makebox[0pt][r]{$#1$}}
        \begin{tabular}{R|RRRRRRRRR}
            8   & 6 & 6 & 6 & 6 & 6 & 7 & 7 & 8 & 8 \\
            7   & 5 & 5 & 5 & 5 & 5 & 6 & 6 & 7 & 8 \\
            6   & 4 & 4 & 4 & 4 & 4 & 5 & 5 & 6 & 7 \\
            5   & 3 & 3 & 3 & 4 & 4 & 5 & 5 & 6 & 7 \\
            4   & 2 & 2 & 2 & 3 & 3 & 4 & 4 & 5 & 6 \\
            3   & 2 & 2 & 2 & 3 & 3 & 4 & 4 & 5 & 6 \\
            2   & 1 & 1 & 1 & 2 & 2 & 3 & 4 & 5 & 6 \\
            1   & 1 & 1 & 1 & 2 & 2 & 3 & 4 & 5 & 6 \\
            0   & 0 & 1 & 1 & 2 & 2 & 3 & 4 & 5 & 6 \\
            \hline
                & 0 & 1 & 2 & 3 & 4 & 5 & 6 & 7 & 8 \\
        \end{tabular}
        \hskip 1cm
        \begin{tabular}{R|RRRRRRRRR}
            8 & 4 & 3   & 2           & 1          & 0           & 1          & 0           & 1   & 0 \\
            7 & 3 & 2   & 1           & 0          &\w-1         & 0          &\w-1         & 0   & 1 \\
            6 & 2 & 1   & 0           & \w-1       &\w\Boxed{-2} &\w-1        &\w\Boxed{-2} &\w-1 & 0 \\
            5 & 1 & 0   &\w-1         & 0          &\w-1         &\circled{0} &\w-1         & 0   & 1 \\
            4 & 0 &\w-1 &\w\Boxed{-2} & \w-1       &\w\Boxed{-2} &\w-1        &\w\Boxed{-2} &\w-1 & 0 \\
            3 & 1 & 0   &\w-1         &\circled{0} &\w-1         & 0          &\w-1         & 0   & 1 \\
            2 & 0 &\w-1 &\w\Boxed{-2} &\w-1        &\w\Boxed{-2} &\w-1        & 0           & 1   & 2 \\
            1 & 1 & 0   &\w-1         & 0          &\w-1         & 0          & 1           & 2   & 3 \\
            0 & 0 & 1   & 0           & 1          & 0           & 1          & 2           & 3   & 4 \\
            \hline
              & 0 & 1   & 2           & 3          & 4           & 5          & 6           & 7   & 8 \\
        \end{tabular}
    \end{center}
    \nni
    The circle points the generators of $\lcoh^1_{\Gamma}(R,w)=\laurpols_{-2}(1)^2$, while the boxes the local minima of $w_0$. Hence $\lcoh^0_{\Gamma}(R,w)=\laurpols_{-4}^+\oplus \laurpols_{-4}(1)^6\oplus\laurpols_0(1)^2$ and the combinatorial graded root is
    
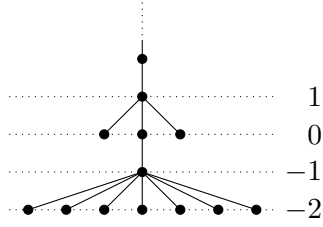
\begin{center}
        \begin{tikzpicture}[scale=0.5]
            \draw[fill] (0,2) circle [radius=0.12];
            \draw[fill] (0,1) circle [radius=0.12];
            \draw[fill] (-1,0) circle [radius=0.12];
            \draw[fill] (0,0) circle [radius=0.12];
            \draw[fill] (1,0) circle [radius=0.12];
            \draw[fill] (0,-1) circle [radius=0.12];
            \draw[fill] (-3,-2) circle [radius=0.12];
            \draw[fill] (-2,-2) circle [radius=0.12];
            \draw[fill] (-1,-2) circle [radius=0.12];
            \draw[fill] (0,-2) circle [radius=0.12];
            \draw[fill] (1,-2) circle [radius=0.12];
            \draw[fill] (2,-2) circle [radius=0.12];
            \draw[fill] (3,-2) circle [radius=0.12];

            \draw[dotted] (0,3.5) -- (0,2.5);
            \draw (0,2.5) -- (0,-2);
            \draw (0,1) -- (-1,0);
            \draw (0,1) -- (1,0);
            \draw (0,-1) -- (-3,-2);
            \draw (0,-1) -- (-2,-2);
            \draw (0,-1) -- (-1,-2);
            \draw (0,-1) -- (1,-2);
            \draw (0,-1) -- (2,-2);
            \draw (0,-1) -- (3,-2);

            \draw[dotted] (-3.5,1) -- (3.5,1);
            \draw[dotted] (-3.5,0) -- (3.5,0);
            \draw[dotted] (-3.5,-1) -- (3.5,-1);
            \draw[dotted] (-3.5,-2) -- (3.5,-2);

            \node[left] at (5,1) {$1$};
            \node[left] at (5,0) {$0$};
            \node[left] at (5,-1) {$-1$};
            \node[left] at (5,-2) {$-2$};
        \end{tikzpicture}
    \end{center}
    \nni
    $\min w_0=-2$, $\rank\lcoh_{\Gamma,red}^0(R,w)=8$, $\rank\lcoh_{\Gamma}^1(R,w)=2$, and $\eu(\lcoh_{\Gamma}^*(R,w))=8=\dim (M)=h(c_{\Gamma})=\delta(C,o)$. In particular, we constructed a set of cohomology groups $\lcoh_{\Gamma}^*(R, w)$ whose Euler characteristic is the delta-invariant. This is a combinatorial `categorification of $\delta$'.

    From analytic point of view, we can consider the valuative Hilbert function on $R$ as well, and the corresponding weight function and analytic lattice cohomology $\lcoh^*(C,o)$. For the concrete example above it turns out that the $h$-function determined by $\Gamma$ and the $h$-function determined by the normalization procedure (Hilbert function) coincide. Hence the weight functions coincide as well (being defined by symmetrization of the $h$-function). In particular, $\lcoh^*(C,o)=\lcoh^*_{\Gamma}(R,w)$.

    It is not hard to see that such coincidences hold for any Newton diagram with $r=s=1$, hence in all these cases $\lcoh^*(C,o)=\lcoh^*_{\Gamma}(R,w)$.

    However, for general $\Gamma$ with $r=s>1$, the two $h$-functions defined on $R(0,c_{\Gamma})=R(0,c_{\semigp})$ do not coincide, and the weight functions differ too. (In fact, $\filt(\ell)\supset \filt_{\Gamma}(\ell)$, hence $h_{(C,o)}(\ell)\leq h_{\Gamma}(\ell)$,
    but the equality might fail.) But, according to \cite{NS}, the corresponding two lattice cohomologies agree again: $\lcoh^*(C,o)=\lcoh^*_{\Gamma}(R,w)$ (and the graded roots as well). For a concrete example when $h_{(C,o)}\not=h_{\Gamma}$ the reader might write the $h$-tables of the germ $(x^3+y^5)(x^5+y^3)$. For more see \cite{NS}.
\end{ex}

\subsection{More examples}\label{ss:moreexamples}\hfill\smallbreak

\begin{ex}[label={ex:smooth}]
    If $(C,o)$ is smooth then $r=1$, $h(\ell)=w_0(\ell)=\ell$ for any $\ell\geq 0$. Hence each $S_n$ is contractible and ${\Got R}(C,o)={\Got R}_{(0)}$ and $\lcoh^*(C,o)=\laurpols_{0}^+$.

    If $(C,o)$ is not smooth then ${\Got R}(C,o)\not={\Got R}_{(0)}$ hence $\lcoh^0(C,o)\not=\laurpols_{0}^+$. Indeed, if $r=1$ then $w_0(0)=0$, $w_0(1)=1$, $w_0(2)=0$, hence $S_0$ has at least two connected components. Similarly, if $r\geq 2$ then $w_0(0)=0$, $w_0(E_i)=1$ for any $i$, and $w_0(1, \ldots, 1)=2-r\leq 0$ (since $\filt(1,\ldots, 1)$ is the maximal ideal). Hence, again $S_0$ has more components.
\end{ex}
\begin{ex}[label={ex:axes}]
    Assume that $(C,o)$ is an ordinary $r$-tuple ($r\geq 2$), that is, $\delta (C,o)=r-1$. In this case we can assume that $(C,o)$ is the union of the coordinate axes of $(\bbc^r,0)$. A computation shows that $c_{\semigp}=(1,\ldots, 1)$, hence all the important information is in the unit cube. For any $I\subset \setof{1,\ldots r}$, $I\neq\emptyset$, $\filt(E_I)$ is the maximal ideal, hence $h(E_I)=1$ and $w_0(E_I)=2-|I|$. Therefore, ${\Got R}(C,o)$ has the form
    
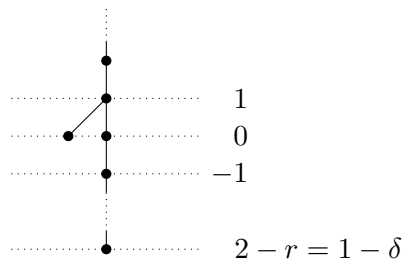
\begin{center}
        \begin{tikzpicture}[scale=0.5]
            \draw[fill] (0,2) circle [radius=0.12];
            \draw[fill] (0,1) circle [radius=0.12];
            \draw[fill] (-1,0) circle [radius=0.12];
            \draw[fill] (0,0) circle [radius=0.12];
            \draw[fill] (0,-1) circle [radius=0.12];
            \draw[fill] (0,-3) circle [radius=0.12];

            \draw (0,1) -- (0,2.5);
            \draw[dotted] (0,2.5) -- (0,3.5);
            \draw (0,1) -- (-1,0);
            \draw (0,1) -- (0,-1.5);
            \draw[dotted] (0,-1.5) -- (0,-2.5);
            \draw (0,-2.5) -- (0,-3);

            \draw[dotted] (-2.5,1) -- (2.5,1);
            \draw[dotted] (-2.5,0) -- (2.5,0);
            \draw[dotted] (-2.5,-1) -- (2.5,-1);
            \draw[dotted] (-2.5,-3) -- (2.5,-3);

            \node[left] at (4.0,1) {$1$};
            \node[left] at (4.0,0) {$0$};
            \node[left] at (4.0,-1) {$-1$};
            \node[right] at (3.1,-3) {$2-r=1-\delta$};
        \end{tikzpicture}
    \end{center}
    with $\min(w)=2-r=1-\delta$, $\lcoh^0(C,o)=\laurpols_{2(2-r)}^+\oplus\laurpols_{0}(1)$, $\lcoh^{\geq 1}(C,o)=0$.
\end{ex}

\begin{ex}[label={ex:semigr345}]
    Consider the irreducible curve $(C,o)$ given by parametrization $t\mapsto (t^3,t^4, t^5)$. It has semigroup $\semigp=\langle 3,4,5\rangle \subset \bbz_{\geq 0}$. It has $\delta=2$ and its graded root and lattice cohomology coincide with the root and lattice cohomology of the $3$-tuple (which has $r=3$). In particular ${\Got R}(C,o)$ and $\lcoh^*(C,o)$ do not determine the number of irreducible components of $(C,o)$.

    Moreover, in Examples \thmref{ex:Teissier} and \thmref{ex:Teissier2} we will also see that ${\Got R}(C,o)$ and $\lcoh^*(C,o)$ do not determine the embedding dimension of $(C,o)$ either.
\end{ex}
\begin{ex}[label={ex:semigrgen}]
    Let $\semigp\subset \bbz_{\geq 0}$ be a subsemigroup with $0\in \semigp$ (i.e. monoid) with $\bbz_{\geq 0}\setminus\semigp$ finite. Then there exists an irreducible curve $(C,o)$ whose semigroup is $\semigp$. Indeed, let $\bar{\beta}_0,\ldots,\bar{\beta}_s$ be a minimal set of generators of $\semigp$. Let $C^\semigp$ be the affine curve defined via the parametrization $t\mapsto (t^{\bar{\beta}_0},\ldots, t^{\bar{\beta}_s})$. Then the affine coordinate ring $\bbc[C^\semigp]$ of $C^\semigp$ is the image in $\bbc[t]$ of the morphism $\phi\colon\bbc[u_0,\ldots, u_s]\to\bbc[t]$, $\phi(u_i)=t^{\bar{\beta}_i}$. That is, the coordinate ring of the affine curve $C^{\semigp}$ equals the monoid ring of $\semigp$. The germ $(C^\semigp,0)$ is irreducible, and its semigroup is $\semigp$.

    In particular, any monoid of $\bbz_{\geq 0}$ (with $\bbz_{\geq 0}\setminus\semigp$ finite) admits a graded root and a lattice cohomology via $(C^{\semigp},o)$.
\end{ex}

\begin{ex}[label={ex:realizability}]
    We believe that the following question is crucial regarding the classification of curve singularities.

    \begin{qtn}[nonum]
        \begin{enumerate}[label={(\alph*)}]
            \item Which graded roots can be realized as ${\Got R}(C,o)$ ?
            \item Which graded $\bbz[U]$-modules can be realized as $\lcoh^*(C,o)$?
        \end{enumerate}
    \end{qtn}

    \nni
    For example, the following graded root cannot be realized as ${\Got R}(C,o)$ by some irreducible $(C,o)$.
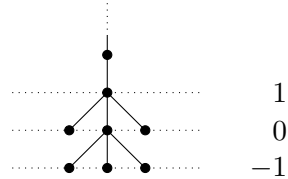
    \begin{center}
        \begin{tikzpicture}[scale=0.5]
            \draw[fill] (0,2) circle [radius=0.12];
            \draw[fill] (0,1) circle [radius=0.12];
            \draw[fill] (-1,0) circle [radius=0.12];
            \draw[fill] (0,0) circle [radius=0.12];
            \draw[fill] (1,0) circle [radius=0.12];
            \draw[fill] (1,-1) circle [radius=0.12];
            \draw[fill] (0,-1) circle [radius=0.12];
            \draw[fill] (-1,-1) circle [radius=0.12];

            \draw (0,1) -- (0,2.5);
            \draw[dotted] (0,2.5) -- (0,3.5);
            \draw (0,1) -- (-1,0);
            \draw (0,1) -- (1,0);
            \draw (0,1) -- (0,-1);
            \draw (0,0) -- (1,-1);
            \draw (0,0) -- (-1,-1);

            \draw[dotted] (-2.5,1) -- (2.5,1);
            \draw[dotted] (-2.5,0) -- (2.5,0);
            \draw[dotted] (-2.5,-1) -- (2.5,-1);

            \node[left] at (5.0,1) {$1$};
            \node[left] at (5.0,0) {$0$};
            \node[left] at (5.0,-1) {$-1$};
        \end{tikzpicture}
    \end{center}
\end{ex}

\subsection{First connections with flat deformation of curve singularities}\label{ss:exdef}\hfill\smallbreak

If $(C,o)$ is an isolated curve singularity with $r$ irreducible components, let its lattice $\bbz^r$ be denoted by $L(C,o)$.

\begin{thm}[label={th:DEF}]
    Assume that we have a flat deformation of isolated plane curve singularities $\setof{(C_t,o)}_{t\in(\bbc,0)}$ with $(C_{t=0},o)$ irreducible. Then the linear map connecting the corresponding lattices $\phi\colon L(C_{t\not=0})=\bbz^r\to L(C_{t=0})=\bbz$, $\ell\mapsto \phi(\ell)=\abs\ell=\sum_i\ell_i$ induces a degree zero graded $\bbz[U]$-module morphism $\lcoh^*(C_{t=0},o)\to \lcoh^*(C_{t\not=0},o)$, and similarly a graded (graph) map of degree zero at the level of graded roots ${\Got R}(C_{t\not=0},o)\to{\Got R}(C_{t=0},o)$.
\end{thm}
\begin{pf}
    In \cite{GorNem2015} it is proved that the map $\phi$ from the statement satisfies $\hh_{t\not=0}(\ell)\geq \hh_{t=0}(\abs\ell)$. In particular, $w_{t\not=0}(\ell)\geq w_{t=0}(|\ell|)$ too. Hence, for any $n\in\bbz$, if $\ell\in S_{n,t\not=0}$, i.e. $w_{t\not=0}(\ell)\leq n$, then $w_{t=0}(\abs\ell)\leq n$ too and $\abs\ell\in S_{n, t=0}$. Moreover, if $(\ell, \setof{i})$ is an 1-cube in $L(C_{t\not=0})$, then it is sent to the 1-cube $[|\ell|,|\ell|+1]$ of $L(C_{t=0})$. By the above inequality at the level of vertices we also get that if $w_{t\not=0}((\ell, \setof{i}))\leq n$ then $w_{t=0}([|\ell|,|\ell|+1])\leq n $ as well. In particular, we get a well-defined map $H_0(S_{n, t\not=0},\bbz)\to H_0(S_{n,t=0},\bbz)$. This induces the wished maps at the level of graded roots and $\lcoh^0$. Note that the maps for $\lcoh^{>0}$ are trivial since $\lcoh^{>0}(C_{t=0},o)=0$ by the last paragraph of Example \thmref{lat-coh-restriction}.
 \end{pf}

\begin{rem}[label={rem:GN}]
    We formulated the above theorem for plane curves since the inequality $\hh_{t\not=0}(\ell)\geq \hh_{t=0}(\abs\ell)$ in \cite{GorNem2015} was stated and proved for plane curves. However, the proof from \cite{GorNem2015} naturally extends to arbitrary deformations (whenever $(C_{t=0},o)$ is irreducible), a fact which will be explained later in Remark \thmref{rem:ketfele}. For the proof of the above inequality see \eqref{eq:hdefdef}.
\end{rem}

One of our main goals is to generalize this theorem for the case when the curves are not necessarily plane curves, and the central curve is not necessarily irreducible. As we will see, the proof of the extended version is technically much harder. The main difficulty is to find the map connecting the lattices (and the tower of spaces $\setof*{S_n}_n$), the analogue of $\ell\mapsto\abs{\ell}$ from above. If the central curve is not irreducible, the existence of a similar map is not evident at all.

\vspace{2mm}

There is a central interest to analyse deformations along which certain invariants stay constant. In our case, we can ask for a characterization of those flat deformations when ${\Got R}$ or/and $\lcoh^*$ stay constant. If along a deformation the lattice cohomology $\lcoh^*$ is constant then by the Euler characteristic formula we obtain that $\delta$ is constant too. Then we can concentrate on the converse question as well: is it true that along a flat $\delta$-constant deformation the lattice cohomology (or/and the graded root) stay constant?

Let us analyse some examples.

\begin{ex}[label={ex:Teissier}]
    Let $(C,o)$ be an irreducible isolated plane curve singularity with semigroup $\semigp$. Then $\gr_\filt\fnring$ is isomorphic, as a graded algebra, with $\bbc[C^\semigp]$. (For some notation see \thmref{ex:semigrgen}.) Moreover, there exists a one-parameter complex analytic deformation of $(C^\semigp,0)$ such that the nearby curves $(C_{t\not=0},o)$ are equisingular with $(C,o)$ (hence share the same Hilbert function) \cite{TZ}. (For a concrete equation see \thmref{ex:Teissier2}.)

    Since along such deformation of irreducible germs $\semigp$ is stable, they are $\delta$-constant deformation. Then $\hh(\ell)=\#\setdef{s\in\semigp}{s<\ell}$ and $w_0$ are stable too. In particular, $\lcoh^*(C,o)=\lcoh^*(C^{\semigp},o)$.
\end{ex}

\begin{ex}[label={ex:Teissier2}]
    Consider the plane curve singularity $(C,o)$ given by $f=(x^3-y^2)^2-x^5y$. Its Newton pairs are $(2,3),(2,1)$. The conductor $c_{\semigp}=c$ is $c=\nu=2\delta=16$, and $\semigp=\gen{4,6,13}$. The semigroup $\semigp$ (marked by filled circles) and the $\hh$ and $w_0$ values for $0\leq\ell\leq 16$ are:
    \medskip
    \begin{center}
        \begin{tikzpicture}[scale=0.8]
            \draw[fill] (0,0) circle [radius=0.10];
            \draw (1,0) circle [radius=0.10];
            \draw (2,0) circle [radius=0.10];
            \draw (3,0) circle [radius=0.10];
            \draw[fill] (4,0) circle [radius=0.10];
            \draw (5,0) circle [radius=0.10];
            \draw[fill] (6,0) circle [radius=0.10];
            \draw (7,0) circle [radius=0.10];
            \draw[fill] (8,0) circle [radius=0.10];
            \draw (9,0) circle [radius=0.10];
            \draw[fill] (10,0) circle [radius=0.10];
            \draw (11,0) circle [radius=0.10];
            \draw[fill] (12,0) circle [radius=0.10];
            \draw[fill] (13,0) circle [radius=0.10];
            \draw[fill] (14,0) circle [radius=0.10];
            \draw (15,0) circle [radius=0.10];
            \draw[fill] (16,0) circle [radius=0.10];

            \draw (0,0) -- (16,0);

            \node at (-1,-0.8) {$\hh$};
            \node at (0,-0.8) {$0$};
            \node at (1,-0.8) {$1$};
            \node at (2,-0.8) {$1$};
            \node at (3,-0.8) {$1$};
            \node at (4,-0.8) {$1$};
            \node at (5,-0.8) {$2$};
            \node at (6,-0.8) {$2$};
            \node at (7,-0.8) {$3$};
            \node at (8,-0.8) {$3$};
            \node at (9,-0.8) {$4$};
            \node at (10,-0.8) {$4$};
            \node at (11,-0.8) {$5$};
            \node at (12,-0.8) {$5$};
            \node at (13,-0.8) {$6$};
            \node at (14,-0.8) {$7$};
            \node at (15,-0.8) {$8$};
            \node at (16,-0.8) {$8$};

            \node at (-1,-1.6) {$w_0$};
            \node at (0,-1.6) {$0$};
            \node at (1,-1.6) {$1$};
            \node at (2,-1.6) {$0$};
            \node at (3,-1.6) {$-1$};
            \node at (4,-1.6) {$-2$};
            \node at (5,-1.6) {$-1$};
            \node at (6,-1.6) {$-2$};
            \node at (7,-1.6) {$-1$};
            \node at (8,-1.6) {$-2$};
            \node at (9,-1.6) {$-1$};
            \node at (10,-1.6) {$-2$};
            \node at (11,-1.6) {$-1$};
            \node at (12,-1.6) {$-2$};
            \node at (13,-1.6) {$-1$};
            \node at (14,-1.6) {$0$};
            \node at (15,-1.6) {$1$};
            \node at (16,-1.6) {$0$};
        \end{tikzpicture}
    \end{center}
    \nni
    The graded root is
    \begin{center}
        \begin{tikzpicture}[scale=0.5]
            \draw[fill] (0,2) circle [radius=0.12];
            \draw[fill] (0,1) circle [radius=0.12];
            \draw[fill] (-1,0) circle [radius=0.12];
            \draw[fill] (0,0) circle [radius=0.12];
            \draw[fill] (1,0) circle [radius=0.12];
            \draw[fill] (0,-1) circle [radius=0.12];
            \draw[fill] (-2,-2) circle [radius=0.12];
            \draw[fill] (-1,-2) circle [radius=0.12];
            \draw[fill] (0,-2) circle [radius=0.12];
            \draw[fill] (1,-2) circle [radius=0.12];
            \draw[fill] (2,-2) circle [radius=0.12];

            \draw (0,1) -- (0,2.5);
            \draw[dotted] (0,2.5) -- (0,3.5);
            \draw (0,1) -- (-1,0);
            \draw (0,1) -- (1,0);
            \draw (0,1) -- (0,-2);
            \draw (0,-1) -- (-2,-2);
            \draw (0,-1) -- (-1,-2);
            \draw (0,-1) -- (1,-2);
            \draw (0,-1) -- (2,-2);

            \draw[dotted] (-2.5,1) -- (2.5,1);
            \draw[dotted] (-2.5,0) -- (2.5,0);
            \draw[dotted] (-2.5,-1) -- (2.5,-1);
            \draw[dotted] (-2.5,-2) -- (2.5,-2);

            \node[left] at (5.0,1) {$1$};
            \node[left] at (5.0,0) {$0$};
            \node[left] at (5.0,-1) {$-1$};
            \node[left] at (5.0,-2) {$-2$};
        \end{tikzpicture}
    \end{center}
    By \thmref{ex:Teissier} this is the root of $(C^\semigp,o)$ too. The deformation of $(C^\semigp,o)$ is $x^3=y^2+tz, z^2=x^5y$, where $t$ is the deformation parameter. Hence, by the above discussion, the (rank of the) lattice is constant, the semigroup is constant, hence ${\Got R}$ and $\lcoh^*$ remain constant.

    \vspace{2mm}

    Consider now another deformation of $f$, namely $f_t=(x^3+y^2)(x^3+ty^2)+x^5y$ with $\abs{t-1}$ small. Then $(C_{t\neq 1},o)$ is Newton non-degenerate, with two irreducible components, $\mu=15$, $\delta=8$. In this case $c_{\semigp}=(8,8)$ and the $\hh$ and $w_0$ tables are the following:

    \begin{center}
        \setlength{\tabcolsep}{0.6em}
        \def\w#1 {\makebox[0pt][r]{$#1$}}
        \begin{minipage}[b]{0.45\textwidth}
            \[\hh\]
            \begin{center}
                \begin{tabular}{R|RRRRRRRRR}
                    8   & 7 & 7 & 7 & 7 & 7 & 7 & 7 & 8 & 8 \\
                    7   & 6 & 6 & 6 & 6 & 6 & 6 & 6 & 7 & 8 \\
                    6   & 5 & 5 & 5 & 5 & 5 & 5 & 5 & 6 & 7 \\
                    5   & 4 & 4 & 4 & 4 & 4 & 4 & 5 & 6 & 7 \\
                    4   & 3 & 3 & 3 & 3 & 3 & 4 & 5 & 6 & 7 \\
                    3   & 2 & 2 & 2 & 2 & 3 & 4 & 5 & 6 & 7 \\
                    2   & 1 & 1 & 1 & 2 & 3 & 4 & 5 & 6 & 7 \\
                    1   & 1 & 1 & 1 & 2 & 3 & 4 & 5 & 6 & 7 \\
                    0   & 0 & 1 & 1 & 2 & 3 & 4 & 5 & 6 & 7 \\
                    \hline
                        & 0 & 1 & 2 & 3 & 4 & 5 & 6 & 7 & 8 \\
                \end{tabular}
            \end{center}
        \end{minipage}
        \begin{minipage}[b]{0.45\textwidth}
            \[w_0\]
            \begin{center}
                \begin{tabular}{R|RRRRRRRRR}
                    8   & 6 & 5 & 4 & 3 & 2 & 1 & 0 & 1 & 0 \\
                    7   & 5 & 4 & 3 & 2 & 1 & 0 &\w-1   & 0 & 1 \\
                    6   & 4 & 3 & 2 & 1 & 0 &\w-1   &\w-2   &\w-1   & 0 \\
                    5   & 3 & 2 & 1 & 0 &\w-1   &\w-2   &\w-1   & 0 & 1 \\
                    4   & 2 & 1 & 0 &\w-1   &\w-2   &\w-1   & 0 & 1 & 2 \\
                    3   & 1 & 0 &\w-1   &\w-2   &\w-1   & 0 & 1 & 2 & 3 \\
                    2   & 0 &\w-1   &\w-2   &\w-1   & 0 & 1 & 2 & 3 & 4 \\
                    1   & 1 & 0 &\w-1   & 0 & 1 & 2 & 3 & 4 & 5 \\
                    0   & 0 & 1 & 0 & 1 & 2 & 3 & 4 & 5 & 6 \\
                    \hline
                        & 0 & 1 & 2 & 3 & 4 & 5 & 6 & 7 & 8 \\
                \end{tabular}
            \end{center}
        \end{minipage}
    \end{center}
    \bigskip
    \nni
    Thus, the graded root and $\lcoh^*$ of $f_{t}$ stay stable along the deformation.

    In this case, $\phi\colon L(C_{t\not=0})=\bbz^2\to L(C_{t=0})=\bbz$, $\ell\mapsto\phi(\ell)=\abs\ell$, has the following property: $\min w_{t\not=1}|_{\phi^{-1}(|\ell|)}=w_{t=1}(\abs\ell)$. And, evidently, this $\phi$ is a linear map.

    On the other hand, the `linear' (rank one) $\hh$ and $w_0$-tables of $f_1$ cannot be embedded `linearly' into the $r=2$-tables of $f_{t\neq 1}$. However, a non-linear weighted embedding exists into the path with entries $(i,i)$, $0\leq i\leq 8$ and $(i,i+1)$, $0\leq i\leq 7$ (and the rank two $w$-table contracts to this path).

    In fact, these two maps connecting the two $w$-tables induce the isomorphisms between the corresponding roots and lattice cohomologies.

    Note also that there exists another choice for such embedding provided by that entries $(i,i)$, $0\leq i\leq 8$ and $(i,i-1)$, $1\leq i\leq 8$, hence for the `combinatorial inverse' of $\phi$ at the level of lattices there is no canonical choice.
\end{ex}

\begin{ex}[label={ex:delta2}]
    Consider the plane curve singularities
    \[
        E_6:\ \setof{x^3+y^4=0}\ \text{ and }\ A_5:\ \setof{(x+y^3)(x-y^3)=0}.
    \]
    For both of these we have $\delta=3$ (with $\mu=6$ and $5$ respectively), and there is a deformation $f_t=x^3+(y^2+tx)^2$ between them (cf.~\cite{ArnoldCurves}).

    In the $t=0$ case, $E_6$, the semigroup $\semigp=\gen{3,4}$ with conductor $c_{\semigp}=6$ gives us
    \bigskip
    \begin{center}
        \setlength{\tabcolsep}{0.6em}
        $\vcenter
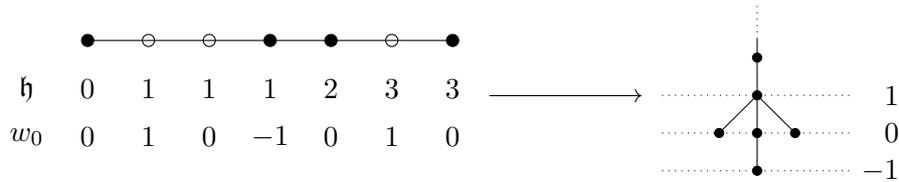
{\hbox{
            \begin{tikzpicture}[scale=0.8]
                \draw[fill] (0,0) circle [radius=0.10];
                \draw (1,0) circle [radius=0.10];
                \draw (2,0) circle [radius=0.10];
                \draw[fill] (3,0) circle [radius=0.10];
                \draw[fill] (4,0) circle [radius=0.10];
                \draw (5,0) circle [radius=0.10];
                \draw[fill] (6,0) circle [radius=0.10];

                \draw (0,0) -- (6,0);

                \node at (-1,-0.8) {$\hh$};
                \node at (0,-0.8) {$0$};
                \node at (1,-0.8) {$1$};
                \node at (2,-0.8) {$1$};
                \node at (3,-0.8) {$1$};
                \node at (4,-0.8) {$2$};
                \node at (5,-0.8) {$3$};
                \node at (6,-0.8) {$3$};

                \node at (-1,-1.6) {$w_0$};
                \node at (0,-1.6) {$0$};
                \node at (1,-1.6) {$1$};
                \node at (2,-1.6) {$0$};
                \node at (3,-1.6) {$-1$};
                \node at (4,-1.6) {$0$};
                \node at (5,-1.6) {$1$};
                \node at (6,-1.6) {$0$};
            \end{tikzpicture}
        }}\vcenter{\hbox{
            \begin{tikzpicture}[scale=1.0]
                \draw[->] (0,0) -- (2,0);
            \end{tikzpicture}
        }}\vcenter{\hbox{
            \begin{tikzpicture}[scale=0.5]
                \draw[fill] (0,2) circle [radius=0.12];
                \draw[fill] (0,1) circle [radius=0.12];
                \draw[fill] (-1,0) circle [radius=0.12];
                \draw[fill] (0,0) circle [radius=0.12];
                \draw[fill] (1,0) circle [radius=0.12];
                \draw[fill] (0,-1) circle [radius=0.12];

                \draw (0,1) -- (0,2.5);
                \draw[dotted] (0,2.5) -- (0,3.5);
                \draw (0,1) -- (0,-1);
                \draw (0,1) -- (-1,0);
                \draw (0,1) -- (1,0);

                \draw[dotted] (-2.5,1) -- (2.5,1);
                \draw[dotted] (-2.5,0) -- (2.5,0);
                \draw[dotted] (-2.5,-1) -- (2.5,-1);

                \node[left] at (4.0,1) {$1$};
                \node[left] at (4.0,0) {$0$};
                \node[left] at (4.0,-1) {$-1$};
            \end{tikzpicture}
        }}$
    \end{center}
    \bigskip
    For the $A_5$-singularity, we get $c_{\semigp}=(3,3)$. The Hilbert function $\hh$ can be computed via the formula from \thmref{ss:planecurves} based an embedded resolution and Poincaré series (cf. \cite{GorNem2015}). Thus we obtain
    \bigskip
    \begin{center}
        $\vcenter{\hbox{
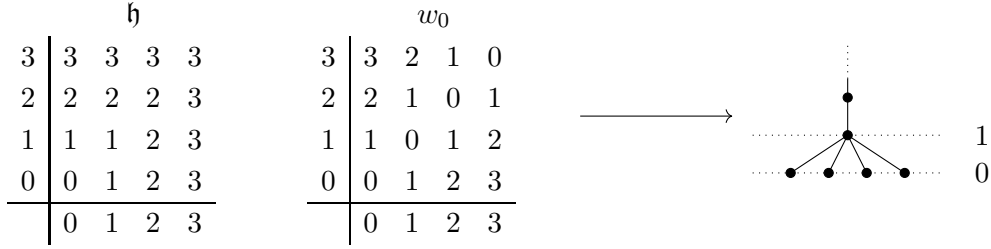
            \begin{minipage}[b]{0.2\textwidth}
                \[\hh\]
                \smallskip
                \begin{tabular}{R|RRRR}
                    3   & 3 & 3 & 3 & 3 \\
                    2   & 2 & 2 & 2 & 3 \\
                    1   & 1 & 1 & 2 & 3 \\
                    0   & 0 & 1 & 2 & 3 \\
                    \hline
                        & 0 & 1 & 2 & 3 \\
                \end{tabular}
            \end{minipage}
            \hskip0.5cm
            \begin{minipage}[b]{0.2\textwidth}
                \[w_0\]
                \smallskip
                \begin{tabular}{R|RRRR}
                    3   & 3 & 2 & 1 & 0 \\
                    2   & 2 & 1 & 0 & 1 \\
                    1   & 1 & 0 & 1 & 2 \\
                    0   & 0 & 1 & 2 & 3 \\
                    \hline
                        & 0 & 1 & 2 & 3 \\
                \end{tabular}
            \end{minipage}
        }}\vcenter{\hbox{
            \begin{tikzpicture}[scale=1.0]
                \draw[->] (0,0) -- (2,0);
            \end{tikzpicture}
        }}\vcenter{\hbox{
            \begin{tikzpicture}[scale=0.5]
                \draw[fill] (0,2) circle [radius=0.12];
                \draw[fill] (0,1) circle [radius=0.12];
                \draw[fill] (-1.5,0) circle [radius=0.12];
                \draw[fill] (-0.5,0) circle [radius=0.12];
                \draw[fill] (0.5,0) circle [radius=0.12];
                \draw[fill] (1.5,0) circle [radius=0.12];

                \draw (0,1) -- (0,2.5);
                \draw[dotted] (0,2.5) -- (0,3.5);
                \draw (0,1) -- (-1.5,0);
                \draw (0,1) -- (-0.5,0);
                \draw (0,1) -- (0.5,0);
                \draw (0,1) -- (1.5,0);

                \draw[dotted] (-2.5,1) -- (2.5,1);
                \draw[dotted] (-2.5,0) -- (2.5,0);

                \node[left] at (4.0,1) {$1$};
                \node[left] at (4.0,0) {$0$};
            \end{tikzpicture}
        }}$
    \end{center}
    \bigskip
    Hence, in this case, though this is a $\delta$-constant deformation, $\Got R$ does not stay constant.

    Since $E_6=(C_{t=0},o)$ is irreducible, we still have $\phi\colon L(C_{t\not=0})=\bbz^2\to L(C_{t=0})=\bbz$, $\ell\mapsto \phi(\ell)=\abs\ell$, hence a graded (graph) map of degree zero at the level of graded roots ${\Got R}(C_{t\not=0},o)\to{\Got R}(C_{t=0},o)$. However, at this time $\min w_{t\not=1}|_{\phi^{-1}(\abs\ell)}\neq w_{t=1}(\abs\ell)$ for $\abs\ell=3$. And an `inverse' weighed embedding is also obstructed.
\end{ex}

\begin{ex}[label={ex:delta3}]
    A slightly more complicated $\delta$-constant deformation when the graded root and the lattice cohomologies do not stay constant is the deformation given by $f_t=x^4+(y^3+tx)^3$ between
    \[
        W_{24}\,:\, \setof{x^4+y^9=0}\text{ and } \ J_{4,0}\,:\, \setof{x^3+y^{12}=0}.
    \]
    Here we have $\delta=12$ for both germs, and in the case of the
    irreducible one it is easy to compute the root from the semigroup $\semigp=\gen{4,9}$ (with $\Got c_{\semigp}=24$). The other one splits into $3$ components, and as such we end up with having to compute $\hh$ and $w_0$ over the cube $[0,c_{\semigp}]=[(0,0,0),(8,8,8)]$. Clearly, with more components, the computation quickly becomes difficult (and even harder to illustrate). In this case, a computation shows that the resulting graded roots are the following:
    \begin{center}
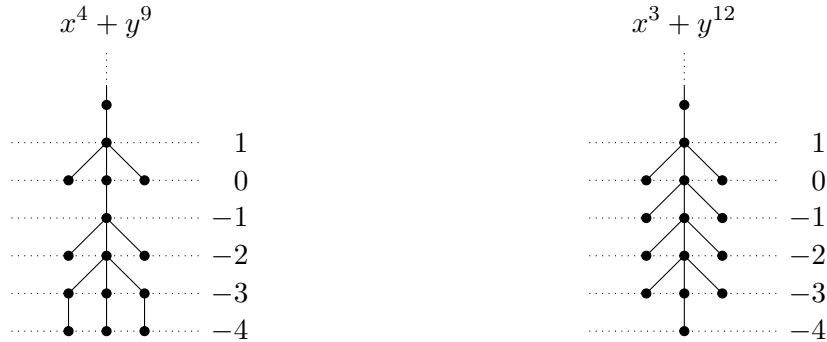
        \begin{minipage}[b]{0.45\textwidth}
            \[x^4+y^9\]
            \begin{center}
                \begin{tikzpicture}[scale=0.5]
                    \draw[fill] (0,2) circle [radius=0.12];
                    \draw[fill] (0,1) circle [radius=0.12];
                    \draw[fill] (-1,0) circle [radius=0.12];
                    \draw[fill] (0,0) circle [radius=0.12];
                    \draw[fill] (1,0) circle [radius=0.12];
                    \draw[fill] (0,-1) circle [radius=0.12];
                    \draw[fill] (-1,-2) circle [radius=0.12];
                    \draw[fill] (0,-2) circle [radius=0.12];
                    \draw[fill] (1,-2) circle [radius=0.12];
                    \draw[fill] (-1,-3) circle [radius=0.12];
                    \draw[fill] (0,-3) circle [radius=0.12];
                    \draw[fill] (1,-3) circle [radius=0.12];
                    \draw[fill] (-1,-4) circle [radius=0.12];
                    \draw[fill] (0,-4) circle [radius=0.12];
                    \draw[fill] (1,-4) circle [radius=0.12];

                    \draw (0,1) -- (0,2.5);
                    \draw[dotted] (0,2.5) -- (0,3.5);
                    \draw (0,1) -- (0,-4);
                    \draw (0,1) -- (-1,0);
                    \draw (0,1) -- (1,0);
                    \draw (0,-1) -- (-1,-2);
                    \draw (0,-1) -- (1,-2);
                    \draw (0,-2) -- (-1,-3) -- (-1,-4);
                    \draw (0,-2) -- (1,-3) -- (1,-4);

                    \draw[dotted] (-2.5,1) -- (2.5,1);
                    \draw[dotted] (-2.5,0) -- (2.5,0);
                    \draw[dotted] (-2.5,-1) -- (2.5,-1);
                    \draw[dotted] (-2.5,-2) -- (2.5,-2);
                    \draw[dotted] (-2.5,-3) -- (2.5,-3);
                    \draw[dotted] (-2.5,-4) -- (2.5,-4);

                    \node[left] at (4.0,1) {$1$};
                    \node[left] at (4.0,0) {$0$};
                    \node[left] at (4.0,-1) {$-1$};
                    \node[left] at (4.0,-2) {$-2$};
                    \node[left] at (4.0,-3) {$-3$};
                    \node[left] at (4.0,-4) {$-4$};

                    \node[right] at (-4.0,0) {};
                \end{tikzpicture}
            \end{center}
        \end{minipage}
        \begin{minipage}[b]{0.45\textwidth}
            \[x^3+y^{12}\]
            \begin{center}
                \begin{tikzpicture}[scale=0.5]
                    \draw[fill] (0,2) circle [radius=0.12];
                    \draw[fill] (0,1) circle [radius=0.12];
                    \draw[fill] (-1,0) circle [radius=0.12];
                    \draw[fill] (0,0) circle [radius=0.12];
                    \draw[fill] (1,0) circle [radius=0.12];
                    \draw[fill] (-1,-1) circle [radius=0.12];
                    \draw[fill] (0,-1) circle [radius=0.12];
                    \draw[fill] (1,-1) circle [radius=0.12];
                    \draw[fill] (-1,-2) circle [radius=0.12];
                    \draw[fill] (0,-2) circle [radius=0.12];
                    \draw[fill] (1,-2) circle [radius=0.12];
                    \draw[fill] (-1,-3) circle [radius=0.12];
                    \draw[fill] (0,-3) circle [radius=0.12];
                    \draw[fill] (1,-3) circle [radius=0.12];
                    \draw[fill] (0,-4) circle [radius=0.12];

                    \draw (0,1) -- (0,2.5);
                    \draw[dotted] (0,2.5) -- (0,3.5);
                    \draw (0,1) -- (0,-4);
                    \draw (0,1) -- (-1,0);
                    \draw (0,1) -- (1,0);
                    \draw (0,0) -- (-1,-1);
                    \draw (0,0) -- (1,-1);
                    \draw (0,-1) -- (-1,-2);
                    \draw (0,-1) -- (1,-2);
                    \draw (0,-2) -- (-1,-3);
                    \draw (0,-2) -- (1,-3);

                    \draw[dotted] (-2.5,1) -- (2.5,1);
                    \draw[dotted] (-2.5,0) -- (2.5,0);
                    \draw[dotted] (-2.5,-1) -- (2.5,-1);
                    \draw[dotted] (-2.5,-2) -- (2.5,-2);
                    \draw[dotted] (-2.5,-3) -- (2.5,-3);
                    \draw[dotted] (-2.5,-4) -- (2.5,-4);

                    \node[left] at (4.0,1) {$1$};
                    \node[left] at (4.0,0) {$0$};
                    \node[left] at (4.0,-1) {$-1$};
                    \node[left] at (4.0,-2) {$-2$};
                    \node[left] at (4.0,-3) {$-3$};
                    \node[left] at (4.0,-4) {$-4$};

                    \node[right] at (-4.0,0) {};
                \end{tikzpicture}
            \end{center}
        \end{minipage}
    \end{center}
\end{ex}

\section{Functoriality under flat deformations. Introductory words.}\label{sec:functorintr}

\subsection{Preliminary words}\label{ss:prel5}
Above, starting from an isolated curve singularity, we defined a CW complex $\cw$ together with a set of compatible weight functions $w$, namely $(\cw, w)=(R(0,c), w)$. In thus procedure we used the Hilbert function $\hh(\ell)$ of $(C,o)$ and the induced weight function $w_0(\ell)=2\hh(\ell)-|\ell|$.

Next, we wish to associate with a flat deformation $(C_t,o)_{t\in(\bbc,0)}$ a (cellular) map ${\mathfrak d} :(\cw,w)_{C_{t\not=0}}\to (\cw, w)_{C_{t=0}}$, compatible with the weights, as in Theorem \thmref{th:DEF}.

The first try to construct ${\mathfrak d}$ is obstructed by the fact that there is no natural geometric principle which would guide what to associate with a lattice point $\ell\in R(0,c)_{C_{t\not=0}}$. The problem is not only the fact that the number of irreducible components (i.e. the rank of the lattices) jumps (in any direction), but also that usually there are `too many irrelevant' lattice points from the point of view of the `essential geometrical properties' of the curves (and the abundance of the lattice points can be rather different for $t=0$ and $t\not=0$). Interestingly enough, the object which plays a key role in the definition of the lattice cohomology, namely the set of lattice points, should be decreased, partly eliminated.

Let us focus on another object, which has a more controlled behaviour under deformations. Let us consider the linear subspaces $\setof{\filt(\ell)}_{\ell}$ of $\fnring_{(C,o)}$.

Assume that the flat deformation $(C_t,o)_t$ is realized in some $(\bbc^N,o)$, that is $(C_t,o)\subset(\bbc^N,o)$ for any $t$. Using this, the comprising spaces for the linear subspaces (when $t$ varies) can be identified. Indeed, for any fixed $t$, if $q\colon\fnring_{(\bbc^N,o)}\to\fnring_{(C_t,o)}$ is the natural projection let us set $\widetilde{\filt}_t(\ell):=q^{-1}(\filt_t(\ell))$. Then $\codim( \widetilde{\filt}_t(\ell)\subset\fnring_{(\bbc^N,o)})=\codim (\filt_t(\ell)\subset\fnring_{(C_t,o)})=\hh_t(\ell)$. Hence, for any $t$, the system of subspaces $\setof*{\widetilde{\filt}_t(\ell)}_{\ell}$ preserves all the information about the inclusion structure of $\setof{\filt(\ell)}_\ell$ and about the Hilbert coefficients $\hh_t(\ell)$, and additionally they sit in (the $t$-independent) $\fnring_{(\bbc^N,o)}$. Furthermore, $\lim_{0\not=t\to 0}\widetilde{\filt}_t(\ell)$ in some Grassmannian of $\fnring_{(\bbc^N,o)}$ exists (cf. \ref{ss:defmap}). Usually it is not of the form $\widetilde{\filt}_{t=0}(\ell')$, but (as we will see) it can be compared with such an object. Therefore, the linear subspace arrangement $\setof{\filt(\ell)}_{\ell}$ (or $\setof*{\widetilde{\filt}(\ell)}_{\ell}$) has several good stability properties under deformations, which can be (and will be) exploited. Though the numerical invariant $\hh(\ell)$ can be recovered from the linear subspace arrangement as the corresponding codimensions, $w_0(\ell)=2\hh(\ell)-\abs\ell$ usually not; hence we need to keep these weights as well. This shows that we need to analyse the structure of $\setof*{(\filt(\ell), w_0(\ell))}_{\ell}$, try to recover the lattice cohomology from it, and determine its behaviour under deformations. This in particular will provide a new definition of the lattice cohomology, which helps to prove the functoriality under deformations.

Before we continue our overview of the new constructions (cf. \ref{ss:prel5b}) let us list some properties of the subspaces $\setof{(\filt(\ell), w_0(\ell))}_{\ell} $ and their compatibility with the cubical decomposition.

Let us start with the following simple observation:

\begin{lem}[label=uniquemax_0]
    For any singularity $(C,o)$ and for any $\ell\in\bbz_{\geq 0}^r$, there is a unique maximal $\ell^*\in \bbz_{\geq 0}^r$ such that $\filt(\ell^*)=\filt(\ell)$. In particular, for any $\ell\leq\ell'\leq\ell^*$, we have $\filt(\ell)=\filt(\ell')=\filt(\ell^*)$.
\end{lem}
\begin{pf}
    Use the property $\filt(\max\setof{\ell_1,\ell_2})=\filt(\ell_1)\cap\filt(\ell_2)$.
\end{pf}

Next, for any fixed subspace $V\subset\fnring_{(C,o)}$, let us consider the following collection of cubes.
\[
    \levelcells^q(V) = \setdef{C\in\cellset_q}{\forall\ell\in C\cap \bbz^r:\,\filt(\ell)=V}, \quad
    \levelcells^*(V) = \bigcup_q\levelcells^q(V), \quad
    \levelset(V) = \bigcup_{C\in\levelcells^*(V)}C.
\]
Furthermore, for any $n\in\bbz$, let
\[
    \levelcells^*_n(V) = \bigcup_q\setdef{C\in\cellset_q}{\forall\ell\in C\cap \bbz^r:\,\filt(\ell)=V, \ w(C)\leq n}, \quad \levelset_n(V) =\bigcup_{C\in\levelcells^*_n(V)}C.
\]
By Lemma \thmref{uniquemax_0}, there is a unique maximal lattice point of $\levelset(V)$: this shall be denoted by $\maxpoint_V$. Moreover, since $\filt\bigl(\max\setof{\ell_1,\ell_2}\bigr)=\filt(\ell_1)\cap\filt(\ell_2)$, for any $W\subset V$ we have $\maxpoint_{W}\geq\maxpoint_V$.

\begin{lem}[label={ss-set-props}]
    For any $V$ all non-empty sets $\levelset_n(V)$ are contractible, and contain the unique common point $\maxpoint_V$. In fact, the inclusion $\levelset_{w(\maxpoint_V)}(V)=\setof{\maxpoint_V}\subset\levelset_n(V)$ admits a strong deformation retraction for any $n$.

    The lattice points $\setof{\maxpoint_V}_V$ are exactly the elements of the semigroup $\semigp$ of $(C,o)$.
\end{lem}
\begin{pf}
    Since $w_0(\ell)=2\hh(\ell)-\abs{\ell}$ is strictly decreasing on $\levelcells^0(V)$, it follows from \thmref{uniquemax_0} that $\maxpoint_V$ has the strictly smallest weight in any $S_n(V)$. Regarding the contractibility, one sees that if $\ell\in S_n(V)$ then $\ell\leq \maxpoint _V$ and the real rectangle $R(\ell,\maxpoint_V)$ belongs to $S_n(V)$.

    The last statement follows from \eqref{eq:Sfromh}.
\end{pf}

This means that within any level set $\levelset_n=\cup_{w(\topcell)\leq n}\topcell$ --- which is the union of various sets $\levelset_n(V)$, plus cells connecting them --- we can always contract these subsets of type $\levelset_n(V)$ onto points without changing the singular (co)homology of $S_n$. (The precise construction will be presented below.)

This serves as the motivation to construct a weighted CW complex with vertices the previously defined $\maxpoint_V$'s, which will give us the same lattice cohomology as the original one. Hence, the lattice points $\bbz^r$ will be replaced by a restricted class of special lattice points, the semigroup elements. On the other hand, the structure of the higher dimensional cubes/cells is more complicated.

The following elementary lemma shows how can we recover in the language of the linear subspace system $\spsystem :=\setof{\filt(\ell)}_{\ell}$ a cube of type $(\maxpoint_V,I)$. This will then help to define the cells in the new situation.

\begin{lem}[label={ss-set-props2}]
    \begin{enumerate}[label={(\alph*)}]
        \item If $V=\filt(\ell_1)$, $W=\filt(\ell_2)$, and $W\subset V$ then there exists $\ell_3\geq \ell_1$ such that $W=\filt(\ell_3)$.
        
        \item If $V,W\in \spsystem$ with $W\subset V$, then there exists a sequence $W=V_0\subset V_1\subset \ldots \subset V_s=V$, $V_i\in \spsystem$, such that ${\rm codim} (V_i\subset V_{i+1})=1$.
        
        \item For any $V\in\spsystem$ the $\maxpoint_V$-adjacent vertices of the $r$-dimensional cube $(\maxpoint_V, \grverts)$ of $R(0,c)=\cw$ correspond bijectively with $\setdef{W\in\spsystem}{\codim(W\subset V)=1}$. Cubes of type $(\maxpoint_V,I)$ $(I\subset \grverts)$ correspond with subsets of $\setdef{W\in\spsystem}{\codim(W\subset V)=1}$ of cardinality $\abs{I}$.
    \end{enumerate}
\end{lem}

\subsection{The strategy.}\label{ss:prel5b}

The following diagram shows the first main steps of our next construction.
\[
    \begin{array}{ccc}
        (C,o) & \rightsquigarrow & (R(0,c),w)=(\cw,w) \\
        \downarrow & & \downarrow\vcenter{\rlap{$\alpha$}} \\
        (\spsystem,w) & \rightsquigarrow & (\tilde{\cw},\tilde{w})
    \end{array}
\]
The first row indicates the original procedure: the construction of the weighted CW complex
$(R(0,c),w)$ from the curve $(C,o)$ (cf. section \ref{sec:curveslattice}). The next steps are the following.

\begin{enumerate}[label={(\arabic*)}]
    \item We define the abstract notion of `weighted system of subspaces';
    
    \item by a general construction we associate to any abstract weighted system of subspaces $(\spsystem,w)$ a weighted CW complex $(\tilde{\cw},\tilde{w})$;
    
    \item we show that $\setof{(\filt(\ell), w_0(\ell)}_{\ell}$ associates a weighted system of subspaces to the curve $(C,o)$;
    
    \item then, for a fixed curve singularity $(C,o)$ we compare the output of the composition $(C,o)\rightsquigarrow\setof{(\filt(\ell), w_0(\ell)}_{\ell}= (\spsystem,w)\rightsquigarrow (\tilde{\cw},\tilde{w})$ with $(R(0,c),w)=(\cw,w)$. We show that there exists a surjective continuous map $\alpha: (\cw,w)\to (\tilde{\cw},\tilde{w})$ with all fibers contractible, and which is a homotopy equivalence (compatible with the weights). In fact, this map $\alpha$ is the universal cellular map such that $\alpha (\ell_1)=\alpha(\ell_2)$ if and only if $\filt (\ell_1)=\filt(\ell_2)$. The main point is that $\alpha$ induces isomorphisms $\lcoh^*(\cw,w)=\lcoh^*(\tilde{\cw},\tilde{w})$ and ${\Got R}(\cw,w)={\Got R}(\tilde{\cw},\tilde{w})$.
\end{enumerate}

\vspace{2mm}

Note that the above steps provide a different construction of $\lcoh^*(C,o)$ (compared with the original one via $(R(0,c),w)$). This new version uses the linear subspaces $\setof{\filt(\ell)}_{\ell}$ with certain weights. We hope that this can serve as a model to extend the definition of the lattice cohomology to several other geometric situations when linear subspace arrangements appear.

\pagebreak
\section{Functoriality under flat deformations. The construction}\label{sec:functor}

\subsection{Weighted system of subspaces and their CW complexes}\label{ss:linsubspacescw}\hfill\smallbreak

\begin{defn}[label={V-cell}]
    Fix a (usually infinite dimensional) vector space $\contspace$. We say that a set of subspaces $\spsystem$ of $\contspace $ form a `system of subspaces' if $\contspace\in\spsystem$, $\spsystem$ is closed under intersection (i.e., if $V_1, V_2\in \spsystem$ then $V_1\cap V_2 \in \spsystem $), and $\codim(V\subset\contspace)<\infty$ for any $V\in\contspace$.
\end{defn}

Our first goal is to define a CW complex associated with any such $\spsystem$. This will be done in three steps. First, we define and analyse the set of `combinatorial cells', then we associate with any combinatorial cell a `topological cell' (in fact, a cube), and finally we glue the topological cells into a CW complex. The next abstract definition is motivated by the statement of Lemma \thmref{ss-set-props2}.

\begin{defn}[label={V-comb-cell}]
    Given a vector space $\contspace$ and a system of subspaces $\spsystem$ of it, we call a subset $\combcell\subseteq\spsystem$ a combinatorial $\spsystem$-cell if it is of the form
    \[
        \combcell=\setdef*[\Big]{\bigcap_{\cellgen\in\cellgenss}\cellgen}{\emptyset\neq\cellgenss\subseteq\cellgens}\cup\setof{V}
    \]
    for some $V\in\spsystem$ and a finite set $\cellgens\subseteq\spsystem$ with $\codim(\cellgen\subseteq V)=1$ for all $\cellgen\in\cellgens$.

    We define the \It{dimension} of the combinatorial $\spsystem$-cell as $\dim\combcell=\abs{\cellgens}$. We call the spaces $V$ and $\cellgen\in\cellgens$ the \It{root} and the \It{generators} respectively, while $\cap_{\cellgen\in\cellgens}\cellgen$ (or $V$ if $\cellgens=\emptyset$) will be referred to as the \It{tip} of $\combcell$. The root, tip and the set of generators are respectively denoted as $\cellroot(\combcell)$, $\celltip(\combcell)$, and $\cellgens(\combcell)$, while the combinatorial $\spsystem$-cell obtained from a root $\cellroot$ and set of generators $\cellgens$ as $\combcell(\cellroot,\cellgens)$.

    We further define the \It{height} of the $\spsystem$-cell $\combcell$ as $\height\combcell=\codim\bigl(\celltip(\combcell)\subseteq\cellroot(\combcell)\bigr)$.
\end{defn}

\begin{rem}[label={V-cell-welldefined}]
    \begin{enumerate}[label={(\alph*)}]
        \item Note that the root, tip, and generators are uniquely determined by the set $\combcell$, so the above definition \thmref{V-comb-cell} makes sense. Indeed,
        \[
            \cellroot(\combcell)=\textstyle{\sum}_{V\in \combcell}V,\quad
            \celltip(\combcell)=\cap_{V\in \combcell}V,\quad
            \cellgens(\combcell)=\setdef*[\big]{W\in\combcell}{\codim\bigl(W\subseteq\cellroot(\combcell)\bigr)=1}.
        \]
        The last identity holds since all intersections of some distinct $\cellgen,\cellgen'\in\cellgens$ have codimension $\geq 2$.

        \item We always have $\height\combcell\leq\dim\combcell.$ The inequality can be strict e.g. in the case when all the generators are distinct subspaces of some $V=\cellroot(\combcell)$ containing a fixed codimension-$2$ subspace $W=\celltip(\combcell)$, and their number is $\dim\combcell\geq 3$.
    \end{enumerate}
\end{rem}

\begin{defn}[label={V-cell-faces-def}]
    For a given combinatorial $\spsystem$-cell $\combcell$, its \It{faces} are the elements of
    \[
        \setdef*[\big]{\combcellface\subseteq\combcell}{\combcellface\text{ is a combinatorial $\spsystem$-cell}}.
    \]
\end{defn}

Next, we start the recursive construction a CW complex $\tilde{\cw}$ associated with $\spsystem$.

\begin{ntn}[label={top-V-cell}]
    For any combinatorial $\spsystem$-cell $\combcell$, let
    \[
        \topcell^\circ (\combcell)=\setof*[\big]{\cellroot(\combcell)}\times (0,1)^{\cellgens(\combcell)}, \quad\text{and}\quad
        \topcell(\combcell)=\setof*[\big]{\cellroot(\combcell)}\times [0,1]^{\cellgens(\combcell)}.
    \]
    We will be gluing the boundaries of the cubes $\topcell(\combcell)$ to the lower dimensional skeleta when constructing the CW complex. Hence, the sets $\topcell^\circ (\combcell)$ will serve as the open cells of $\ti\cw$. We call $\topcell^\circ (\combcell)$ the \It{topological open $\spsystem$-cell} associated to the combinatorial cell $\combcell$. We may view this as just the open cube $(0,1)^{\cellgens(\combcell)}$, or in the case of $\dim\combcell=0$, this vertex as simply representing $\cellroot(\combcell)\in\spsystem$. However, in the above notation of $\topcell^\circ (\combcell)$ we wished to insert the root of the cube as well.

    The closed cubes $\topcell(\combcell)$ will be called \It{topological closed $\spsystem$-cells}.

    The exact method of gluing will be defined shortly, but the general principle is that in $[0,1]^{\cellgens(\combcell)}$, each vertex $u\colon\cellgens(\combcell)\to\setof{0,1}$ will be glued to the $0$-cell representing $\cap\,\setdef*[\big]{\cellgen}{G\in\cellgens(\combcell),\ u(\cellgen)=1}$ for $u\neq 0$, and $\cellroot(\combcell)$ for $u=0$.
\end{ntn}

For the $0$-dimensional skeleton, set
\[
    \sk_0\ti\cw=\spsystem\times\setof{\emptyset},
\]
i.e.~the vertices $\sk_0\ti\cw$ of $\ti\cw$ are essentially the elements of $\spsystem$.

Assume that we have already defined the $(q-1)$-dimensional skeleton $\sk_{q-1}\ti\cw$. For any $\topcell(\combcell)$ of dimension $<q$ let $\kappa_{\combcell}\colon\topcell(\combcell)\to\sk_{q-1}\ti\cw$ be the characteristic map of $\topcell(\combcell)$ into $\sk_{q-1}\ti\cw$ with attaching map $\kappa_{\combcell}|_{\D\topcell(\combcell)}$.

Next, we wish to define the attaching map of a fixed topological $\spsystem$-cell $\combcell$ of dimension $q$ (i.e. with $|\cellgens(\combcell)|=q$). We need to give a continuous map
\[
    \partial\kappa_{\combcell}\colon\D\topcell(\combcell)=\setof*[\big]{\cellroot(\combcell)}\times\D[0,1]^{\cellgens(\combcell)}\to\sk_{q-1}\ti\cw.
\]
Let us fix $u\in\D[0,1]^{\cellgens(\combcell)}=\setdef*[\big]{u\colon\cellgens(\combcell)\to [0,1]}{\exists\,\cellgen\in\cellgens(\combcell): u(\cellgen)\in\setof{0,1}}$, and set
\[
    \cellgens_0(u)=u^{-1}(0), \ns
    \cellgens_1(u)=u^{-1}(1), \ns
    \cellgens_*(u)=u^{-1}\bigl((0,1)\bigr).
\]
Define also $\cellroot_u=\cap_{W\in \cellgens_1(u)\cup\setof*{\cellroot(\combcell)}}\, W$. (If $\cellgens_1(u)\not=\emptyset$ then $R_u=\cap_{W\in \cellgens_1(u)}\, W$; otherwise $R_u=\cellroot(\combcell)$.) Finally, set
\[
    \ti\cellgens_*(u)=\setdef{\cellroot_u\cap\cellgen}{\cellgen\in\cellgens_*(u)}\setminus\setof{\cellroot_u} \subset \spsystem.
\]
Since $\cellgens_0(u)\cup \cellgens_1(u)\not=\emptyset$, $\abs{\cellgens_*(u)}<q$. Furthermore, the map $G\mapsto R_u\cap G$ ($G\in\cellgens_*(u)$) is usually not injective, and the intersection $R_u\cap G$ can also be equal to $R_u$. In particular, $|\ti\cellgens_*(u)|\leq | \cellgens_*(u)|$, but the inequality in some cases can be strict. The element $u$ defines the proper face $\combcell_u$ of $\combcell$, with root $R_u$ and generator set $\ti\cellgens_*(u)$. Furthermore, the point $(\cellroot(\combcell), u)$ from the boundary of the topological cell $\topcell(\combcell)$ is sent to the topological cell $\topcell(\combcell_u)=\topcell(R_u, \ti\cellgens_*(u))$ by the map
\[
    (R(\combcell),u)\mapsto (\cellroot_u,\ti u)
    \in\setof*{\cellroot_u}\times(0,1)^{\ti\cellgens_*(u)}
    =\topcell^\circ \bigl(\cellroot_u,\ti\cellgens_*(u)\bigr),
\]
where $\ti u\colon\ti\cellgens_*(u)\to(0,1)$ is defined as
\begin{equation}\label{eq:partialkappa2}
    \ti u(\ti\cellgen)=
        \max\setdef*[\big]{u(\cellgen)}{\cellgen\in\cellgens_*(u), R_u\cap\cellgen=\ti\cellgen}, \ \ti\cellgen\in \ti\cellgens_*(u).
\end{equation}
Then, finally, the value $\D\kappa_{\combcell}((R(\combcell), u))$ via the attaching map of $\D\topcell(\combcell)$ is given by
\begin{equation}\label{eq:partialkappa}
    \partial\kappa_{\combcell}\bigl((\cellroot(\combcell), u)\bigr)=
    \kappa_{\combcell_u}( (\cellroot_u,\ti u))\in\sk_{q-1}\ti\cw.
\end{equation}
This in the case $\ti\cellgens_*(u)=\emptyset$ says that $\partial\kappa_{\combcell}$ sends $(\cellroot(\combcell), u)$ to the vertex $\cellroot_u\in \sk_0\ti\cw$.
\begin{rem}[label=kappacube1] Let us summarize and refine some details.
    The element $u\in\D[0,1]^{\cellgens(\combcell)}$ defines two topological cubes at two different levels and also a cell of $\sk_{q-1}\ti\cw$.

    The first cube is a $\abs{\cellgens_*(u)}$-face sitting in the boundary of $\topcell(\combcell)=[0,1]^{\cellgens(\combcell)}$. Its interior is
    \[
        \ti D^\circ (u) :=\setdef{v\colon\cellgens(\combcell)\to [0,1]}{\cellgens_0(v)=\cellgens_0(u),\ \cellgens_1(v)=\cellgens_1(u)},
    \]
    while its closure is
    \[
        \ti D (u)=\setdef{v\colon\cellgens(\combcell)\to [0,1]}{\cellgens_0(v)\supseteq \cellgens_0(u),\ \cellgens_1(v)\supseteq \cellgens_1(u)}.
    \]

    The second cube is the topological cube/cell $\topcell (R_u,\ti \cellgens_*(u))$ associated with a face of $\combcell$, namely with the combinatorial $|\ti\cellgens_*(u)|$-cube $\combcell_u=(R_u,\ti \cellgens_*(u) )$. Its interior is $\topcell^\circ (R_u,\ti \cellgens_*(u))$.

    The map from \refeqn{eq:partialkappa2} (with the substitutions $u\rightsquigarrow v$) defines a map $v\in \ti D^\circ (u)\mapsto \ti v \in \topcell^\circ (R_u,\ti \cellgens_*(u))$. This extends naturally to a continuous map $\delta_{\ti D(u)}:\ti D(u)\to\topcell (R_u,\ti \cellgens_*(u))$ by the very same formula \refeqn{eq:partialkappa2} (we simply allow $v$ to take values of $0$ or $1$ for elements in $\cellgens(\combcell)\setminus(\cellgens_0(u)\cup\cellgens_1(u))$ too).

    If $\ti\cellgens_*(u)= \cellgens_*(u)$ then $\delta_{\ti D(u)}$ basically is the identity, however, if $|\ti\cellgens_*(u)|<|\cellgens_*(u)|$ then $\delta_{\ti D(u)}$ is a non-linear topological surjective contraction with all fibers contractible. For a concrete case see Example \thmref{kep}.

    The third object is the image of $\topcell(R_u,\ti \cellgens_*(u))$ via $\kappa_{\combcell_u}$ in $\sk_{q-1}\ti\cw$, it is a $|\ti\cellgens_*(u)|$-cell of this CW complex. Its interior is homeomorphic with $\topcell^\circ (R_u,\ti \cellgens_*(u))$, but the boundary of $\topcell(R_u,\ti \cellgens_*(u))$ might be contracted under the gluing procedure.

    Regarding the gluing map $\partial\kappa_{\combcell}$, for any fixed $u$, the gluing principle is the following: among the coordinates of $u$, the $1$'s determine the root $R_u$ of the open cell we are gluing into, the $0$'s determine which codimension-$1$ subspaces we omit from among the old generators, and the values in $(0,1)$ identify the exact point in that face $\topcell^\circ (R_u,\ti \cellgens_*(u))$ to which we attach the point $\bigl(\cellroot(\combcell),u\bigr)$. Finally, the map $\cellgens_*(u)\to\ti \cellgens_*(u)$, $G\mapsto R_u\cap G$, guides the cube-contraction $\delta_{\ti D (u)}$ via $v\mapsto \ti v$ (cf. \refeqn{eq:partialkappa2}).
\end{rem}

\begin{rem}[label={gluing-extension}]
    Once the continuity of $\partial\kappa_{\combcell}:\D\topcell_q(\combcell)\to\sk_{q-1}\ti\cw$ is verified (cf. Proposition \thmref{gluing-cont}), we can also consider the inclusion $\iota_{\combcell}:\D\topcell_q(\combcell)\to \topcell_q(\combcell)$ too. The fibered coproduct of these two maps gives the following commutative diagram (which defines/identifies $\kappa_{\combcell}$ as well)
    \[\begin{array}{ccc}
        \D\topcell_q(\combcell) & \stackrel{\partial\kappa_{\combcell}}{\longrightarrow} & \sk_{q-1}\ti\cw \\
        \downarrow\vcenter{\rlap{$\iota_{\combcell}$}} & & \downarrow \\
        \topcell_q(\combcell) & \stackrel{\kappa_{\combcell}}{\longrightarrow} & \sk_{q-1}\ti\cw\sqcup \topcell_q(\combcell)
    \end{array}\]
    From the continuity of $\D\kappa_{\combcell}$ the continuity of $\kappa_{\combcell}$ follows automatically.

    Regarding the map $\kappa_{\combcell}$ note the following fact as well. The formula \refeqn{eq:partialkappa2} can be considered even if $u\in\topcell^\circ(\combcell)$, i.e. for $u$ not on the boundary. In that case, we simply have $\cellgens_0(u)=\cellgens_1(u)=\emptyset$, and $\cellroot_u=\cellroot(\combcell)$, $\ti\cellgens_*(u)=\cellgens(\combcell)$, with the map leaving such points fixed. This extension of $\D\kappa_{\combcell}$ to $\topcell_q(\combcell)$ is the map $\kappa_{\combcell}$ in the diagram.

    Once the inductive construction of $\ti\cw$ is completed, we can consider the composition
    \[
        \topcell_q(\combcell)\stackrel{\kappa_{\combcell}}{\longrightarrow}\sk_{q-1}\ti\cw\sqcup\topcell_q(\combcell)\hookrightarrow\sk_{q}\ti\cw\hookrightarrow\ti\cw.
    \]
    This is the characteristic map of the cell $\topcell_q(\combcell)$ in $\ti\cw=\ti\cw(\spsystem)$, still denoted by $\kappa_{\combcell}$.
\end{rem}

Hence, in order to finish the inductive construction of the CW complex $\ti\cw$ , we need the following:

\begin{prop}[label=gluing-cont]
    Assume that $\sk_{q-1}\ti \cw$ was already constructed (hence all the maps $\D\kappa_{\combcell}$ and $\kappa_{\combcell}$ are continuous for any $\spsystem$-cell $\combcell$ of dimension $<q$). Then, for any $\combcell$ with $|\cellgens(\combcell)|=q$, $\D\kappa_{\combcell}$ is continuous.
\end{prop}
\begin{pf}
    Let us write $\cellroot=\cellroot(\combcell)$ and $\cellgens=\cellgens(\combcell)$.

    In order to prove the statement, it is enough to prove that the restriction of $\D\kappa_{\combcell}$ to any closed face of $\topcell(\combcell)$ is continuous. The interior of such a face can be defined as
    \[
        \ti D^\circ=\ti D^\circ_{\cellgens_0,\cellgens_1}=\setdef{u\in[0,1]^{\cellgens}}{\cellgens_0(u)=\cellgens_0,\,\cellgens_1(u)=\cellgens_1}\subseteq\D[0,1]^{\cellgens},
    \]
    for some fixed $\cellgens_0,\cellgens_1\subseteq\cellgens$, not both of them empty. Let $\ti D$ be its closure in the cube $[0,1]^{\cellgens}$. They can be identified with $\ti D^\circ(u)$ and $\ti D(u)$ for a certain $u$ considered in Remark \thmref{kappacube1}.

    Also, let $\cellgens_*=\cellgens\setminus(\cellgens_0\cup\cellgens_1)$.

    If one only considers $\D\kappa_{\combcell}|_{\ti D^\circ}$, we see that $R_u$ and $\ti\cellgens_*(u)$ are fixed, and thus the continuity can be seen directly from the continuity of the max-function in \refeqn{eq:partialkappa2}. In other words, $\D\kappa_{\combcell}$ would be continuous if the topology on $\sk_{q-1}\ti\cw$ was obtained simply via taking the disjoint union of its open cells. Now we need to check that $\D\kappa_{\combcell}$ is in fact compatible with the attaching maps corresponding to the cells in $\sk_{q-1}\ti\cw$. That is to say, for any given such open face $\ti D^\circ=\ti D^\circ _{\cellgens_0,\cellgens_1}$, let $\cellroot_{\ti D}$ and $\ti\cellgens_*(\ti D)$ be the values of $\cellroot_u$ and $\ti\cellgens_*(u)$ for all $u\in \ti D^\circ $, $\combcell_{\ti D}=\combcell\bigl(R_{\ti D},\ti\cellgens_*(\ti D)\bigr)$, and extend continuously the formula in \refeqn{eq:partialkappa2} to a map
    \[
        \delta_{\ti D}\colon {\ti D}\to\topcell(\combcell_{\ti D})
    \]
    as in Remark \thmref{kappacube1} (we simply allow $v$ to take values of $0$ or $1$ for elements in $\cellgens\setminus(\cellgens_0\cup\cellgens_1)$ too as in the definition of $\delta_{\ti D(u)}$). Then, we claim that this map satisfies the identity
    \begin{equation}\label{eq:claimidentity}
        \kappa_{\combcell_{\ti D}}\circ\delta_{\ti D}=\kappa_{\combcell}|_{ \ti D}.
    \end{equation}
    This can be checked directly: for any element $(\cellroot,v)\in \ti D=\setof{\cellroot}\times\setdef*[\big]{v\in[0,1]^{\cellgens}}{v|_{\cellgens_0}=0,\,v|_{\cellgens_1}=1}$,
    \[
        \delta_{\ti D}\bigl((\cellroot,v)\bigr)=(\cellroot',v')\in\setof*{\cellroot'}\times[0,1]^{\cellgens'}, \ns
        \cellroot'=\cellroot_{\ti D}, \ns
        \cellgens'=\ti\cellgens_*(\ti D), \ns
        v'(\cellgen')=\max\setdef*[\big]{v(\cellgen)}{\cellroot'\cap\cellgen=\cellgen'}.
    \]
    That point $(\cellroot',v')$ is in turn taken by $\kappa_{\combcell_{\ti D}}$ to some $(\cellroot'',v'')\in[0,1]^{\cellgens''}$, where
    \begin{align*}
        \cellroot''
            & =\cap\,\bigl(\setof*{\cellroot'}\cup\setdef*{\cellgen'}{v'(\cellgen')=1}\bigr)
            =\cap\,\bigl(\setof*{\cellroot'}\cup\setdef*{\cellgen'\in\cellgens'}{\exists\cellgen: v(\cellgen)=1, \cellroot'\cap\cellgen=\cellgen', \cellgen\in\cellgens_*}\bigr)= \\
        & =\cap\,\bigl(\setof*{\cellroot'}\cup\setdef*{\cellroot'\cap\cellgen}{\cellgen\in\cellgens_*, v(\cellgen)=1, \cellgen\not\supseteq\cellroot'}\bigr)
            =\cap\,\bigl(\setof{\cellroot}\cup\cellgens_1\cup\setdef{\cellgen\in\cellgens_*}{v(G)=1}\bigr)
            =\cellroot_v, \\
        \cellgens''
            & =\setdef*{\cellroot''\cap\cellgen'}{\cellgen'\in\cellgens', v'(\cellgen')\in(0,1)}\setminus\setof*{\cellroot''}= \\
        & =\setdef*{\cellroot_v\cap(\cellroot_{\ti D}\cap\cellgen)}{\cellgen\in\cellgens_*,\cellgen\not\supseteq\cellroot_{\ti D}, v(\cellgen)\in(0,1)}\setminus\setof*{\cellroot_v}
        =\setdef{\cellroot_v\cap\cellgen}{\cellgen\in\cellgens_*(v)}\setminus\setof{\cellroot_v}= \\
        & =\ti\cellgens_*(v),
    \end{align*}
    and for any $\cellgen''\in\cellgens''=\ti\cellgens_*(v)$, we have
    \begin{align*}
        u''(\cellgen'')
            & =\max\setdef{v'(\cellgen')}{\cellgen'\in\cellgens',\cellroot''\cap\cellgen'=\cellgen''}
            =\max\setdef{v(\cellgen)}{\cellgen\in\cellgens,\cellroot''\cap\cellgen=\cellgen''}= \\
        & =\max\setdef{v(\cellgen)}{\cellgen\in\cellgens,\cellroot_v\cap\cellgen=\cellgen''}.
    \end{align*}
    So, the equality \refeqn{eq:claimidentity} is indeed true.

    From \refeqn{eq:claimidentity}, the continuity of $\kappa_{\combcell}|_{\ti D}$ follows because both maps on the left hand side are continuous: $\kappa_{\combcell_{\ti D}}$ because of the induction hypothesis (since $\combcell_{\ti D}$ is lower dimensional), $\delta_{\ti D}$ by its definition (see above). Thus, $\kappa_{\combcell}$ is continuous when restricted to any closed face, and as these sets cover $\D\topcell(\combcell)$, the map $\kappa_{\combcell}$ is continuous on the whole.
\end{pf}

\begin{ex}[label=kep]
    To get a better intuition about the above picture, it is worth considering an example. Let all elements of $\spsystem$ be subspaces of a $2$-dimensional space, with $\setof{e_1,e_2}$ the canonical basis, and $\cellroot=\gen{e_1,e_2}$, $\cellgen_1=\gen{e_1}$, $\cellgen_2=\gen{e_2}$, $\cellgen_3=\gen{e_1+e_2}$, $\celltip=\setof{0}$ are all elements of $\spsystem$.

    These vertices together form a $3$-dimensional combinatorial cell $\combcell$ of height $2$: the names of vertices reflect that. Between them we also have $6$ edges in total: one connecting each ($1$-dimensional) generator to both $\cellroot$ and $\celltip$. There will also be $3$ different $2$-dimensional cells.

    Since $\height\combcell<\dim\combcell$, we can see that some of the contraction maps onto the boundary topological cells will be nontrivial: some nontrivial configurations of points on the boundary of $\topcell(\combcell)=[0,1]^3$ will be mapped to the same point. The following diagram illustrates this:

    \begin{center}
        \(\vcenter
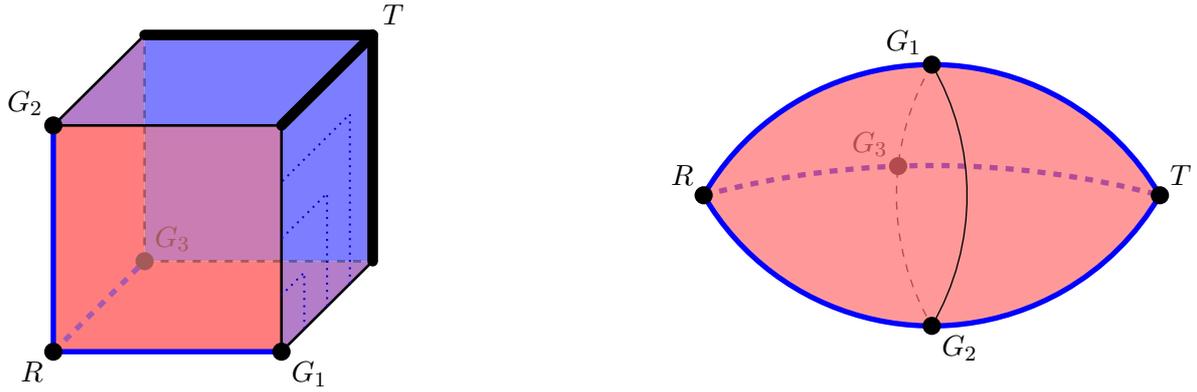
{\hbox{
        \begin{tikzpicture}[scale=3.0]
            \draw [draw=none, fill=red, fill opacity=0.3] (0,0) -- (0.4,0.4) -- (0.4,1.4) -- (0,1) -- (0,0);
            \draw [draw=none, fill=red, fill opacity=0.3] (0,0) -- (0.4,0.4) -- (1.4,0.4) -- (1,0) -- (0,0);
            \draw [draw=none, fill=blue, fill opacity=0.3] (0.4,0.4) -- (1.4,0.4) -- (1.4,1.4) -- (0.4,1.4) -- (0.4,0.4);

            \draw[dashed, color={black!50!white}, line width=1pt] (0.4,0.4) -- (1.4,0.4);
            \draw[dashed, color={black!50!white}, line width=1pt] (0.4,0.4) -- (0.4,1.4);
            \draw[line width=2pt, color={blue!50!white}, dashed] (0,0) -- (0.4,0.4);

            \draw[draw=none, fill={black!50!white}] (0.4,0.4) circle [radius=0.04];
            \node [color={black!50!white}] at (0.4,0.4) [above right] {$\cellgen_3$};

            \draw [draw=none, fill=red, fill opacity=0.3] (0,0) -- (1,0) -- (1,1) -- (0,1) -- (0,0);

            \draw [draw=none, fill=blue, fill opacity=0.3] (1,0) -- (1.4,0.4) -- (1.4,1.4) -- (1,1) -- (1,0);
            \draw [draw=none, fill=blue, fill opacity=0.3] (0,1) -- (0.4,1.4) -- (1.4,1.4) -- (1,1) -- (0,1);
            \draw[line width=2pt, blue] (0,0) -- (1,0);
            \draw[line width=2pt, blue] (0,0) -- (0,1);
            \draw (1,0) -- (1.4,0.4) -- (1.4,1.4) -- (1,1);
            \draw [line width=1pt] (1,0) -- (1.4,0.4);
            \draw [line width=1pt] (1,0) -- (1,1);
            \draw [line width=1pt] (0,1) -- (0.4,1.4);
            \draw [line width=1pt] (0,1) -- (1,1);

            \draw [draw=none, fill=black] (0,0) circle [radius=0.04];
            \draw [draw=none, fill=black] (1,0) circle [radius=0.04];
            \draw [draw=none, fill=black] (0,1) circle [radius=0.04];
            \draw [line width=4pt, line cap=round] (1.4,1.4) -- (1,1);
            \draw [line width=4pt, line cap=round] (1.4,1.4) -- (0.4,1.4);
            \draw [line width=4pt, line cap=round] (1.4,1.4) -- (1.4,0.4);

            \node at (0,0) [below left] {$\cellroot$};
            \node at (1,0) [below right] {$\cellgen_1$};
            \node at (0,1) [above left] {$\cellgen_2$};
            \node at (1.4,1.4) [above right] {$\celltip$};

            \draw [color={blue!80!black}, dotted, line width=0.75pt] (1.10,0.10) -- (1.10,0.35) -- (1,0.250);
            \draw [color={blue!80!black}, dotted, line width=0.75pt] (1.20,0.20) -- (1.20,0.70) -- (1,0.500);
            \draw [color={blue!80!black}, dotted, line width=0.75pt] (1.30,0.30) -- (1.30,1.05) -- (1,0.750);
        \end{tikzpicture}
        }}\hskip3cm\vcenter{\hbox{
        \begin{tikzpicture}[scale=3.0]
            \draw [line width=2pt, color={blue!50!white}, dashed] (-1,0) arc (105:75:{1/sin(15)});
            \draw [line width=0.5pt, color={black!50!white}, dashed] (0,{tan(30)}) arc (150:210:{1/cos(30)});
            \draw [draw=none, fill={black!50!white}] (-0.1475,0.1288) circle [radius=0.04];

            \node [color={black!50!white}] at (-0.1475,0.1288) [above left] {$\cellgen_3$};

            \draw [draw=none, fill={red}, fill opacity=0.4] (-1,0) arc (150:30:{1/sin(60)}) arc (-30:-150:{1/sin(60)});

            \draw [line width=2pt, color={blue}] (-1,0) arc ( 150: 30:{1/sin(60)});
            \draw [line width=2pt, color={blue}] (-1,0) arc (-150:-30:{1/sin(60)});
            \draw [line width=0.5pt, color={black}] (0,{tan(30)}) arc (30:-30:{1/cos(30)});

            \draw [draw=none, fill={black}] (-1,0) circle [radius=0.04];
            \draw [draw=none, fill={black}] ( 1,0) circle [radius=0.04];
            \draw [draw=none, fill={black}] (0,{ tan(30)}) circle [radius=0.04];
            \draw [draw=none, fill={black}] (0,{-tan(30)}) circle [radius=0.04];

            \node at (-1,0) [above left] {$\cellroot$};
            \node at (1,0) [above right] {$\celltip$};
            \node at (0,{ tan(30)}) [above left ] {$\cellgen_1$};
            \node at (0,{-tan(30)}) [below right] {$\cellgen_2$};
        \end{tikzpicture}
        }}\)
    \end{center}
    The black dots and the set of $3$ thick black edges on the first diagram represent what will collapse into a single vertex each upon the gluing: the latter ones will be attached to the vertex $T$. The $3$ blue edges and $3$ blue faces (one in the back) will collapse onto a single edge each. The dotted lines on the rightmost blue face show some points that will be glued to the same point on that edge. Lastly, the red faces will become the $2$-cells.

    What we thus end up with is a single $3$-cell looking like what is seen on the second diagram. This is the image of $\topcell(\ti C)$ in $\ti\cw$. (The black edges $G_1G_2, \ G_2G_3, \ G_3G_1$ in the middle are not cells in $\cw$, they serve only to illustrate the spatial arrangement of the faces.)
\end{ex}

Repeating the gluing for each integer $q\geq 0$ (taking the limit if the set of dimensions is not bounded), we get a CW complex $\ti\cw(\spsystem)$.

\begin{defn}[label=cw-def]
    For (an intersection-closed) system of vector spaces $\spsystem$, its associated CW complex $\ti\cw(\spsystem)$ is that obtained in the above-described procedure.
\end{defn}

\begin{rem}[label=cw-dim-finite]
    \begin{enumerate}[label={(\alph*)}]
        \item For general $\spsystem$, the dimension $\dim\ti\cw(\spsystem)$ may well be infinite. However, if $\spsystem=\setof{\filt(\ell)}_{\ell}$, where $\ell$ belongs to some lattice (as in \ref{ss:prel5}) then $\dim\ti\cw(\spsystem)$ is less than or equal to the rank of the lattice (cf. \thmref{lattice-to-sscw}).
        
        \item Though the above construction $\spsystem \rightsquigarrow \ti \cw(\spsystem)$ can be applied for any $\spsystem$, the combinatorial/topological cells were defined with one eye on the properties of the system of subspaces associated with curve singularities, namely on Lemma \thmref{ss-set-props2}\It{(c)}. In the language of system of subspaces this means that the generators have codimension one (in Hilbert function language $\hh(\ell+E_i)-\hh(\ell)\in\setof{0,1}$). This in particular means that the above construction is adapted to the curve singularities and their lattice cohomology (or any other subspace systems with similar properties), and it should be modified for more general subspace systems, e.g. in the case of higher dimensional singularities. The corresponding general construction will be the subject of a subsequent manuscript.
    \end{enumerate}
\end{rem}

Though the intersection-closed system $\spsystem$ can be very general, the homotopy type of the topological space $\ti\cw(\spsystem)$ is not very complicated (the additional system of weight makes it really interesting). Recall that $\cw=R(0,c)$ considered in earlier sections is contractible. Rather similar property holds for $\ti\cw(\spsystem)$ as well.

\begin{lem}[label=contractible]
    \begin{enumerate}[label={(\alph*)}]
        \item A pair of vertices given by
         $V,W\in\spsystem$ belongs to the same connected component of $\ti\cw(\spsystem)$ if and only if there exists $U\in\spsystem$, $U\subset V\cap W$, such that the inclusions $U\subset V$ and $U\subset W$ can be filled by a sequences of codimension one subspaces from $\spsystem$ (i.e.~there exists $U=V_0\subset$ $ V_1\subset \cdots \subset V_n=V$, $V_i\in\spsystem$, $\codim(V_i\subset V_{i+1})=1$, and a similar sequence for $U\subset W$).
        \item All the connected components of $\ti\cw(\spsystem)$ are contractible.
    \end{enumerate}
\end{lem}
\begin{pf} (a)
    Two vertices $V$ and $W$ in $\ti\cw(\spsystem)$ belong to the same connected component if and only if they are connected by a sequence of edges, that is, if they can be connected by the concatenation of an arbitrary number of sequences of type given in the statement. Hence we need only to prove that any such general sequence can be replaced by one as in the statement. Let $V=V_0,V_1,\ldots,V_t=W$ be an arbitrary sequence of edges. Any edge $(V_{i-1},V_i)$ may either go `upward' (when $V_{i-1}\subset V_i$) or `downward' (when $V_{i-1}\supset V_i$). We intend to get to a new path where we first `go down' from $V$ to a subspace $U\in\spsystem$, then up to $W$, i.e. where all downward moves precede all upward moves.

    Suppose we are not yet in that state, i.e.\ $V_{i-1}\subset V_i\supset V_{i+1}$ for some $i$. Since $\spsystem$ is closed under intersection, $V_i^*=V_{i-1}\cap V_{i+1}\in\spsystem$ holds, and we can change $V_i$ to $V_i^*$ instead. In fact, the paths before and after are homotopic due to the existence of the $2$-cell $\topcell(V_i^*,\setof{V_{i-1}, V_{i+1}})$. The exception is when $V_{i-1}=V_{i+1}$, but in this case we can simply omit $V_i$, and we will still get a homotopy (due to the edge $\topcell(V_i,\setof{V_{i-1}})$.

    Repeating this process, thus rearranging the `up' and `down' moves, we can obtain a path of the desired type, and a homotopy (with fixed endpoints) connecting it to the original.

    (b) We note that a similar argument proves that all connected components are simply connected. We start with any loop in $\ti\cw(\spsystem)$, which is known to be homotopic to one within $\sk_1\cw(\spsystem)$: some sequence $V_0,V_1,\ldots,V_t=V_0$. On this, we perform the same type of process: we pick an orientation in the loop, and any two consecutive steps of type `up' and the `down' either get removed, or are replaced with a `down' and an `up' step. Since we take intersections, we always contain the intersection of all original $V_i$'s, thus this process will end in a finite number of steps; and that can only happen when the length of the loop is $0$ (since it can not contain only steps going up or only those going down). Indeed, we \It{always} end up in the vertex $\bigcap_i V_i$.

    Thus, every map $S^1\to\ti\cw(\spsystem)$ can be contracted to a constant map. In general, for maps $S^k\to\ti\cw(\spsystem)$ $(k\geq 2)$ we proceed in a slightly different way though (a proof which works for $k=1$ as well).

    Since for any map $\psi\colon S^k\to\ti\cw(\spsystem)$ ($k\geq 1$), the set $\Im\psi$ is compact, that image is contained in a finite number of open cells. Consider the vertices of all these open cells, and also all intersections formed by them: we denote this finite set by $\spsystem^*$ ($\spsystem^*\subset\spsystem$). In fact, $\ti\cw^*=\ti\cw(\spsystem^*)$ is a connected subcomplex of $\ti\cw(\spsystem)$ containing $\Im\psi$. This has a minimal vertex $U=\cap_{W\in \spsystem^*}W$, which we will contract it to. We prove this by induction, at first step we contract $\ti\cw^*$ to a subcomplex spanned by a smaller set of vertices. Repeating this process, due to the $\abs{\spsystem^*}<\infty$, by eliminating all but $U$.

    Suppose we have more than one vertex in $\ti\cw^*$. Pick any vertex $V\in\spsystem^*$ that is maximal within $\spsystem^*$ (there may be more than one), and let $\cellgens=\setdef{W \in\spsystem^*}{W\subset V,\,\codim(W\subset V)=1}$. Since this set is a finite, and also non-empty (due to there existing a sequence of edges within $\ti\cw^*$ from $V$ to $U$), we can consider $\topcell=\topcell\combcell(V,\cellgens)$, the unique largest cell in $\ti\cw^*$ containing $V$. Since all cells of $\ti\cw^*$ that contain $V$ (and are thus rooted at $V$) are faces of $\topcell$, we can just contract it to its faces not containing $V$. More precisely, let $\topcell=\setof{V}\times [0,1]^{\cellgens}$, and for any $u\colon\cellgens\to[0,1]$, we set
    \[
        u'=u+u_{const},
    \]
    where $u_{const}\colon\cellgens\to [0,1]$ is the constant function with value
    $(1-\max u)$. Note that
    \begin{enumerate}[label={(\arabic*)}]
        \item $u'(\cellgen)=1$ for at least one $\cellgen\in\cellgens$, thus $\kappa_{\combcell}(u')$ is on some smaller dimensional face of $\topcell$, but not in any open face rooted at $V$;
        \item $u'=u$ for all such points (where a coordinate is $1$);
        \item $u\mapsto u'$ is continuous.
    \end{enumerate}
    Thus, we can obtain a homotopy $H\colon\ti\cw^*\times[0,1]\to\ti\cw^*$ where
    \begin{alignat*}{3}
        H\bigl(\kappa_{\combcell}(u), t\bigr) & =\kappa_{\combcell}\bigl((1-t)u+tu'\bigr), && \quad u\in[0,1]^{\cellgens} \\
        H\bigl(p,t) & = p, && \quad p\in\ti\cw(\spsystem^*\setminus\setof{V}).
    \end{alignat*}
    By the previous observations, this is indeed the desired strong deformation retract as it is well-defined and continuous on both $\topcell\times[0,1]$ and $\ti\cw(\spsystem^*\setminus\setof{V})\times[0,1]$ (both being closed sets whose union is $\ti\cw^*$), and for $H_t\colon p\mapsto H(p,t)$, $H_0$ is the identity, $H_1$ maps to $\ti\cw(\spsystem^*\setminus\setof{V})$, while $H_t$ is always the identity when restricted to $\ti\cw(\spsystem^*\setminus\setof{V})$.

    Repeating this step, we can contract $\ti\cw^*$ to its minimal vertex. Thus in particular the original map $\psi$ is always homotopic to a constant map, hence all homotopy groups save for $\pi_0$ are trivial. Hence, given that $\ti\cw$ is a CW complex, the connected components are contractible.
\end{pf}
\begin{rem}[label=notcontractible]
    The above contractibility of the connected components of $\ti\cw(\spsystem)$ basically follows from the fact that $\spsystem$ is intersection-closed. If we drop this assumption, then we still can define the set of combinatorial cells $\combcell(\cellroot,\cellgens)$, where $R$ is a root $V$, the generator $\cellgens$ is a set of subspaces $W$ of codimension one in $V$, but in this case in the definition of the cell we impose additionally that any intersection of type $\cap_{G\in {\mathcal H}} G$, for ${\mathcal H}\subset \cellgens$, should belong to $\spsystem$. Then, in a similar way as above, we glue the associated topological cells into a CW complex $\ti\cw(\spsystem)$. However, in this new situation, $\ti\cw(\spsystem)$ might have non-trivial homotopy type. E.g., assume that $V=\gen{e_1,e_2,e_3}$ is 3-dimensional, and (the non-intersection closed) $\spsystem$ consists of the subspaces $0$, $\gen{e_1}$, $\gen{e_1,e_3}$, $V$, $\gen{e_2,e_3}$, $\gen{e_2}$. Then $\ti\cw(\spsystem)\sim S^1$.
\end{rem}

Next we introduce weights to this complex as well.

\begin{defn}[label=wss-def]
    We call $(\spsystem,w)$ a \It{weighted system of vector spaces} (WSS) if $\spsystem$ is an intersection closed system of vector spaces as in Definition \thmref{V-cell}, and $w\colon\spsystem\to\bbz$ is a function bounded from below. Additionally, for any $V,W\in\spsystem$ with $W\subseteq V$, we require
    \begin{equation}\label{eq:weightfun}
        w(W)\leq w(V)+\codim(W\subseteq V).
    \end{equation}
\end{defn}
We invite the reader to recall the definition of weighted CW complex from Definition \thmref{def:compat-weight}.

Next, for any WSS $(\spsystem,w)$, we introduce a weight function on its associated CW complex.

\begin{defn}[label=V-cell-weight]
    Let $(\spsystem, w)$ be a weighted system of vector spaces, and consider the associated complex $\ti\cw=\ti\cw(\spsystem)$. For a $\spsystem$-cell $\combcell$, let
    \[
        \ti w\bigl(\topcell(\combcell)\bigr)=w\bigl(\cellroot(\combcell)\bigr)+\height\combcell.
    \]
    \begin{prop}[label=V-cell-weight-compat]
        The defined $\ti w$ forms a set of compatible weight function on $\ti\cw(\spsystem)$.
    \end{prop}
    \begin{pf}
        This is a consequence of the definition \thmref{wss-def}, since for any face $\combcellface$ of a $\spsystem$-cell $\combcell$,
        \begin{align*}
            \ti w\bigl(\topcell(\combcellface)\bigr)
                & =w\bigl(\cellroot(\combcellface)\bigr)+\height\combcellface
                    \leq w\bigl(\cellroot(\combcell)\bigr)+\codim\bigl(\cellroot(\combcellface)\subseteq \cellroot(\combcell)\bigr)+\codim\bigl(T(\combcellface)\subseteq R(\combcellface)\bigr)\leq \\
                & \leq w\bigl(\cellroot(\combcell)\bigr)+\codim\bigl(\celltip(\combcell)\subseteq \cellroot(\combcell)\bigr)
                    =\ti w\bigl(\topcell(\combcell)\bigr). \qedhere
        \end{align*}
    \end{pf}
    The resulting pair $(\ti\cw,\ti w)$ is called the \It{weighted CW complex associated to $(\spsystem, w)$}.
\end{defn}

\subsection{The weighted system of subspaces associated with a curve singularity}\label{ss:linsubspacesimf}\hfill\smallbreak

Let $(C,o)$ be an isolated curve singularity and consider the associated weighted complex $(\cw, w)=(R(0,c),w)$ together with its cubes $\setof{\cellset_q}_q$, and the filtration $\setof{\filt(\ell)}_{\ell}$ introduced in section \ref{sec:curveslattice}. Note that $\setof{\filt(\ell)}_{\ell}$ is a linear subspace system of $\fnring=\fnring_{C,o}$. We denote it by $\Im\filt$. We write $\grverts$ for the index set of the irreducible components of $(C,o)$, i.e. $r=\abs{\grverts}$.

Next we compare the cubes of $R(0,c)$ and of $\Im\filt$ and of the CW complex associated with $\Im\filt$. Furthermore, we enhance $\Im\filt$ with a weight system induced by the weights of $R(0,c)$ (compatibly with Definition \thmref{V-cell-weight}).
We start with the following notation:
\[
    \ti{\cellset}_q = \setof{\text{combinatorial $\Im\filt$-cells of dimension $q$}}, \quad
    \ti{\cellset}_* = \cup_q\ti{\cellset}_q.
\]

\begin{prop}[label={lattice-to-sscw}]
    \begin{enumerate}[label={(\alph*)}]
        \item For any $q$-cube $\cell\in\cellset_q$, the set
        \[
            \combcell=\setdef{\filt(\ell)}{\ell\in \cell\cap L}
        \]
        is a combinatorial $\Im\filt$-cell of dimension at most $q$ (where $\dim C>\dim\combcell$ can happen).

        Let $\phi:\cellset_*\to\ti{\cellset}_*$ be the map that takes a cube $\cell$ to the thus associated $\combcell$.
        \item Any combinatorial $\Im\filt$-cell $\combcell$ can be realized as $\phi(C)$ for some $C\in \cellset_*$.

        In fact, for each $\Im\filt$-cell $\combcell$ of dimension $q$, there exists (usually not unique) $q$-cube $\cell\in\cellset_q$ with $\phi(\cell)=\combcell$. We call such a cube $C$ \It{minimal representative} of $\combcell$.
    \end{enumerate}

    In particular, all the $\Im\filt$-cells have dimension no greater than the rank of the lattice $\bbz^r$, namely $r=\text{number of components of }(C,o)$.
\end{prop}
\begin{pf}
    (a)
    Fix $C=(\ell, I)$. Then consider $V=\filt(\ell)$ and $V_i=\filt(\ell+E_i)$ for all $i\in I$. Then $\cap_{i\in J}V_i=\filt(\ell+E_J)$ for any $\emptyset\neq J\subseteq I$. Hence, $\combcell$ as a combinatorial cell has root $V$ and generators $\setdef{V_i}{V_i\neq V}$. Clearly $\dim(\combcell)\leq\abs{I}$.

    (b)
    For the other direction, the case $q=0$ is trivial, so we will assume $q>0$. With the notation introduced after Lemma \thmref{uniquemax_0}, consider the lattice points $\maxpoint_V$ and $\maxpoint_\cellgen$ ($\cellgen\in\cellgens(\combcell)$). Since $\cellgen\subsetneq V$, we get $\maxpoint_\cellgen\gneqq\maxpoint_V$. Define $I_G(\ell^*_V)\subset\grverts$ as the collection of those indices $i\in \grverts$ such that $\maxpoint_V<\maxpoint_V+E_{i}\leq\maxpoint_\cellgen$. Then $I_G(\ell^*_V)\not=\emptyset$. Moreover, since ${\rm codim}(G\subset V)=1$, we have $\filt(\maxpoint_V+E_{i})=\cellgen$ for any $i\in I_G(\ell^*_V)$. (Note that $I_G(\ell^*_V)$ can be defined also as the set of indices $i\in\grverts$ such that $\filt(\maxpoint_V+E_{i})=\cellgen$, since in this case $\maxpoint_V<\maxpoint_V+E_{i}\leq\maxpoint_\cellgen$ also holds.) Now for any collection of non-empty subsets $J_G(\ell^*_V)\subset I_G(\ell^*_V)$ (indexed by $G\in\cellgens(\combcell)$) set
    \[
        \cell=(\maxpoint_V,\cup_{\cellgen\in\cellgens(\combcell)}J_G(\ell^*_V)).
    \]
    Clearly $\phi(\cell)=\combcell$. Moreover, if each $J_G(\ell^*_V)$ contains exactly one element, say $i_G$, then all the indices $i_\cellgen$ must be distinct (since the spaces $\filt(\maxpoint_V+E_{i_\cellgen})=\cellgen$ are distinct). As such, the dimension of $\cell$ is exactly $q=|\cellgens(\combcell)|$.
\end{pf}

Next, we wish to determine all the representatives $C$ of a fixed $\combcell=\combcell(V,\cellgens)$, $\cellgens\not=\emptyset$. Clearly, if $C=(\ell, I)\in \phi^{-1}(\combcell)$, then $\filt(\ell)=V$. Moreover, for any $G\in\cellgens$ there must exist $i\in I$ such that $\filt(\ell+E_i)=G$, and for any $i\in I$ the space $\filt(\ell+E_i)$ either belongs to $\cellgens$ or it is $V$.

Let us fix $\ell$ with $\filt(\ell)=V$. Let $\phi^{-1}(\combcell)_{\ell}$ be the set of cells of type $(\ell, I)\in\phi^{-1}(\combcell)$. For any $G\in \cellgens$ set $I_G(\ell):= \setdef{i\in \grverts}{\filt(\ell+E_i)=G}$, and write $I^*(\ell)$ for the support $\abs{\ell^*_V-\ell}=\setdef{i}{E_i\leq \ell^*_V-\ell}$. Then $C=(\ell,I)\in\phi^{-1}(\combcell)_{\ell}$ if and only if $I\cap I_G(\ell)\neq\emptyset$ and $I\subset I^*(\ell)\cup(\cup_G I_G(\ell))$.

Note that the sets $\setof{I_G(\ell)}_G$ are all pairwise disjoint. Indeed, if $i\in I_G(\ell)\cap I_{G'}(\ell)$ then $G=\filt(\ell+E_i)=G'$. Furthermore, they are all determined by $I_G(\ell^*_V)$ (introduced in the proof of Lemma \thmref{lattice-to-sscw}) as follows.

\begin{lem} [label={lemma-ig}]
    $I_G(\ell)=I_G(\ell^*_V)\setminus\abs{\ell^*_V-\ell}$ for any $G\in\cellgens$.
\end{lem}
\begin{pf}
    If $i\in I_G(\ell)$ then $\filt(\ell+E_i)=G$, hence $i\not\in\abs{\ell^*_V-\ell}$ (otherwise, $\ell+E_i\leq \ell^*_V$, which would imply $\filt(\ell+E_i)=V$). But then $(*)$ $\filt(\ell^*_V+E_i)=\filt(\ell+E_i)\cap \filt(\ell^*_V)=G\cap V=G$.

    If $i\in I_G(\ell^*_V)\setminus\abs{\ell^*_V-\ell}$, then write $W:=\filt(\ell+E_i)$. Then $W\subset\filt(\ell)=V$. But $W=V$ cannot hold since that would imply $\ell+E_i\leq \ell^*_V$, or $i\in\abs{\ell^*_V-\ell}$. Hence $W$ has codimension 1 in $V$. But then, again from $(*)$, $W\cap V=G$, hence $W=G$.
\end{pf}
\begin{cor}[label={cor-felosrolas}] Fix $\combcell=\combcell(V,\cellgens)$ with $\cellgens\not=\emptyset$. Then the following facts hold:
    \begin{enumerate}[label={(\arabic*)}]
        \item $\phi^{-1}(\combcell)_{\ell}\not=\emptyset \ \Leftrightarrow\ I_G(\ell^*_V)\setminus|\ell^*_V-\ell| \not=\emptyset$ for any $G\in\cellgens$.
        
        \item $I_G(\ell)\subset I_G(\bar{\ell})$ for any $\bar{\ell}$ with $\ell\leq \bar{\ell}\leq \ell^*_V$.
        
        \item $\phi^{-1}(\combcell)_{\ell}\not=\emptyset \ \Rightarrow\ \phi^{-1}(\combcell)_{\bar{\ell}}\not=\emptyset$ for any $\bar{\ell}$ with $\ell\leq \bar{\ell}\leq \ell^*_V$.
        
        \item If $(\ell, I)\in \phi^{-1}(\combcell)$ and $\ell\leq \bar{\ell}\leq \ell^*_V$ then $(\bar{\ell},J)\in \phi^{-1}(\combcell) $ for any $J$ with $I\subset J\subset I\cup|\ell^*_V-\bar{\ell}|$.
        
        \item $(\ell, I)\in \phi^{-1}(\combcell)$ if and only if $\filt(\ell)=V$ and $I\subset |\ell^*_V-\ell|\cup (\cup_GI_G(\ell^*_V)) $, with $I\cap(I_G(\ell^*_V)\setminus |\ell^*_V-\ell|)\not=\emptyset$ for any $G$.
    \end{enumerate}
\end{cor}

 The above statements imply the following facts as well.
\begin{cor}[label={V-cell-faces}]
    \begin{enumerate}[label={(\alph*)}]
        \item Any cell $\cell\in \phi^{-1}(\combcell)$ (a not necessarily minimal representative of $\combcell$) contains a minimal representative as a face (not necessarily of codimension $1$).
        
        \item Given an $\Im\filt$-cell $\combcell$ and any representative $\cell$ of $\combcell$, we have
        \[
            \setdef{\phi(\cellface)}{\text{$\cellface$ is a face of $\cell$}}=\setdef*[\big]{\combcellface\subseteq\combcell}{\combcellface\text{ is an $\Im\filt$-cell}}.
        \]
        In particular, the left hand side is independent of the choice of $\cell$.
    \end{enumerate}
\end{cor}

Next we compare the weights appearing at different levels. Firstly, the CW complex $\cw=R(0,c)$ ($c$ finite or infinite) has a weight function defined as $w_0(\ell)=2\hh(\ell)-\abs{\ell}$ for any lattice point $\ell$, and $w_q(\topcell)=\max\setdef{w_0(\ell)}{\text{$\ell$ vertex of $\topcell$}}$.

On the other hand, on $\spsystem=\Im\filt$ associated with $(C,o)$, we set the following weights. First note that $V$ determines uniquely $\maxpoint_V$ (it is maximal with $V=\filt(\maxpoint_V)$). Then, we set $w_{\spsystem}(V):= w_0(\maxpoint_V)$.

Let us verify that $(\spsystem, w_{\spsystem})$ defines a weighted system of vector spaces, i.e. $ w_{\spsystem}$ satisfies \refeqn{eq:weightfun}. Indeed, for any $V,W\in\Im\filt$ with $W\subseteq V$, we have $\filt\bigl(\max\setof{\maxpoint_V,\maxpoint_W}\bigr)=W$, so $\maxpoint_W\geq\maxpoint_V$. Following an increasing path in the lattice from $\maxpoint_V$ to $\maxpoint_W$, the value of $\filt$ can only decrease in codimension-$1$ steps, and $w_0(\ell)$ increases by $1$ at each such step (and otherwise decreases by $1$).

Next, by Definition \thmref{V-cell-weight}, for any combinatorial or topological cell of $\spsystem$ or $\ti\cw(\spsystem)$ respectively, we have $\ti w_{\spsystem}(\topcell(\combcell))=w_{\spsystem}(R(\combcell))+\height\combcell$. This defines a set of compatible weight functions on $\ti\cw(\spsystem)$, or on the set of combinatorial cells of $\spsystem$.

It is worth to emphasize that $\ti w_{\spsystem}(\combcell)=w_{\spsystem}(R(\combcell))+\height\combcell$ defined in this way is \underline{not} the maximum of the $\ti w_{\spsystem}$-weights of the vertices of $\combcell$. Consider e.g. the situation from Example \thmref{ex:delta2}, the segment $[0,3]$: $\filt(0)=\fnring$, $\hh(0)=0$, $\maxpoint_{\fnring}=0$, $\filt(3)=\mathfrak{m}$ (the maximal ideal), $\hh(3)=1$, $\maxpoint_{\mathfrak{m}}=3$. The segment $[0,1]$ is sent to a 1-cell $\combcell$ in $\ti\cw$ with end-vertices $\maxpoint_{\fnring}$ and $\maxpoint_{\mathfrak{m}}$ (and the segment $[1,3]$ to $\maxpoint_{\mathfrak{m}}$). Moreover, $\ti w(\fnring)=0$, $\ti w(\mathfrak{m})=-1$, but $\ti w(\combcell)=\ti w(\fnring)+\codim (\mathfrak{m}\subset \fnring)=1$.

The next statement connects the set of weight functions of $(\cw,w)$ and $(\ti\cw(\spsystem), \ti w_{\spsystem})=(\ti\cw, \ti w)$ in the spirit of the comparisons of statements from \thmref{lattice-to-sscw}--\thmref{V-cell-faces}.

\begin{prop}[label={Im-F-cell-weight}]
    Let $V$ and $W$ be the root and tip of $\combcell\in\ti\cellset_q$.
    Then the following facts hold:
    \begin{enumerate}[label={(\alph*)}]
        \item For any representative $\cell=(\ell,I)$ of $\combcell$, we have
        \[
            w(\cell)=w(\ell)+\height\combcell= w(\ell)+{\rm codim}(W\subset V).
        \]
        
        \item
        \[
            \min\setdef*[\big]{w(\cell)}{\cell\in\cellset_*, \phi(\cell)=\combcell}=
            w(\maxpoint_V)+\height\combcell= \ti w_{\spsystem}(\combcell).
        \]
    \end{enumerate}
\end{prop}
\begin{pf}
    (a) Take a sequence of lattice points $\setof{x_i}_{i=0}^t$ such that $x_0=\ell$, $x_t=\ell+E_I$, and for any $i$ one has $x_{i+1}=x_i+E_{v(i)}$ for a certain $v(i)$. Hence $\filt(x_0)=V$ and $\filt(x_t)=W$. Then $w(x_{i+1})-w(x_i)\in\setof{1,-1}$, and it is 1 exactly when ${\rm codim}(\filt(x_{i+1})\subset \filt(x_i))=1$. Hence $w(x_i)\leq w(\ell)+{\rm codim}(W\subset V)$ for any $x_i$. Since any vertex of $C$ appears as such an $x_i$ we get $w(C)\leq w(\ell)+{\rm codim}(W\subset V)$.

    Next, we analyse three consecutive points $x_{i-1}, \ x_i,\ x_{i+1}$. Assume that $w(x_{i})-w(x_{i-1})=-1$ and $w(x_{i+1})-w(x_i)=1$. Then if we replace $x_i$ by $\bar{x}_i=x_{i-1}+E_{v(i)}$, by the matroid inequality we will have $w(\bar{x}_i)+w(x_i)\geq w(x_{i-1})+w(x_{i+1})$. Hence if we replace the path $x_{i-1}, \ x_i,\ x_{i+1}$ by $x_{i-1}, \ \bar{x}_i,\ x_{i+1}$, the new one will satisfy $w(\bar{x}_{i})-w(x_{i-1})=1$ and $w(x_{i+1})-w(\bar{x}_i)=-1$. That is, in the sequence of $w$-jumps, the sequence $-1$, 1 is replaced by 1, $-1$. In particular there exists a sequence $\setof*{x_i}_i$ such that we have all the 1's at the first $k={\rm codim}(W\subset V)$ steps. In particular, $w(x_k)=w(\ell)+k$, hence $w(C)\geq w(\ell)+{ \rm codim}(W\subset V)$.

    (b) Assume that $\phi(C)=\combcell$, $C=(\ell, I)$. Then $\filt(\ell)=V$, hence $\ell\leq \maxpoint_V$. Hence $w(\ell)\geq w(\maxpoint_V)$ and $w(C)\geq \ti w(\combcell)$. On the other hand, by the proof of part \It{(b)} of \thmref{lattice-to-sscw}, there exists some $C=(\maxpoint_V,I)$ with $\phi(C)=\combcell$, which realizes the equality $w(C)=\ti w(\combcell)$.
\end{pf}

\subsection{\texorpdfstring{The weighted homotopy equivalence of $(\cw,w)$ and $(\ti \cw, \ti w)$ for a curve singularity}{The weighted homotopy equivalence of (X,w) and (X,w)}
}\label{ss:wehomcurves}\hfill\smallbreak

For an isolated curve singularity $(C,o)$ we considered two weighted CW complexes. The first one is $(\cw,w)=(R(0,c),w)$ ($c\geq c_{\semigp}$), the other one in $(\ti\cw,\ti w)$ associated with $\Im\filt$. In this section we prove that the induced lattice cohomologies and graded roots are canonically isomorphic. To this end, we will find a mapping that induces a homotopy equivalence on all level sets $\setof*{S_n}_n$ and $\setof*{\ti S_n}_n$ in the two complexes.

Note that any point of $\cw$ belongs to a unique (relative) open cell $(\ell, I)^\circ$.

\begin{defn}[label={CW-hom-equiv-def}]
    We define $\alpha\colon \cw\to\ti\cw$ as follows. Fix a cell $\topcell=C=(\ell, I)\in \cellset_*$, and set $\combcell=\combcell(V, \cellgens):=\phi(C)$. Regard $\topcell$ embedded in $(\bbr_{\geq 0})^r$ with its natural coordinates. Then the restriction of $\alpha$ to $C^\circ=\topcell^\circ$ is the composition $\kappa_{\combcell}\circ \beta_{C^\circ}$, where $\beta_{C^\circ}:\topcell^\circ\to \topcell^\circ (\combcell(V,\cellgens))$ is given by
    \begin{equation}\label{eq:alpha}
        \beta_{C^\circ}(p)=\Bigl(V,\bigl(\cellgens\ni\cellgen\mapsto\max\setdef{p_i-\ell_i}{i\in I, \filt(\ell+E_i)=\cellgen}\bigr)\Bigr)\in
        \topcell^\circ(\combcell(V,\cellgens)).
    \end{equation}
\end{defn}

\begin{thm}[label={CW-hom-equiv}]
    For all level sets $\levelset_n\subseteq\cw$ and $\ti\levelset_n\subseteq\ti\cw$ $(n\in \bbz$), we have
    \begin{equation}\label{eq:Imalpha}
        \Im\bigl(\alpha|_{\levelset_n}\bigr)\subseteq\ti\levelset_n,
    \end{equation}
    and $\alpha|_{\levelset_n}\colon\levelset_n\to\ti\levelset_n$ is a homotopy equivalence.
\end{thm}

\begin{pf} The proof consists of several steps. We start with the following.
    \begin{lem}[label={alpha-is-cont}]
        The map $\alpha$ is continuous.
    \end{lem}
    \begin{pf}
        It is sufficient to prove that $\alpha|_\cell$ is continuous for any closed cell $\cell=(\ell,I)$ in $\cw$. We prove this by induction on $\dim\cell=\abs{I}$. For $\dim\cell=0$, there is nothing to prove. Otherwise, assume that we know this continuity whenever $\dim\cell<q$ for some integer $q>0$. Then fix some $C=(\ell, I)$ of dimension $q$. At a point $p\in(\ell, I)^\circ$, a small neighborhood $U$ is mapped into the open cell $\kappa_{\combcell}(\topcell^\circ (\combcell(V,\cellgens)))$, so $\alpha(x)=\kappa_{\combcell}(V,u)$ for all $x\in U$, with $u\in(0,1)^\cellgens$ defined as $u(\cellgen)$ being the maximum in the \thmref[text={definition}]{CW-hom-equiv-def}. This is clearly continuous in $x$. If $p\in\D\cell$, then by the induction hypothesis, $\alpha|_\cellface$ is continuous at $p$ for all proper faces $\cellface\subseteq\cell$ containing $p$, hence $\alpha|_{\D\cell}$ is continuous at $p$.

        We then only need to verify the following fact. Let $\bar{\alpha}:C\to \ti\cw$ be the natural continuous extension of $\alpha|_{C^\circ}:C^\circ \to \ti\cw$ to $C$, that is, $\bar{\alpha}=\kappa_{\combcell}\circ \beta_{C}$, where $\beta_{C}:C\to \topcell (\combcell(V,\cellgens))$ is given by the very same formula as in \eqref{eq:alpha} as a continuous extension of $\beta_{C^\circ}$. Since both $\beta_{C}$ and $\kappa_{\combcell}$ are continuous, $\bar{\alpha}$ is continuous. Then we have to show that $\bar{\alpha}|_{\partial C}=\alpha|_{\partial C}$.
    
        Fix $p\in \D C$. Using the notations
        \[
            I_0=\setdef{i\in I}{p_i-\ell_i=0}, \quad
            I_1=\setdef{i\in I}{p_i-\ell_i=1}, \quad
            I_*=\setdef{i\in I}{0<p_i-\ell_i<1},
        \]
        we have $p\in D^\circ $ for $D=(\ell+E_{I_1},I_*)$. Let us also introduce
        \[
            \cellgens_1=\setdef{\filt(\ell+E_i)}{i\in I_1}, \quad
            \ti\cellgens_*=\setdef{\filt(\ell+E_{I_1}+E_i)}{i\in I_*}\setminus\setof{\filt(\ell+E_{I_1})}.
        \]
        In a small neighborhood $U$ of $p$, we again have $\alpha(x)=\kappa_{\combcell}(V,u_x)$ for all $x\in U\cap C^\circ $ where $u_x$ is defined as the maximum from above. Thus,
        \[
            u_p=\lim_{\substack{x\to p \\ x\in\ C^\circ}}u_x=\bigl(\cellgens\ni\cellgen\mapsto\max\setdef{p_i-\ell_i}{i\in I,\filt(\ell+E_i)=\cellgen}\bigr),
        \]
        corresponding to the point $\beta_{C}(p)=(V,u_p)\in\topcell(\combcell(V,\cellgens))$. Hence $\bar{\alpha}(p)=\kappa_{\combcell}(V,u_p)$. Meanwhile, by its very definition, $\alpha(p)$ is obtained by evaluating $\kappa_{\combcell_{\ti D}}$ at $\beta_{D^\circ}(p)$, namely at
        \[
            \Bigl(V\cap(\cap_{W\in \cellgens_1}W),\bigl(\ti\cellgens_*\ni\cellgen\to\max\setdef{p_i-\ell_i}{\filt(\ell+E_{I_1}+E_i)=\cellgen}\bigr)\Bigr).
        \]
        But $\beta_{D^\circ}(p)=\delta_{\ti D}(\beta_{C}|_{D_0}(p))$ (for the definition of $\delta_{\ti D}$, $\combcell_{\ti D}$ and for the identity $\kappa_{\ti C_{\ti D}}\circ \delta_{\ti D}=\kappa_{\ti C}|_{\ti D}$ see the proof of Proposition \thmref{gluing-cont}), hence $\bar{\alpha}(p)=\alpha(p)$ follows from the next commutative diagram
    
        \begin{picture}(300,200)(0,0)
            \put(200,180){\makebox(0,0){$D^0$}}
            \put(100,100){\makebox(0,0){$\ti D\ (\mbox{as face of $\topcell(\combcell)$})$}}
            \put(290,100){\makebox(0,0){$\topcell^\circ(\combcell_{\ti D})$}}
            \put(200,20){\makebox(0,0){$\ti\cw$}}
            \put(185,165){\vector(-1,-1){50}}
            \put(215,165){\vector(1,-1){50}}
            \put(135,85){\vector(1,-1){50}}
            \put(265,85){\vector(-1,-1){50}}
            \put(150,150){\makebox(0,0){$\beta_{C}|_{D^\circ}$}}
            \put(150,50){\makebox(0,0){$\kappa_{\combcell}|_{\ti D}$}}
            \put(250,150){\makebox(0,0){$\beta_{D^\circ}$}}
            \put(250,50){\makebox(0,0){$\kappa_{\combcell_{\ti D}}$}}
            \qbezier(180,180)(-100,100)(180,20)
            \put(175,22){\vector(4,-1){5}}
            \qbezier(220,180)(500,100)(220,20)
            \put(225,22){\vector(-4,-1){5}}
            \put(335,160){\makebox(0,0){$\alpha|_{D^\circ}=\kappa_{\combcell_{\ti D}}\circ \beta_{D^\circ}$}}
            \put(50,160){\makebox(0,0){$\bar{\alpha}|_{D^\circ}=\kappa_{\combcell}|_{\ti D}\circ \beta_{C}|_{D^\circ}$}}
            \put(160,100){\vector(1,0){80}}
            \put(200,110){\makebox(0,0){$\delta_{\ti D}$}}
        \end{picture}
    \end{pf}
    
    Next, we verify the inclusion \eqref{eq:Imalpha}.
    
    Fix $C=(\ell, I)$. If $\phi(C)=\combcell$, then by the continuity of $\alpha$ we also have $\alpha((\ell, I))\subset \topcell(\combcell)$. But, $\ti w(\combcell)\leq w(C)$ by Proposition \thmref{Im-F-cell-weight}. Hence, if $C\subset S_n$ then $\topcell(\combcell)\subset\ti S_n$ and $\alpha(C)\subset\topcell(\combcell)\subset \ti S_n$.
    
    \vspace{2mm}
    
    Next, we show that $\alpha|_{\levelset_n}:S_n\to \ti S_n$ is a homotopy equivalence. We do this by checking that $\alpha|_{\levelset_n}$ is onto $\ti\levelset_n$, and the fiber $\alpha^{-1}(\ti p)\cap\levelset_n$ is contractible for all $\ti p\in\ti S_n$ and $n\in\bbz$. Since any $\ti p\in \ti S_n$ has a small neighbourhood $U$ in $\ti S_n$ such that the inclusion $(\alpha|_{S_n})^{-1}(\ti p)\hookrightarrow(\alpha|_{S_n})^{-1}(U)$ is a homotopy equivalence, the map $\alpha|_{S_n}$ is a quasifibration in the sense of Dold and Thom \cite{DoldThom}, hence the contractibility of the fibers imply the homotopy equivalence.
    
    First note that if $\phi(C)=\combcell$ then $\alpha$ maps the cube $C$ surjectively onto $\topcell(\combcell)$. Therefore, if $\ti p\in \ti S_n$ belongs to $\topcell^\circ (\combcell)$, then from Corollary \thmref{cor-felosrolas} and Proposition \thmref{Im-F-cell-weight} follows that there exists $C$ with $\phi(C)=\combcell$ and $w(C)=\ti w(\combcell)$. Hence $\topcell(\combcell)$ (and $\ti p$ too) belongs to the image of $\alpha|_{S_n}$.
    
    Next we fix $\ti p\in \ti S_n\cap \topcell^\circ (\combcell)$ and we wish to describe the fiber $\alpha^{-1}(\ti p)\cap S_n$. Set $V:=R(\combcell)$.
    
    If $p\in\alpha^{-1}(\ti p)\cap S_n $ and $p\in (\ell, I)^\circ$ then $\phi((\ell,I))=\combcell$, $\alpha((\ell,I))= \topcell(\combcell)$, and $w((\ell, I))\leq n$. Then $\ell\in S_n(V)$ (cf. \ref{ss:prel5}) and, in fact, (using the notation $T^*(\ell)=|\maxpoint_V-\ell|$ from the discussions in front of Lemma \thmref{lemma-ig}) $(\ell, I^*(\ell))\subset S_n(V)$. Let $S_n(V,\combcell)$ be the subspace of $S_n(V)$ consisting of cubes of type $(\ell, J)$ such that $J\subset I^*(\ell)$ and there exists a certain cube $(\ell, I)\in\phi^{-1}(\combcell)$ with $J\subset I$ and $w((\ell, I))\leq n$. From Proposition \thmref{Im-F-cell-weight} follows that if $\ell\in S_n(V,\combcell)$ then the cube $(\ell, I^*(\ell))\subset S_n(V,\combcell)$, hence the whole real rectangle $R(\ell, \maxpoint_V)$ belongs to $S_n(V, \combcell)$. In particular, the space $S_n(V,\combcell)$ is contractible.
    
    Then consider the projection $pr:\phi^{-1}(\combcell)\cap S_n\to S_n(V, \combcell)$ induced by the natural coordinate-wise projections $(\ell, I)\to (\ell, J)$, where $J=I\cap I^*(\ell)$. One verifies that these `local cube projections' glue together into a continuous map. In particular, the composition $\alpha^{-1}(\ti p)\cap S_n\hookrightarrow \phi^{-1}(\combcell)\cap S_n\stackrel{pr}{\longrightarrow} S_n(V, \combcell)$, denoted by $pr(\ti p)$, is also continuous. But, by a local study at the level of cubes shows that each fiber of $pr(\ti p)$ is contractible. Since $pr(\ti p)$ is a quasifibration with contractible fibers and contractible base space $S_n(V, \combcell)$, its total space $(\alpha|_{S_n})^{-1}(\ti p)=\alpha^{-1}(\ti p)\cap S_n$ is contractible too. But this fact is what we wish to proof.
    
    This finishes the proof of Theorem \thmref{CW-hom-equiv}.
\end{pf}

\begin{cor}
    For the above weighted CW complexes, we have an isomorphism
    \[
        \lcoh^*(\cw,w)\simeq\lcoh^*(\ol\cw,\ol w), \quad
         {\Got R}(\cw,w)\simeq {\Got R}(\ol\cw,\ol w)
    \]
\end{cor}

\subsection{The deformation map}\label{ss:defmap}\

Let $\setof{(C_t,o)}_{t\in(\bbc,0)}$ be a flat deformation of isolated curve singularities. We assume that $(C_t,o)\subset (\bbc^N,o)$ for any $t$, and we consider the ideal filtration $\widetilde{ \filt}_t $ of $\fnring_N=\fnring_{(\bbc^N,0)} $ as in \ref{ss:prel5}. We can assume that for $t\not=0$ the lattice $L(C_t,o)=\bbz^{r_t}$ of $(C_t,o)$ is constant, and for fixed $\ell\in L(C_t,o) $ and $t\not=0$ the Hilbert function $\hh_{t}(\ell)$ is constant. Hence, the corresponding weighted CW complex $\ti \cw_{t\not=0}$ is constant as well. Let $\ti \cw_{t=0}$ be the weighted complex of $(C_0,o)$ too.

 First, we wish to define a map $\sk_0\ti\cw_{t\not=0}\to \sk_0\ti\cw_{t=0}$. Recall that if $\ti\cw$ is the complex associated with a curve $(C,o)$ then its 0-skeleton $\sk_0\ti\cw$ is in bijection with the set of lattice points of type $\maxpoint _V\in L(C,o)$, and these correspond precisely to the semigroup elements of $(C,o)$, cf. Lemma \thmref{ss-set-props}.

 Let us fix such a lattice point $\maxpoint=\maxpoint_V\in L(C_{t\not=0})$ (or, differently said, $\maxpoint \in \sk_0\ti\cw_{t\not=0}$).

For each $t\not=0$ consider the $\bbc$-subvector space $\widetilde{\filt}_{t\not=0}(\maxpoint)\subset\fnring_N $ of finite and constant codimension. This family considered in a convenient Grassmannian has a limit $\widetilde{\filt}_{lim}(\maxpoint)$, where $\widetilde{\filt}_{lim}(\maxpoint)\subset \fnring_N$ is a subvector space of the same codimension $\hh_{t\not=0}(\maxpoint)$. In fact, the above infinite Grassmannian can be replaced by a finite Grassmannian as follows. Let $J$ be an ideal of $\fnring_N$ (e.g. a power $\mathfrak{m}_N^d$ of the maximal ideal with sufficiently large $d$) such that $J\subset \widetilde{\filt}_{t\not=0}(\maxpoint)$ for any $t\not=0$. Then $\widetilde{\filt}_{t\not=0}(\maxpoint)/J\subset \fnring_N/J$ has a limit in the compact Grassmannian, which can be lifted as $\widetilde{\filt}_{lim}(\maxpoint)$ into $ \fnring_N$ via $\fnring_N\to \fnring_N/J$.

Even more, by the Curve Selection Lemma (see e.g.\ \cite{MBook}), applied in this finite model, for a generic element $g\in \widetilde{\filt}_{lim}(\maxpoint)$ (regarded as a germ of $\fnring_N$), there exists a family of germs $\setof*{g_t}_{t\in (\bbc,0)}$ such that $g_0=g$ and $g_t\in\widetilde{\filt}_{t}(\maxpoint)$ for $t\not=0$, and $g_t$ is generic in $\widetilde{\filt}_{t}(\maxpoint)$.

For $t\neq 0$, since $\maxpoint\in\semigp_{(C_t,o)}$ and $g_t$ is generic in $\widetilde{\filt}_{t}(\maxpoint)$, for its class $[g_t]\in\filt_t(\maxpoint)\subset\fnring_{(C_t,o)}$ we get
\begin{equation}\label{eq:vege1}
    \Got v_{(C_t,o)}[g_t]=\maxpoint \in\sk_0\ti\cw_{t\not=0}=\semigp_{(C_t,o)}.
\end{equation}
Let us define $\maxpoint_{lim}$ as $\Got v_{(C_0,o)}[g]\in \sk_0\ti\cw_{t=0}=\semigp_{(C_0,o)}$.

The main result of this subsection is the following.

\begin{thm}[label={th:DEF2}]
    Let $\setof{(C_t,o)}_{t\in(\bbc,0)}$ be a flat deformation of isolated curve singularities. Then the map $\sk_0\ti\cw_{t\not=0}\to \sk_0\ti\cw_{t=0}$, $\ell^*\mapsto \ell^*_{lim}$, induces a degree zero graded $\bbz[U]$-module morphism $\lcoh^0(C_{t=0},o)\to \lcoh^0(C_{t\not=0},o)$, and similarly a graded (graph) map of degree zero at the level of graded roots ${\Got R}(C_{t\not=0},o)\to{\Got R}(C_{t=0},o)$.
\end{thm}
\begin{pf}
    Note that $|\maxpoint|$ (respectively $|\maxpoint_{lim}|$) is the intersection multiplicity in $(\bbc^N,o)$ of $(C_t,o)$ with $\setof{g_t=0}$, $t\neq 0$, (respectively of $(C_0,o)$ with $\setof{g=0}$). Therefore, by semicontinuity,
    \begin{equation}\label{eq:vege2}
           \abs{\maxpoint}=(C_t, g_t)_{o\in\bbc^N}\leq (C_0, g_0)_{o\in\bbc^N}=\abs{\maxpoint_{lim}}.
    \end{equation}
    \begin{lem}[label={lem_F}]
        \begin{enumerate}[label={(\alph*)}]
            \item $\widetilde{\filt}_{lim}(\maxpoint)\subset \widetilde{\filt}_{t=0}(\maxpoint_{lim})$.
        
            \item $\widetilde{\filt}_{t=0}(\maxpoint_{lim})$ is the smallest subspace of type $\widetilde{\filt}_{t=0}(\ti \ell)\subset \fnring $ ($ \ti \ell\in L(C_{t=0},o)$) with $\widetilde{\filt}_{lim}(\maxpoint)\subset \widetilde{\filt}_{t=0}(\ti \ell)$.
        \end{enumerate}
    \end{lem}
    \begin{pf}
        (a) For the generic element $g_0\in \widetilde{\filt}_{lim}(\maxpoint)$ we have $\Got v_{(C_0,o)}[g_0]=\maxpoint_{lim}$, hence by semicontinuity, for any $f\in \widetilde{\filt}_{lim}(\maxpoint)$ we have $\Got v_{(C_0,o)}[f]\geq \maxpoint_{lim}$.
    
        (b) $\widetilde{\filt}_{lim}(\maxpoint)\subset \widetilde{\filt}_{t=0}(\ti \ell)$ implies $\Got v_{(C_0,o)}[g_0]=\maxpoint_{lim}\geq \ti \ell$. Hence $\widetilde{\filt}_{t=0}(\maxpoint_{lim})\subset \widetilde{\filt}_{t=0}(\ti \ell)$.
    \end{pf}
    In particular, we also have for the codimensions in $\fnring_N$
    \begin{equation}\label{eq:vege3}
           \hh_{t\not=0}(\maxpoint)={\rm codim}\widetilde{\filt}_{t\not=0}(\maxpoint)=
           {\rm codim}\widetilde{\filt}_{lim}(\maxpoint)\geq
           {\rm codim}\widetilde{\filt}_{t=0}(\maxpoint_{lim})= \hh_{t=0}(\maxpoint_{lim}).
    \end{equation}
    Let us introduce the notations $d:=|\maxpoint_{lim}|-|\maxpoint|\in \bbz_{\geq 0}$ (cf. \eqref{eq:vege2}) and
    \[
        e:= \hh_{t\not=0}(\maxpoint)- \hh_{t=0}(\maxpoint_{lim})={\rm codim } (\widetilde{\filt}_{lim}(\maxpoint)\subset \widetilde{\filt}_{t=0}(\maxpoint_{lim}))\in \bbz_{\geq 0}.
    \]
    Then the identities \eqref{eq:vege2} and \eqref{eq:vege3} combined give
    \begin{equation}\label{eq:vege4}
        w_{t\not=0}(\maxpoint)= w_{t=0}(\maxpoint_{lim})+2e+d\geq w_{t=0}(\maxpoint_{lim}).
    \end{equation}
    Therefore, for any $n\in \bbz$, $\maxpoint\mapsto \maxpoint_{lim}$ maps $\sk_0\ti\cw_{t\not=0}\cap \ti S_n$ to $\sk_0\ti\cw_{t=0}\cap \ti S_n$.
    
    Let us consider next an 1-cell $(\maxpoint, \setof{i})\in\cw_{t\neq 0}$, where $E_i$ is a base vector of $L(C_{t\neq 0})$. It induces the 1-cell $\combcell$ with endpoints $\maxpoint$ and $(\maxpoint+E_i)^*$, where $(\maxpoint+E_i)^*$ is associated with $\maxpoint+E_i$ as in Lemma \thmref{uniquemax_0} (it is the smallest semigroup element larger than $\maxpoint+E_i$). By the above map $\maxpoint $ is sent to $\maxpoint_{lim}$ and $(\maxpoint+E_i)^* $ to $(\maxpoint+E_i)_{lim}$. Moreover, by the definition from subsection \ref{ss:linsubspacesimf} we have
    \begin{equation}\label{eq:vege5}
        w_{t\not=0}(\combcell)= w_{t\not =0}(\maxpoint)+1.
    \end{equation}
    Note that all the 1-cells of $\ti\cw_{t\not=0}$ are of this type
    
    For any $n\in\bbz$ we will prove the following fact. If $\maxpoint$ and $(\maxpoint+E_i)^*$ and $\kappa (\topcell(\combcell))$ belong to $\ti\cw_{t\not=0}\cap \ti S_n$, then $\maxpoint_{lim}$ and $(\maxpoint+E_i)_{lim}$ sit in the same connected component of $\ti\cw_{t=0}\cap \ti S_n$.
    
    Indeed, from the assumption we have
    \begin{equation}\label{eq:vege6}
        w_{t\not=0}(\maxpoint )+1\leq n.
    \end{equation}
    From Lemma \thmref{lem_F}\It{(b)} we obtain that $\widetilde{\filt}_{t=0}((\maxpoint+E_i)_{lim})\subset \widetilde{\filt}_{t=0}(\maxpoint_{lim})$, hence $\maxpoint_{lim}\leq (\maxpoint+E_i)_{lim}$ too. In particular, we can consider the inclusions
    \[
        \widetilde{\filt}_{lim}((\maxpoint+E_i)^*)\subset
        \widetilde{\filt}_{t=0}((\maxpoint+E_i)_{lim})\subset
        \widetilde{\filt}_{t=0}(\maxpoint_{lim}).
    \]
    Therefore,
    \[
        \codim(\widetilde{\filt}_{t=0}((\maxpoint+E_i)_{lim})\subset \widetilde{\filt}_{t=0}(\maxpoint_{lim}))\leq
        \codim(\widetilde{\filt}_{lim }(((\maxpoint+E_i)^*)\subset \widetilde{\filt}_{t=0}(\maxpoint_{lim}))
    \]
     where the right hand side equals
     \[
        \codim(\widetilde{\filt}_{lim }(((\maxpoint+E_i)^*)\subset\widetilde{\filt}_{lim }(\maxpoint))+
        \codim(\widetilde{\filt}_{lim }(\maxpoint))\subset \widetilde{\filt}_{t=0}(\maxpoint_{lim}))
        =1+e.
     \]
     Thus, we have two points $\maxpoint_{lim}$ and $(\maxpoint+E_i)_{lim}$ in $\ti S_n$ such that $w_{t=0}(\maxpoint_{lim})+2e+d\leq n-1$ and ${\rm codim}(\widetilde{\filt}_{t=0}((\maxpoint+E_i)_{lim})\subset \widetilde{\filt}_{t=0}(\maxpoint_{lim}))\leq 1+e$. Then for any increasing path in $\cw_{t=0}$ which connects $\maxpoint_{lim}$ and $(\maxpoint+E_i)_{lim}$ we have the following facts: the length of the path is $\leq 1+e$, and at any consecutive step along the path $w_{t=0}$ can increase maximum by 1. Hence along any such path $w_{t=0}$ is $\leq w_{t=0}(\maxpoint_{lim})+e+1$. Since
     \[
        w_{t=0}(\maxpoint_{lim})+e+1\leq n-1-2e-d+e+1=n-e-d\leq n
     \]
     we obtain that $\maxpoint_{lim}$ and $(\maxpoint+E_i)_{lim}$ belong to the same connected component of $\ti\cw_{t=0}\cap \ti S_n$.
    
     Summarized, $\maxpoint\mapsto \maxpoint_{lim}$ induces a well-defined map $\pi_0(\ti\cw_{t\not=0}\cap \ti S_n)\to \pi_0(\ti\cw_{t=0}\cap \ti S_n)$ for any $n$, and they are compatible with respect to the inclusions $\ti S_n\hookrightarrow\ti S_{n+1}$.
\end{pf}

\begin{ex}[label={ex:an}]
    Usually it can happen that $\widetilde{\filt}_{lim}(\maxpoint)\not= \widetilde{\filt}_{t=0}(\maxpoint_{lim})$. It is also possible that $\maxpoint$ and $(\maxpoint+E_i)^*$ are not the endpoints of an 1-cell of $\ti\cw_{t=0}$.
    
    Indeed, consider the deformation of isolated irreducible plane curves $\setof{x^2+ty^3+y^5=0}_{t\in(\bbc,0)}$. For $t\neq 0$ and $\ell=0,1,\ldots,6$ we have the following semigroup elements and ($t$-independent) ideals $\widetilde{\filt}_{t\not=0}(\ell)$, which equal the limits $\widetilde{\filt}_{lim}(\ell)$ as well:
    
    \vspace*{4mm}
    
    \begin{center}
        \begin{tikzpicture}[scale=2]
            \draw[fill] (0,0) circle [radius=0.04];
            \draw (1,0) circle [radius=0.04];
            \draw[fill] (2,0) circle [radius=0.04];
            \draw[fill] (3,0) circle [radius=0.04];
            \draw[fill] (4,0) circle [radius=0.04];
            \draw[fill] (5,0) circle [radius=0.04];
            \draw[fill] (6,0) circle [radius=0.04];

            \draw (0,0) -- (6,0);

            \node at (0,-0.4) {$\fnring_N$};
            \node at (1,-0.4) {$\mathfrak{m}_N$};
            \node at (2,-0.4) {$\mathfrak{m}_N$};
            \node at (3,-0.4) {$(x,y^2)$};
            \node at (4,-0.4) {$\mathfrak{m}_N^2$};
            \node at (5,-0.4) {$(x^2,xy,y^3)$};
            \node at (6,-0.4) {$(x^2,xy^2,y^3)$};
        \end{tikzpicture}
    \end{center}
    \nni
    For $t=0$ the semigroup elements and ideals $\widetilde{\filt}_{t=0}(\ell)$ are:
    
    \vspace*{4mm}

    \begin{center}
        \begin{tikzpicture}[scale=2]
            \draw[fill] (0,0) circle [radius=0.04];
            \draw (1,0) circle [radius=0.04];
            \draw[fill] (2,0) circle [radius=0.04];
            \draw (3,0) circle [radius=0.04];
            \draw[fill] (4,0) circle [radius=0.04];
            \draw[fill] (5,0) circle [radius=0.04];
            \draw[fill] (6,0) circle [radius=0.04];

            \draw (0,0) -- (6,0);

            \node at (0,-0.4) {$\fnring_N$};
            \node at (1,-0.4) {$\mathfrak{m}_N$};
            \node at (2,-0.4) {$\mathfrak{m}_N$};
            \node at (3,-0.4) {$(x,y^2)$};
            \node at (4,-0.4) {$(x,y^2)$};
            \node at (5,-0.4) {$(x,y^3)$};
            \node at (6,-0.4) {$(x^2,xy,y^3)$};
        \end{tikzpicture}
    \end{center}
    \nni
    Hence, $4_{lim}=4$ but $\widetilde{\filt}_{lim}(4)=\mathfrak{m}^2\varsubsetneq (x,y^2)= \widetilde{\filt}_{t=0}(4_{lim})$. Moreover, $5_{lim}=6$, hence $4_{lim}$ and $5_{lim}$ are not endpoints of an 1-cell of $\ti\cw_{t=0}$.
\end{ex}

\begin{rem}[label={rem:ketfele}]
    If $(C_{t=0},o)$ is irreducible then we have constructed two maps $\pi_0(\ti\cw_{t\not=0}\cap \ti S_n)\to\pi_0(\ti\cw_{t=0}\cap \ti S_n)$. The first one, considered for plane curves, is induced by $\maxpoint\mapsto |\maxpoint|$, and in Theorem \thmref{th:DEF} we probed that it induces a graded graph map of degree zero at the level of graded roots ${\Got R}(C_{t\not=0},o)\to{\Got R}(C_{t=0},o)$. In fact, in Remark \thmref{rem:GN} we noticed that the statement can be extended to non-plane curves as well. (This fact will follow from the computations from below.) The second map, constructed in Theorem \thmref{th:DEF2}, is induced by $\maxpoint\mapsto \maxpoint_{lim}$. Recall that by \eqref{eq:vege2} $|\maxpoint|\leq \maxpoint_{lim}$.
    
    Let us consider any lattice point $ \ell\in L(C_{t=0})=\bbz$ with $|\maxpoint|\leq \ell\leq\maxpoint_{lim}$. It satisfies
    \[
        \widetilde{\filt}_{t=0}(|\maxpoint|)\supset
        \widetilde{\filt}_{t=0}( \ell)\supset
        \widetilde{\filt}_{t=0}(\maxpoint_{lim})\supset
        \widetilde{\filt}_{lim}(\maxpoint).
    \]
    Therefore,
    \begin{equation}\label{eq:hdefdef}
        \hh_{t=0}(|\maxpoint|)\leq
        \hh_{t=0}(\ell)\leq
        \hh_{t=0}(\maxpoint_{lim})\leq
        \hh_{t\not=0}(\maxpoint).
    \end{equation}
    Moreover, since $-\ell \leq -|\maxpoint|$ we get $ w_{t=0}(\ell)\leq w_{t\not=0}(\maxpoint)$ for any such $\ell$. Hence, indeed (even if $(C_t,o)$ are not plane curves), both $\ell^*\mapsto |\ell^*|$ and $\ell^*\mapsto \ell^*_{lim}$ induce maps $\sk_0\ti\cw_{t\not=0}\cap \ti S_n\to\sk_0\ti\cw_{t=0}\cap \ti S_n$, and by the original arguments from above, maps $\pi_0(\ti\cw_{t\not=0}\cap \ti S_n)\to\pi_0(\ti\cw_{t=0}\cap \ti S_n)$.
    
    Furthermore, any 0 and 1-dimensional cell contained in the segment $[|\maxpoint|, \maxpoint_{lim}]$ of $\ti\cw_{t=0}$ has weight $\leq w_{t\not=0}(\maxpoint )$. Therefore, if $\maxpoint\in \ti\cw_{t\not=0}\cap\ti S_n$ then the segment $[|\maxpoint|, \maxpoint_{lim}]$ belongs to $\ti\cw_{t=0}\cap\ti S_n$ hence $|\maxpoint|$ and $\maxpoint_{lim}$ belong to the same connected component of $\ti\cw_{t=0}\cap\ti S_n$. In particular, the two maps induces the same graded graph map at the level of graded roots ${\Got R}(C_{t\not=0},o)\to{\Got R}(C_{t=0},o)$ (and similarly at the level of $\lcoh^0$).
\end{rem}
\begin{rem}
    In the case of a Gorenstein curve $(C,o)$ both $\lcoh^*(C,o)$ and ${\Got R}(C,o)$ have an additional $\bbz_2$-symmetry, cf. \ref{ss:GORduality}. One can ask if in the case of a flat deformation of Gorenstein curves the graded graph map at the level of graded roots ${\Got R}(C_{t\not=0},o)\to{\Got R}(C_{t=0},o)$ is $\bbz_2$-equivariant? The answer is no, a fact which can be checked already at the deformation considered in Example \thmref{ex:an}.
    
    However, we formulate the following:
    \begin{con}
        In the case of a $\delta$-constant flat deformation of Gorenstein curves the graded graph map at the level of graded roots ${\Got R}(C_{t\not=0},o)\to{\Got R}(C_{t=0},o)$ (constructed above) is $\bbz_2$-equivariant.
    \end{con}
\end{rem}

As an enhancement of Theorem \thmref{th:DEF2} we also formulate the following general conjecture.
\begin{con}
    Let $\setof*{(C_t,o)}_{t\in(\bbc,0)}$ be a flat deformation of isolated curve singularities. Then the map $\ell^*\mapsto \ell^*_{lim}$ from the level $\sk_0\ti\cw_{t\not=0}\cap \ti S_n\to \sk_0\ti\cw_{t=0}\cap \ti S_n$ can be lifted to $\ti\cw_{t\not=0}\cap \ti S_n\to \ti\cw_{t=0}\cap \ti S_n$ inducing a degree zero graded $\bbz[U]$-module morphism $\lcoh^q(C_{t=0},o)\to \lcoh^q(C_{t\not=0},o)$ for any $q\geq 1$.
\end{con}



\end{document}